\documentclass[sort&compress,3p]{elsarticle}

\usepackage{amssymb}
\usepackage{amsmath}
\usepackage{bm}
\usepackage{amsthm}
\usepackage{caption}
\usepackage{subcaption}
\usepackage{xcolor}
\usepackage{soul}
\usepackage[ruled,vlined]{algorithm2e}
\usepackage{setspace}
\usepackage{multirow}

\usepackage{pgfplots}
\pgfplotsset{compat=1.17}

\def\NN{\mathbb{N}}
\def\RR{\mathbb{R}}

\newcommand{\bfm}[1]{\bm{#1}}
\newcommand{\mc}[1]{\mathcal{#1}}

\newcommand{\defset}[2]{\left\{ #1: \ #2 \right\}}
\newcommand{\sprod}[2]{\left\langle#1, #2\right\rangle}
\newcommand{\sprodS}[2]{\langle#1, #2\rangle}
\newcommand{\norm}[1]{\| #1 \|}
\newcommand{\C}[1]{{C}^{#1}}
\newcommand{\G}[1]{{G}^{#1}}
  
\newcommand{\bV}{\bfm{{V}}}
\newcommand{\bVm}{\bm{V}_m}
\newcommand{\bE}{\boldsymbol{\mc{E}}}

\newcommand{\pd}{d}

\newcommand{\poly}{\mathbb{P}}

\newcommand{\bern}[2]{B^{#1}_{#2}}
\newcommand{\bernT}[2]{B^{\mbox{{\tiny{$\triangle$}}}, #1}_{#2}}
\newcommand{\bernB}[2]{B^{\mbox{{\tiny{$\Box$}}}, #1}_{#2}}

\newcommand{\indSetT}[1]{\mathbb{I}^{\mbox{{\tiny{$\triangle$}}}}_{#1}}
\newcommand{\indSetQ}[1]{\mathbb{I}^{\mbox{{\tiny{$\Box$}}}}_{#1}}

\newcommand{\w}{\omega}

\newcommand{\tr}{\theta}

\newcommand{\set}[1]{\left\{ #1 \right\}}

\newcommand{\parDomT}{\triangle_0}
\newcommand{\parDomQ}{\Box_0}
\newcommand{\parDom}{\mathrm{D}}

\newcommand{\indOmega}{\mathcal{I}_{\Omega}}
\newcommand{\indE}{\mathcal{I}_{\mc E}}
\newcommand{\indV}{\mathcal{I}_{\bfm{V}}}
\newcommand{\indQ}{\mathcal{I}_{\mbox{{\tiny{$\Box$}}}}}
\newcommand{\indT}{\mathcal{I}_{\mbox{{\tiny{$\triangle$}}}}}

\newcommand{\lfunvert}[1]{\lambda^{v}_{#1}}
\newcommand{\lfunedgetr}[1]{\lambda^{e,0}_{#1}}
\newcommand{\lfunedgew}[1]{\lambda^{e, 1}_{#1}}
\newcommand{\lfunTri}[1]{\lambda^{{t}}_{#1}}
\newcommand{\lfunQ}[1]{\lambda^{{q}}_{#1}}

\newcommand{\Bv}[1]{\varphi_{#1}^{v}}
\newcommand{\Bt}[1]{\varphi_{#1}^{t}}
\newcommand{\Bq}[1]{\varphi_{#1}^{q}}

\newcommand{\Betr}[1]{\varphi_{#1}^{e,\tr}}
\newcommand{\Bew}[1]{\varphi_{#1}^{e,\w}}

\newcommand{\xia}{\xi_1}
\newcommand{\xib}{\xi_2}
\newcommand{\xic}{\xi_3}
\newcommand{\xiBo}{\bfm{\xi}}
\newcommand{\xiu}{\xi}

\newcommand{\xa}{x_1}
\newcommand{\xb}{x_2}
\newcommand{\xc}{x_3}
\newcommand{\xBo}{\bfm{x}}
\newcommand{\parV}{\zeta}

\newcommand{\Da}{\partial_1}
\newcommand{\Db}{\partial_2}
\newcommand{\Dt}[1]{\frac{d}{d #1} }
\newcommand{\Dti}[2]{\frac{d^{#2}}{d #1^{#2}} }
\newcommand{\Dv}{D}
\newcommand{\Dx}{\partial_{x_1}}
\newcommand{\Dy}{\partial_{x_2}}
\newcommand{\Dz}{\partial_{x_3}}

\newcommand{\uQ}{\underline{\makebox[0.58em][c]{$Q$}}}
\newcommand{\Qkl}[2]{\uQ^{\,(#1)}_{\,#2}}
\newcommand{\tkl}[2]{\underline{t}^{(#1)}_{\,#2}}
\newcommand{\nkl}[2]{\underline{n}^{(#1)}_{\,#2}}
\newcommand{\dkl}[2]{\underline{d}^{(#1)}_{\,#2}}

\newtheorem{theorem}{Theorem}

\newtheorem{lemma}[theorem]{Lemma}

\newtheorem{remark}{Remark}
\newtheorem{example}{Example}
\newtheorem{definition}{Definition}

\begin{document}

\begin{frontmatter}

\title{$\C{1}$-smooth isogeometric spline functions of general degree over mixed triangular and quadrilateral surface meshes} 

\author[UL,UR]{Jan Gro\v{s}elj}
\ead{jan.groselj@fmf.uni-lj.si}
\author[CUAS]{Mario Kapl}
\ead{m.kapl@cuas.at}
\author[UL,UR]{Marjeta Knez} 
\ead{marjetka.knez@fmf.uni-lj.si}
\author[RICAM]{\corref{cor}Thomas Takacs}
\ead{thomas.takacs@ricam.oeaw.ac.at}
\author[UP]{Vito Vitrih}
\ead{vito.vitrih@upr.si}
\address[UL]{FMF, University of Ljubljana, Jadranska 19, 1000 Ljubljana, Slovenia}
\address[UR]{IMFM, Jadranska 19, 1000 Ljubljana, Slovenia} 
\address[CUAS]{Faculty of Engineering $\&$ IT, Carinthia University of Applied Sciences, Europastra\ss{}e 4, 9524 Villach, Austria}
\address[RICAM]{Johann Radon Institute for Computational and Applied Mathematics, Austrian Academy of Sciences, Altenberger Str. 69, 4040 Linz, Austria}
\address[UP]{UP FAMNIT and UP IAM, University of Primorska, Glagolja\v ska 8, 6000 Koper, Slovenia}

\cortext[cor]{Corresponding author.}

\begin{abstract}
We describe the construction of a $C^1$ super-smooth isogeometric spline space over mixed triangular and quadrilateral surface meshes with mesh elements of an arbitrary degree~$\delta \geq 1$. The constructed $C^1$ super-smooth isogeometric spline space consists of all $C^1$-smooth isogeometric spline functions, which are additionally $C^2$-smooth at all vertices, and whose pull-backs are polynomials of an arbitrary degree~$d \geq \delta$ on the single mesh elements. We develop a general framework to study the $C^1$-smoothness conditions of a function across an edge of two neighboring mesh elements and to construct a basis of the $C^1$ super-smooth isogeometric spline. The resulting basis functions possess local supports and are categorized into four types of functions given by vertex, edge, triangle and quadrilateral basis functions. Several numerical examples are presented for solving interpolation and least-squares approximation problems over different mixed triangular and quadrilateral surface meshes, which demonstrate the potential of the generated $C^1$ super-smooth isogeometric spline space for applications in isogeometric analysis. 
\end{abstract}


\end{frontmatter}

\section{Introduction}
\label{sec:intro}

Mixed triangular and quadrilateral meshes are a promising alternative to the typically used triangular or quadrilateral meshes to model the geometry of complex surface domains. Unlike pure triangular or quadrilateral meshes, the mixed triangular and quadrilateral meshes are not restricted to one type of elements, namely either to triangles or quadrilaterals, but enable the use of both types of elements. This makes the mixed triangular and quadrilateral meshes a powerful and flexible tool in surface design~\cite{Pe95,BoLePiPuSiTaZo13}. Examples of relevant applications of mixed triangular and quadrilateral meshes are surface tiling~\cite{PePoWo18}, surface design in architecture~\cite{PoEiVaWa15}, the representation of trimmed objects~\cite{To22,KiSeYo09}, the merging of quadrilateral meshes with possibly different mesh sizes using triangles in the transition areas~\cite{ThEnSchTeSc22} as considered e.g. in adaptive refinement~\cite{BaShWe1983,FuSch20}, and the combining of structured quadrilateral meshes with unstructured triangular meshes as used e.g. in computational fluid dynamics~\cite{DuFoGuVi09}.

To solve fourth-order partial differential equations (PDEs), such as the biharmonic problem, e.g.~\cite{BaDeQu15,FaJuKaTa22,NgKaPe15}, the Kirchhoff-Love shell problem, e.g.~\cite{FaVeKiKa23,KiBlLi09,KiBaHs10}, the Cahn-Hilliard equation, e.g.~\cite{GoCaHu09,KaMeBo16}, or problems of strain gradient elasticity, e.g.~\cite{MaReBeJu18,KhakaloNiiranenC1}, over surface meshes via the Galerkin method using isogeometric analysis~\cite{HuCoBa05}, the associated surface has to be $G^1$-smooth~\cite{Pe02}, and the functions for discretizing the solution of the considered PDE need to be at least $C^1$-smooth. The use and construction of $C^1$-smooth functions over triangular and quadrilateral meshes for numerical simulation were of interest long before the rise of isogeometric analysis. The first $C^1$-smooth finite elements were developed in the 60s with the rectangular Bogner--Fox--Schmit element~\cite{BoFoSc65} and the triangular Argyris element~\cite{ArFrSc68} and Bell element~\cite{Be69}. The first example of a $C^1$-smooth element over general quadrilaterals is the Brenner--Sung element~\cite{BrSu05}. While in the context of the finite element method $C^1$-smooth functions over triangular meshes played a more prominent role, cf.~\cite{Ci02,BrSc07,LaSc07}, $C^1$-smooth functions over quadrilateral meshes gained increasing attention and importance since the introduction of isogeometric analysis, cf.~\cite{CoSaTa16,HuSaTaTo21,KaSaTa18,VeWeMaTaTo24}. Driven by geometric modeling, the number of existing constructions of $C^1$-smooth functions over triangular and quadrilateral meshes is large. Only a few recent examples, which are also suited for the use in isogeometric analysis, are \cite{GrMa22,GrSp21,SpMaPeSa12,Sp13,JaQi14,GrSaSp26,ZaQi19} in case of triangular meshes and \cite{ChAnRa18,BeMa14,BlMoXu20,FaJuKaTa22,KaBuBeJu16,KaSaTa19,KaNgPe17,ToSpHu17}
in case of quadrilateral meshes.   

In contrast to purely triangular or quadrilateral meshes, the existing literature on $C^1$-smooth (isogeometric) functions over mixed triangular and quadrilateral meshes is very limited. The work~\cite{To22} presents the construction of quadratic isogeometric spline functions over mixed quadratic triangular and quadrilateral meshes with an application to solving second-order PDEs on watertight reconstruction of trimmed surfaces. The constructed quadratic splines are $C^1$-smooth in all structured parts of the mixed mesh but only $C^0$-smooth in unstructured regions. In \cite{ThEnSchTeSc22}, specific isogeometric spline functions, called U-splines, with $C^1$ or even higher smoothness over unstructured quadrilateral meshes are generated. They also allow to deal with mixed triangular and quadrilateral meshes, e.g. obtained from combining different quadrilateral meshes using triangles in the transition area. In this case, however, the constructed splines are only $C^0$-smooth in the transition area. A first technique, which generates globally $C^1$-smooth isogeometric spline functions over mixed triangular and quadrilateral meshes, is presented in~\cite{MoViVi16}. There, a general framework with a specific study of some concrete configurations is described including the computation of the dimension and of a basis of the resulting $C^1$-smooth spline space. In contrast to most other approaches, in~\cite{MoViVi16} a topological mesh is used rather than a physically given one. The method~\cite{GrKaKnTaVi2020} combines the classical Argyris $C^1$ triangular finite element~\cite{ArFrSc68} with the Brenner-Sung $C^1$ quadrilateral finite element~\cite{BrSu05,KaSaTa20}, which represents the corresponding analogue finite element for quadrilaterals. In this way, $C^1$-smooth isogeometric spline functions over planar mixed triangular and quadrilateral meshes with (bi-)linear mesh elements are constructed, which are additionally $C^2$-smooth at the vertices, and whose pull-backs are polynomials of (bi-)degree~$p \geq 5$ on the single mesh elements. The optimal approximation properties of the generated $C^1$ super-smooth isogeometric spline space are proven theoretically and further verified numerically by solving interpolation and $L^2$-approximation problems as well as the biharmonic problem. In~\cite{GrKaKnTaVi2024}, $C^1$-smooth isogeometric spline functions of general degree over planar mixed triangular and quadrilateral meshes with (bi-)quadratic mesh elements are studied with the focus on the case of two neighboring mesh elements. For this specific configuration, a full theoretical framework was developed to investigate the $C^1$-smoothness conditions of a function across the two adjacent mesh elements. This enabled to compute the dimension as well as a basis of the associated $C^1$-smooth spline space. 

The aim of this paper is to construct for a given mixed triangular and quadrilateral surface mesh with mesh elements of general degree~$\delta \geq 1$ $C^1$-smooth isogeometric spline functions of general degree~$d \geq \delta$, which are additionally $C^2$-smooth at all vertices. This will include the study of the $C^1$-smoothness conditions of an isogeometric function across an edge of two adjacent mesh elements, the construction of a basis of the resulting $C^1$ super-smooth isogeometric spline space --- categorized into vertex, edge, triangle and quadrilateral basis functions --- and the application of the constructed $C^1$-smooth spline functions in isogeometric analysis through the solution of interpolation and least-squares approximation problems. Thereby, the present paper will extend and generalize the works~\cite{GrKaKnTaVi2020, GrKaKnTaVi2024} in several directions. First, the novel method is not limited to planar domains but can also handle surface domains. Second, the technique in this paper is not restricted to (bi-)linear or (bi-)quadratic mesh elements but can also deal with mesh elements of an arbitrary degree~$\delta \geq 1$. Third, while the existing work~\cite{GrKaKnTaVi2024} using (bi-)quadratic mesh elements only focuses on the case of two neighboring mesh elements, the novel approach works for any mixed triangular and quadrilateral surface mesh. Fourth, while in~\cite{GrKaKnTaVi2024} no numerical simulation is performed, the method in this paper demonstrates the potential of the constructed $C^1$-smooth spline functions for the use in isogeometric analysis by solving interpolation and least-squares problems. In addition to these new contributions regarding the construction of $C^1$-smooth isogeometric spline functions on mixed triangular and quadrilateral surface meshes, the presented approach is also novel for meshes consisting purely of triangles or quadrilaterals, where mesh elements of arbitrary degree have not previously been investigated in this context.

The remainder of the paper is organized as follows. Section~\ref{sec:C1-smooth-functions} first introduces the class of mixed triangular and quadrilateral surface meshes used throughout the paper. Afterwards, the general concept of $C^1$-smooth isogeometric spline functions over mixed triangular and quadrilateral surface meshes is explained. In Section~\ref{sec:SplineSpace}, the $C^1$-smoothness conditions of a function over an edge of two neighboring mesh elements are studied in detail. These are then used to construct a Hermite-like basis for an edge, which is employed to generate a basis for the entire $C^1$ super-smooth spline space over a mixed triangular and quadrilateral surface mesh. In this context, super-smooth refers to the fact that $C^2$-smoothness at all vertices is additionally enforced. By construction, all basis functions possess a local support and are categorized into four different types of functions, namely into vertex, edge, triangle and quadrilateral basis functions. In addition, an alternative approach for the construction of a basis, which is based on interpolation, is described. Section~\ref{sec:Examples} employs the constructed $C^1$ super-smooth isogeometric spline space to solve interpolation and least-squares approximation problems over different mixed triangular and quadrilateral surface meshes. Finally, we conclude the paper in Section~\ref{sec:Conclusion}.

\section{$\C{1}$-smooth isogeometric functions over mixed surface meshes} \label{sec:C1-smooth-functions}

In this section we recall the notion of $G^1$ geometrically continuous surfaces which are formed by a partition of triangular and quadrilateral surface elements. We then define $C^1$-smooth isogeometric spline spaces over such surface domains.

\subsection{Triangular and quadrilateral polynomial surface elements} 

Let $\poly^1_d$ denote the space of univariate polynomials of degree $d$ and let $\poly^2_d$ and $\poly^2_{d,d}$ denote the spaces of bivariate polynomials of total degree and bi-degree $d$, respectively. For the parametric domains on triangles and quadrilaterals, we fix 
\begin{equation*}
\parDomT :=\defset{(\xia,\xib)\in [0,1]^2}{\xia + \xib \leq 1}
\quad \mbox{ and } \quad \parDomQ:=\left[0,1\right]^2,     
\end{equation*}
respectively. The Bernstein basis polynomials, denoted by  $\bern{d}{i}$, $\bernT{d}{\mathbf i}$ 
and $\bernB{d}{\mathbf i}$, where $\mathbf i = (i_1,i_2)$, and defined as 
\begin{align*}
 \bern{d}{i}(\xiu) &:= \binom{d}{i}\xiu^i (1-\xiu)^{d-i},  &&i\in \set{0,1,\dots, d},\\    
 \bernT{d}{\mathbf i}(\xia,\xib) &:= \frac{d!}{i_1!i_2!(d-i_1-i_2)!} \xia^{i_1} \xib^{i_2} (1-\xia-\xib)^{d-i_1-i_2},  
&&\mathbf i \in \indSetT{d} := \set{(i_1,i_2)\in \NN_0\times \NN_0, \; i_1+i_2\leq d},\\
 \bernB{d}{\mathbf i}(\xia,\xib) &:= \bern{d}{i_1}(\xia) \bern{d}{i_2}(\xib),  &&\mathbf i \in 
\indSetQ{d} := \set{0,1,\dots, d}^2,
\end{align*}
are used to represent the univariate, triangle and tensor-product polynomials of (bi-)degree $d$, respectively.

Triangular and quadrilateral polynomial surface elements are defined as images of bijective and regular parameterizations 
\[
\bfm{F}: \parDom \to \bfm{F}(\parDom) \subset \mathbb{R}^3,
\]
where $\parDom=\parDomT$ and $\bfm{F} \in (\poly^2_{\delta})^3$ for triangular, and 
$\parDom=\parDomQ$ and $\bfm{F} \in (\poly^2_{\delta,\delta})^3$ for quadrilateral elements, for some $\delta \in\mathbb{Z}^+$. 
Throughout the paper, we denote by $\xiBo=(\xia,\xib)^T$ the ``local'' coordinates with respect to a parameter domain~$\parDom$, and by $\xBo=(\xa,\xb,\xc)^T$ the ``global'' coordinates in the ambient space.
\begin{definition}\label{def-isogeometric-fun}
The function on the surface element
\[
\varphi: \bfm{F}(\parDom) \to \RR, \quad \mbox{such that} \quad 
\varphi \circ \bfm{F} = f \in 
\begin{cases}
\poly^2_{\pd} & \mbox{ if } \;\parDom= \parDomT\\
\poly^2_{\pd, \pd} & \mbox{ if } \; \parDom= \parDomQ
 \end{cases},
\]
where $\pd \geq \delta$, is called the \emph{isogeometric function}. 
\end{definition}
While \emph{isoparametric} refers to the case where the element mapping $\bfm{F}$ and the function $f$ on the element are from the same space, i.e., $d=\delta$, the term \emph{isogeometric} is often used more loosely to describe a spline-based discretization.
\begin{definition}\label{def-unit-normal}
The mapping 
\begin{equation*}
\begin{aligned}
    \bfm{N}: \parDom \to \RR^3, \quad
    \xiBo \mapsto (\Da \bfm{F} \times \Db \bfm{F})(\xiBo),
\end{aligned}
\end{equation*}
provides at every point $\xBo=\bfm{F}(\xiBo)$ on the surface element a normal vector.
\end{definition}
Since $\bfm{F}$ is regular, $\norm{ \bfm{N}(\xiBo)} >0$ and the unit normal $\bfm{N}^\ast := \bfm{N}/\norm{\bfm{N}}$ is thus well-defined. If the element is planar (and with a positive orientation), i.e., $\bfm{F}=(F_1, F_2, 0)^T \simeq (F_1, F_2)^T$ and $\bfm{N}^* (\xiBo)= (0,0,1)^T$, we can describe the graph of the isogeometric function over the planar domain via 
$\xiBo \mapsto (F_1(\xiBo), F_2(\xiBo), f(\xiBo))^T = \bfm{F}(\xiBo) + f(\xiBo) \bfm{N}^*(\xiBo)$.

\subsection{$G^1$ surface domain} 

\begin{figure}[hbt]
\begin{center}
\begin{minipage}{0.32\textwidth}
\centering
\includegraphics[width=0.95\textwidth]{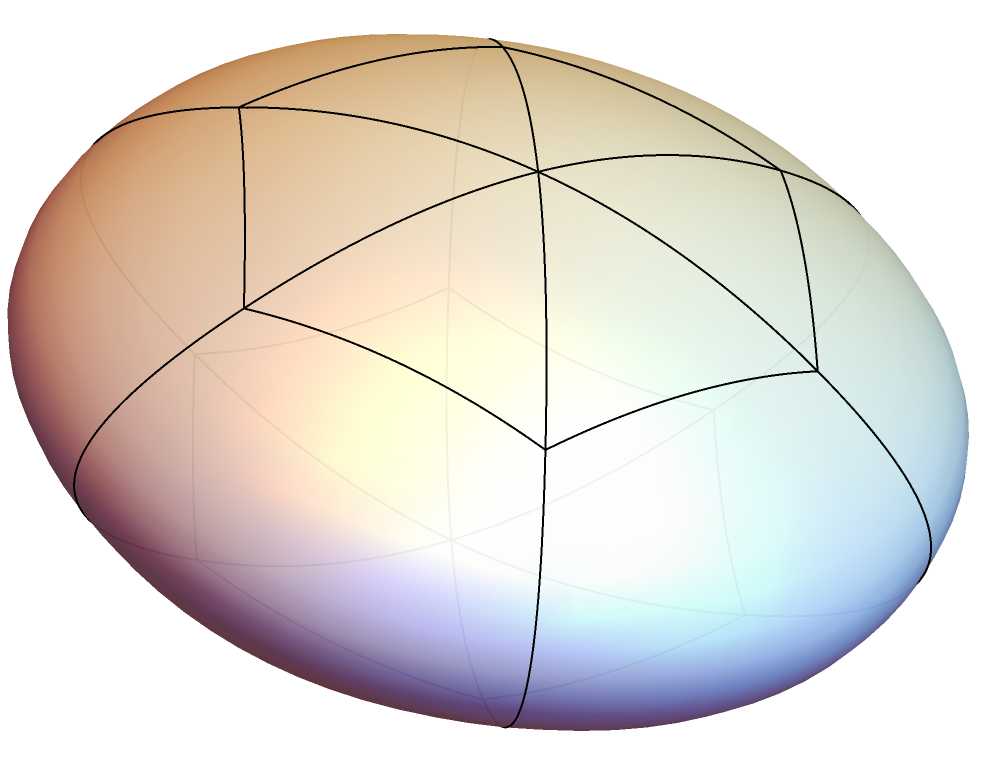}
\end{minipage}
\begin{minipage}{0.32\textwidth}
\centering
\includegraphics[width=1\textwidth]{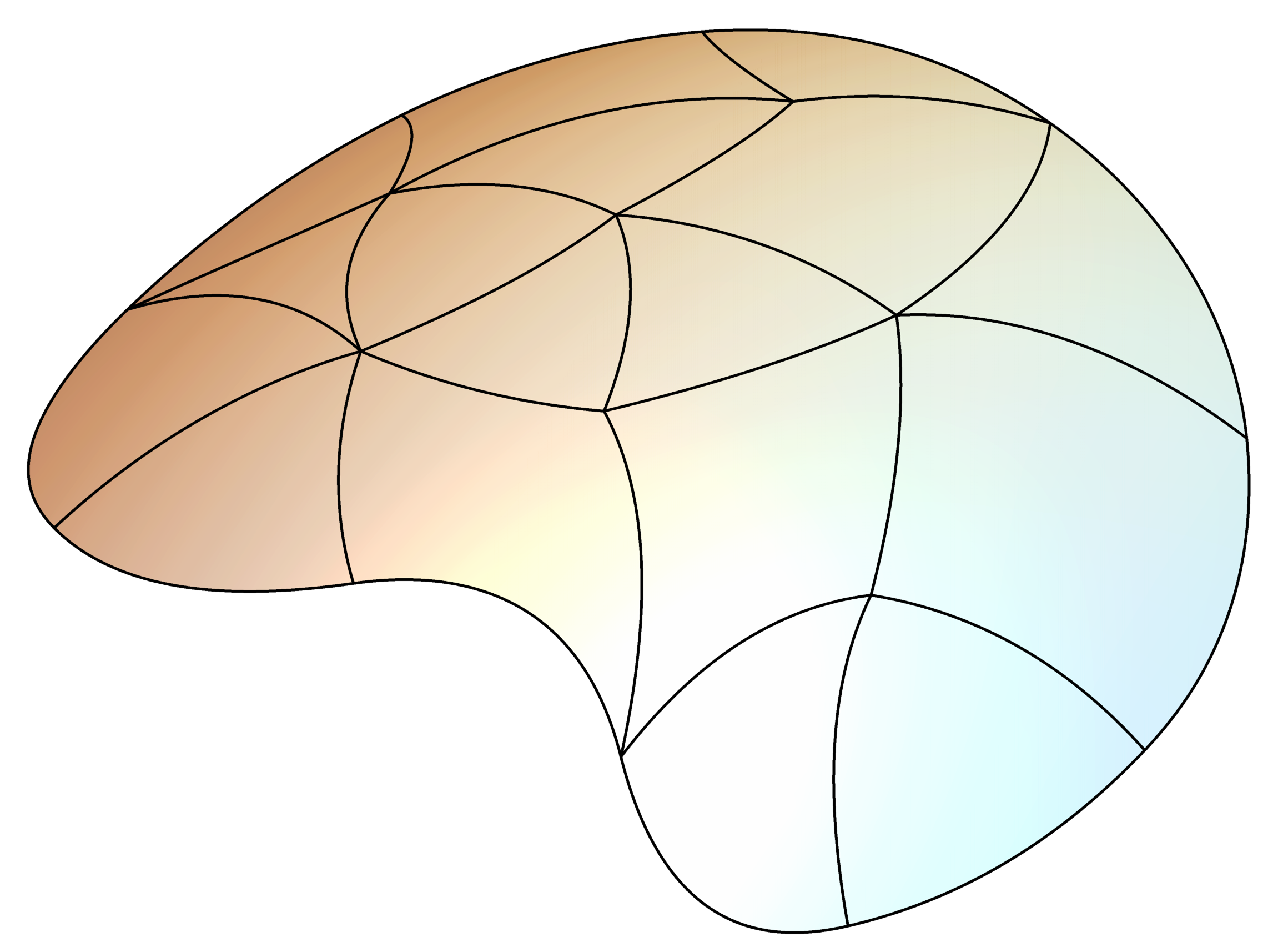}
\end{minipage}
\begin{minipage}{0.32\textwidth}
\centering
\includegraphics[width=1\textwidth]{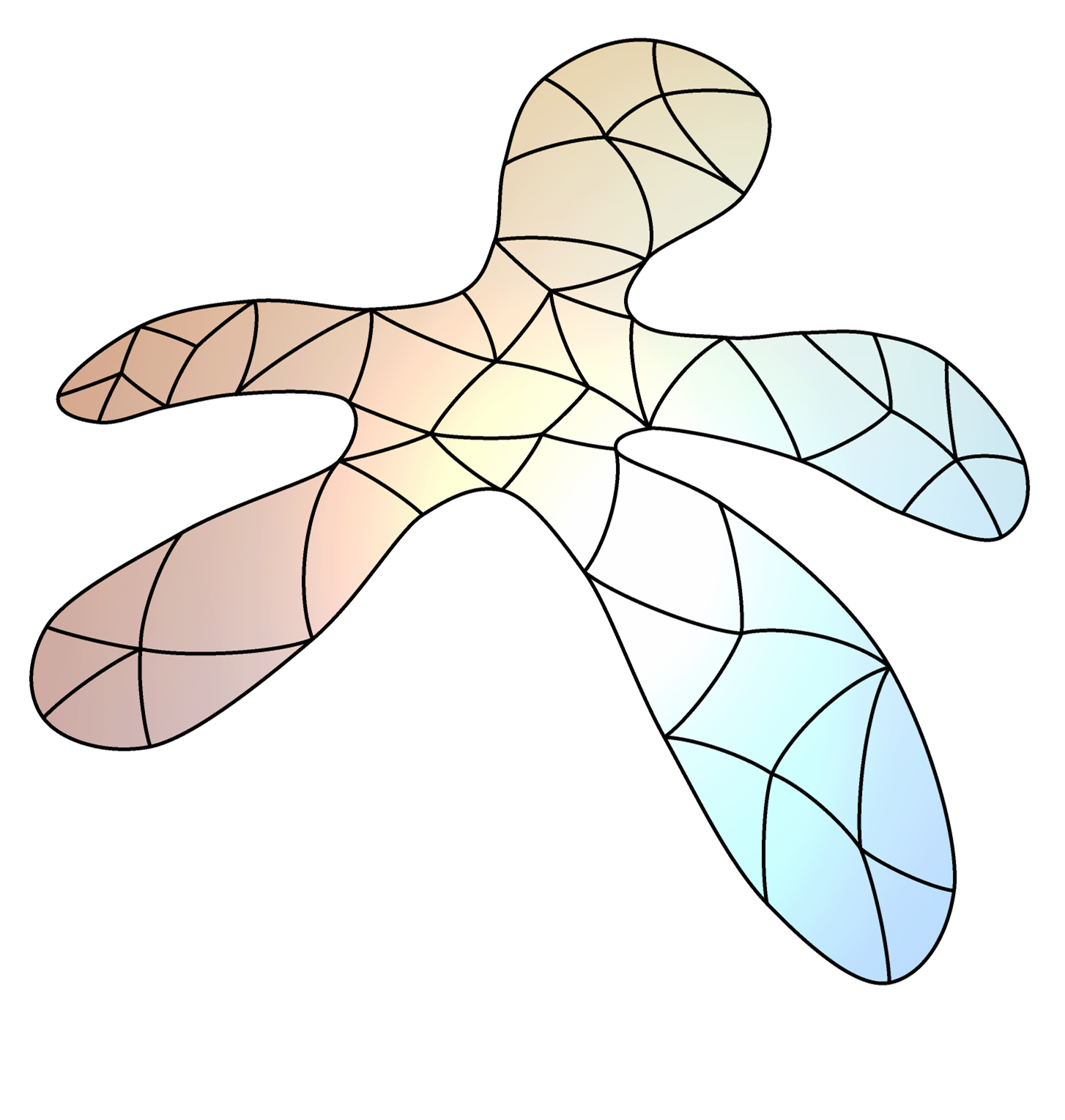}
\end{minipage}
\end{center}
\caption{Examples of mixed surfaces meshes.}
\label{fig:MixedSurfaceMesh}
\end{figure}

Assume that $\Omega \subset \RR^3$ is a connected surface, i.e., a $2$-dimensional manifold either having a boundary or being closed, that is $G^1$-smooth, i.e., it possesses a continuous tangent plane. A {\it mixed surface mesh} is a collection of triangular and quadrilateral polynomial surface elements, defined as follows.
\begin{definition}   
Let $\indT$ and $\indQ$ be two index sets, and let $\indOmega = \indT \cup \,\indQ$. 
For each $k \in \indOmega$, let $\Omega^{(k)}$ be a triangular or quadrilateral polynomial surface element, specified by its unique geometry mapping $\bfm{F}^{(k)}: \parDom^{(k)} \to {\Omega}^{(k)}$, 
$$
\parDom^{(k)} = \begin{cases}
\parDomT & \mbox{ if } k\in \indT\\
\parDomQ & \mbox{ if } k\in \indQ\\
\end{cases}, \quad 
\bfm{F}^{(k)} \in \begin{cases}
\left(\poly^2_{\delta}\right)^3 & \mbox{ if }k\in \indT\\
\left(\poly^2_{\delta, \delta}\right)^3 & \mbox{ if } k\in \indQ\\
\end{cases}.
$$
The collection $\{\Omega^{(k)}: \, k \in \indOmega\}$ is called a (conforming) mixed surface mesh of $\Omega$ if 
\begin{equation*}
{\Omega} = \bigcup_{k\in \indOmega} {\Omega}^{(k)}, 
\end{equation*}
the intersection of any two distinct elements is either empty, a common vertex, or a common edge. An edge is always the image of an entire edge of the parameter domain. In addition, we assume that the element parameterizations match in a $C^0$ sense (up to trivial orientation changes).
\end{definition}
Some examples of mixed surfaces meshes are shown in Figure~\ref{fig:MixedSurfaceMesh}.

Given a mixed surface mesh $\{\Omega^{(k)}: \, k \in \indOmega\}$, we denote  by $\set{\bE_\ell: \,\ell \in \indE}$ the set of edge curves, and by 
$\set{\bVm:\; m \in \indV}$ the set of vertices. Each edge possesses the parameterization $\bfm{E}_\ell: [0,1] \to \bE_\ell$, with $\bfm{E}_\ell\in (\poly^1_\delta)^3$. We call an edge an \emph{interface}, if it is contained in the interior of the domain, i.e., if it is contained in exactly two elements, otherwise it is called a \emph{boundary edge}.

For any $k\in \indOmega$, such that $\bE_\ell \subset {\Omega}^{(k)}$, the edge parameterization agrees with the element parameterization on a suitable restriction of the domain, i.e., 
\[
    \bfm{E}_\ell = \bfm{F}^{(k)} \circ \Qkl{k}{\ell},
\]
where $\Qkl{k}{\ell}:[0,1]\rightarrow \parDom^{(k)}$ is a linear function of the form
\[
    \Qkl{k}{\ell}(\xi) = \underline{v}^{(k)}_{\,\ell,1} (1-\xi) + \underline{v}^{(k)}_{\,\ell,2}\ \xi \quad \mbox{for some}\quad  \underline{v}^{(k)}_{\,\ell,1}, \underline{v}^{(k)}_{\,\ell,2} \in \{(0,0)^T,(0,1)^T,(1,0)^T,(1,1)^T\}.
\]
For the interface ${\bE}_\ell={\Omega}^{(k_1)}\cap{\Omega}^{(k_2)}$, it holds that
\begin{equation} \label{eq-C0-for-Qmega}
\bfm{F}^{(k_1)}\circ \Qkl{k_1}{\ell} = \bfm{F}^{(k_2)}\circ\Qkl{k_2}{\ell}.
\end{equation}
We define the tangent space to the surface as 
\begin{equation*}
    T\Omega := \defset{(\xBo,T_{\xBo}\Omega)}{\xBo \in \Omega},
\end{equation*}
where, for $\xBo = \bfm{F}^{(k)}(\xiBo)$ and $\bfm{N}^{(k)}(\xiBo) = 
(\Da \bfm{F}^{(k)} \times \Db \bfm{F}^{(k)})(\xiBo)$,  
we have
\begin{equation*}
    T_{\xBo}\Omega := \defset{\mathbf{y} \in \RR^3}{\mathbf{y} \perp \bfm{N}^{(k)}(\xiBo)}.
\end{equation*}
The corresponding tangent plane is the affine space $\xBo + T_{\xBo}\Omega$. The $G^1$-smoothness of the surface $\Omega$ implies that the tangent space and tangent plane are well-defined for all points $\xBo\in\Omega$. Thus, along any interface ${\bE}_\ell={\Omega}^{(k_1)}\cap{\Omega}^{(k_2)}$ the unit surface normals coincide up to the sign, i.e.,
\begin{equation} \label{eq-G1-for-Qmega}
\frac{\bfm{N}^{(k_1)}_\ell}{\norm{\bfm{N}^{(k_1)}_\ell}}= \pm \frac{\bfm{N}^{(k_2)}_\ell}{\norm{\bfm{N}^{(k_2)}_\ell}},
\end{equation}
where 
$\bfm{N}^{(k)}_\ell := \bfm{N}^{(k)} \circ \Qkl{k}{\ell}$ are the restrictions of the normals $\bfm{N}^{(k)}$, with $k\in\{k_1,k_2\}$, to the edge ${\bE}_\ell$.
Without loss of generality, we assume that both elements are oriented the same way, i.e., that the sign in~\eqref{eq-G1-for-Qmega} is positive.

The next two lemmas provide further insights needed in the subsequent analysis. 
\begin{lemma} \label{lemma-1}
Suppose that $\bE_\ell$ is the interior edge shared by ${\Omega}^{(k_1)}$ and ${\Omega}^{(k_2)}$. Let $\bfm{t}_\ell:[0,1]\rightarrow \RR^3$ be the tangent vectors at points on the edge ${\bE}_\ell$, with
\begin{subequations} \label{eq-def-tk-nk}
\begin{equation} \label{eq-def-tk}
\bfm{t}_\ell(\xi) := \Dt{\xi} \bfm{F}^{(k_1)}\left(\Qkl{k_1}{\ell}(\xi)\right) , \quad \xi\in [0,1].
\end{equation}
We have
\[
\bfm{t}_\ell(\xi) = J \bfm{F}^{(k_1)} \left( \Qkl{k_1}{\ell}(\xi)\right) \,\tkl{k_1}{\ell} =  J \bfm{F}^{(k_2)} \left( \Qkl{k_2}{\ell}(\xi)\right) \,\tkl{k_2}{\ell},
\]
where 
$J \bfm{F}^{(k)} = \left(\begin{array}{cc}
\partial_1 \bfm{F}^{(k)} &
\partial_2 \bfm{F}^{(k)}
\end{array}\right)$ and 
$\tkl{k}{\ell} = 
\Dt{\xi} \Qkl{k}{\ell}(\xi)$ is the constant parametric direction vector fulfilling
\begin{equation*}
\tkl{k}{\ell} \in  \left\{
(1,0)^T, (-1,0)^T, (0,1)^T, (0,-1)^T, (1,-1)^T, (-1,1)^T
\right\}.
\end{equation*}
Moreover, the in-plane normal vectors $\bfm{n}^{(k)}_\ell:[0,1]\rightarrow \RR^3$ are 
\begin{equation} \label{eq-def-nk}
\bfm{n}^{(k)}_\ell(\xi) := 
J \bfm{F}^{(k)}\left( \Qkl{k}{\ell}(\xi)\right){\nkl{k}{\ell} }, \quad k\in \set{k_1,k_2},
\end{equation}
with the (parametric) normals 
\begin{equation*}
\nkl{k}{\ell} = \left(\begin{array}{cc}
    0 & 1 \\ -1 & 0
\end{array}\right)\tkl{k}{\ell}.
\end{equation*}
\end{subequations}
For any $k\in\set{k_1,k_2}$ and $\xBo=\bfm{E}_\ell(\xi)$, the two vectors $\bfm{n}^{(k)}_\ell(\xi)$ and $\bfm{t}_\ell(\xi)$ span $T_{\xBo}\Omega$, i.e.,
\begin{equation*} 
T_{\xBo}\Omega = {\rm span}\set{\bfm{t}_\ell(\xi), \bfm{n}_\ell^{(k)}(\xi)}.  
\end{equation*}
Let 
\begin{equation}\label{eq-def-of-gluing-functions}
\alpha_\ell(\xi) := \det\big(\bfm{n}^{(k_2)}_\ell , \bfm{n}^{(k_1)}_\ell,\bfm{N}^{(k_1)}_\ell\big)(\xi), \quad \mbox{and}\quad \alpha_\ell^{(k)}(\xi) := \tfrac{1}{\gamma_\ell(\xi)} \widehat{\alpha}_\ell^{(k)}(\xi), \quad k \in \set{k_1,k_2},
\end{equation}
where
\begin{equation}\label{eq-def-of-alpha-hat}
\widehat{\alpha}_\ell^{(k)}(\xi) := 
{\norm{\tkl{k}{\ell}}^2 \,\sprod{\bfm{N}^{(k)}_\ell(\xi)}{\bfm{N}^{(k_1)}_\ell(\xi)}}, \quad k\in \set{k_1, k_2}, \quad\mbox{and} \quad \gamma_\ell(\xi) :=\gcd\left(\widehat{\alpha}_\ell^{(k_1)}(\xi), \widehat{\alpha}_\ell^{(k_2)}(\xi)\right).
\end{equation}
Here $\gcd$ is the polynomial greatest common divisor, which, by convention, is normalized such that its leading coefficient is $1$. 
The $G^1$ continuity condition \eqref{eq-G1-for-Qmega} is equivalent to 
\begin{equation} \label{eq-G1-for-Omega} 
\gamma_\ell(\xi) \ \alpha_\ell^{(k_1)}(\xi) \ \bfm{n}^{(k_2)}_\ell(\xi) - \gamma_\ell(\xi) \ \alpha_\ell^{(k_2)}(\xi) \ \bfm{n}^{(k_1)}_\ell(\xi) + \alpha_\ell(\xi) \ \bfm{t}_\ell(\xi) \equiv 0.
\end{equation}
We further have
\begin{equation} \label{eq-def-of-gluing-functions-conn}
\alpha_\ell(\xi) \beta_\ell(\xi)= \gamma_\ell(\xi)\alpha_\ell^{(k_2)}(\xi) \beta_\ell^{(k_1)}(\xi) - \gamma_\ell(\xi)\alpha_\ell^{(k_1)}(\xi) \beta_\ell^{(k_2)}(\xi), \quad \xi \in [0,1],
\end{equation}
with 
\begin{equation} \label{eq-def-of-gluing-functions-beta}
\beta_\ell(\xi) := \norm{\bfm{t}_\ell(\xi)}^2 \quad \mbox{and}\quad \beta_\ell^{(k)}(\xi) := \sprod{\bfm{n}^{(k)}_\ell(\xi)}{\bfm{t}_\ell(\xi)}, \quad k\in \set{k_1, k_2}.     
\end{equation}
Functions in \eqref{eq-def-of-gluing-functions} and \eqref{eq-def-of-gluing-functions-beta} are called the gluing functions.
\end{lemma}
\begin{proof}
Vectors $\bfm{n}^{(k)}_\ell(\xi)$, $k\in \set{k_1,k_2}$, clearly lie in the tangent space $$T_{\bfm{E}_\ell(\xi)}\Omega$$ for each $\xi \in [0,1]$ 
and they are not parallel to $\bfm{t}_\ell(\xi)$ since 
\begin{equation}\label{eq-nklcrosstkl}
\bfm{n}^{(k)}_\ell(\xi) \times \bfm{t}_\ell(\xi) = \norm{\tkl{k}{\ell}}^2  \left(\Da \bfm{F}^{(k)} \times \Db \bfm{F}^{(k)}\right)\left(\Qkl{k}{\ell}(\xi)\right) = {\norm{\tkl{k}{\ell}}^2  \bfm{N}^{(k)}_\ell(\xi)} \not = \bfm{0}.
\end{equation}
Thus, the $G^1$ continuity condition \eqref{eq-G1-for-Qmega} is equivalent to 
\begin{align*}
0\equiv\det\left(\bfm{n}^{(k_2)}_\ell(\xi), \bfm{n}^{(k_1)}_\ell(\xi), \bfm{t}_\ell(\xi)\right) = 
\widehat\alpha_\ell^{(k_1)}(\xi) \bfm{n}^{(k_2)}_\ell(\xi) - \widehat\alpha_\ell^{(k_2)}(\xi) \bfm{n}^{(k_1)}_\ell(\xi) + \alpha_\ell(\xi) \bfm{t}_\ell(\xi) 
\end{align*}
for some functions $\widehat\alpha_\ell^{(k_1)}$, $\widehat\alpha_\ell^{(k_2)}$ and $\alpha_\ell$. 
Multiplying this equality by $\sprod{\cdot \times \bfm{t}_\ell}{\bfm{N}_\ell^{(k_1)}}$ and by 
$\sprod{\cdot \times \bfm{n}^{(k_1)}_\ell}{\bfm{N}_\ell^{(k_1)}}$, and using \eqref{eq-nklcrosstkl},
it is straightforward to see that these functions are equal to the ones defined in~\eqref{eq-def-of-gluing-functions}--\eqref{eq-def-of-alpha-hat}. 
Multiplying \eqref{eq-G1-for-Omega} by $\sprod{\,\cdot\,}{\bfm{t}_\ell}$ 
yields the equality \eqref{eq-def-of-gluing-functions-conn}. 
\end{proof}

\begin{lemma} \label{lemma-1-1}
Suppose that $\bE_\ell$ is the interior edge shared by ${\Omega}^{(k_1)}$ and ${\Omega}^{(k_2)}$.
Along that interface, we have
\begin{equation}\label{eq:def-direction-d-ell}
\frac{\beta_\ell(\xi) \bfm{n}^{(k_1)}_\ell(\xi) - \beta_\ell^{(k_1)}(\xi)\bfm{t}_\ell(\xi)}{\alpha_\ell^{(k_1)}(\xi) \ } \equiv \frac{\beta_\ell(\xi) \bfm{n}^{(k_2)}_\ell(\xi) - \beta_\ell^{(k_2)}(\xi)\bfm{t}_\ell(\xi)}{\alpha_\ell^{(k_2)}(\xi)} =: \bfm{d}_\ell(\xi), \quad \xi \in [0,1].
\end{equation}
Moreover, the vector $\bfm{d}_\ell(\xi)$ is orthogonal to both $\bfm{t}_\ell(\xi)$ and the unit surface normal  $\frac{\bfm{N}^{(k_1)}_\ell(\xi)}{\norm{\bfm{N}^{(k_1)}_\ell(\xi)}}$. 
\end{lemma}
\begin{proof}
Identity~\eqref{eq:def-direction-d-ell} follows directly from \eqref{eq-G1-for-Omega}--\eqref{eq-def-of-gluing-functions-conn}. Using \eqref{eq-def-of-gluing-functions-beta}, we obtain that  $\sprod{\bfm{d}_\ell(\xi)}{\bfm{t}_\ell(\xi)}=0$ and hence $\bfm{d}_\ell(\xi)$ is orthogonal to
$\bfm{t}_\ell(\xi)$.  Moreover, $\bfm{d}_\ell(\xi)$ lies in the tangent space $T_{\bfm{E}_\ell(\xi)}\Omega$
and is therefore orthogonal to the surface normal.    
\end{proof}

\subsection{Isogeometric spline space} 

The isogeometric spline space $\mathcal{V}_\pd(\Omega)$ of degree $\pd \geq \delta$ over the surface domain $\Omega$ is defined as 
\[
 \mathcal{V}_\pd(\Omega) := 
 \defset{\varphi : {\Omega}\rightarrow \RR}{\varphi \circ \bfm{F}^{(k)} = f^{(k)} \in
 \begin{cases}
\poly^2_{\pd,\pd} & \mbox{ if } 
\parDom^{(k)}= \parDomQ\\
\poly^2_{\pd} & \mbox{ if } 
\parDom^{(k)}= \parDomT
 \end{cases}, \mbox{ for all } k \in \indOmega\;
}.
\]
To distinguish between the triangular and quadrilateral elements and to simplify further notation, we define
\begin{equation} \label{def-varsigma}
\varsigma^{(k)} := \begin{cases}
0,&\parDom^{(k)}=\parDomT\\
1,&\parDom^{(k)}=\parDomQ
\end{cases}.
\end{equation}
In this paper we study the $\C{1}$-smooth isogeometric subspace, defined as $\mathcal{V}^1_d(\Omega) := \mathcal{V}_d(\Omega) \cap \C{1}(\Omega)$. 
An isogeometric function is $C^1$-smooth at an interface $\bE_\ell$ if and only if it is continuous and the derivatives in a suitable crossing direction along the interface are continuous. The next theorem provides these conditions expressed in local coordinates. 
\begin{theorem}\label{thm-C1ForTwoElements}
An isogeometric function $\varphi \in \mathcal{V}_\pd(\Omega)$ is $C^1$-smooth, if and only if 
\begin{equation}\label{eq:C0-condition-interface}
    f^{(k_1)} \circ \Qkl{k_1}{\ell} (\xi) \equiv f^{(k_2)} \circ \Qkl{k_2}{\ell} (\xi)
\end{equation}
and
\begin{equation}\label{eq:C1-condition-interface}
    (\nabla f^{(k_1)} \circ \Qkl{k_1}{\ell})(\xi) \ \dkl{k_1}{\ell} \equiv (\nabla f^{(k_2)} \circ \Qkl{k_2}{\ell})(\xi) \ \dkl{k_2}{\ell},
\end{equation}
for all interfaces $\bE_\ell = \Omega^{(k_1)}\cap\Omega^{(k_2)}$, where $f^{(k)}=\varphi \circ \bfm F^{(k)}$ and the directions $\dkl{k}{\ell}$ are defined as 
\begin{equation}\label{eq:def-d-l-k}
 \dkl{k}{\ell}(\xi) := \frac{\beta_\ell(\xi) \nkl{k}{\ell} - \beta_\ell^{(k)}(\xi) \tkl{k}{\ell}}{\alpha_\ell^{(k)}(\xi) \ }.
\end{equation}
\end{theorem}
\begin{proof}
Consider an interface $\bE_\ell = \Omega^{(k_1)}\cap\Omega^{(k_2)}$.
The $\C{0}$-smoothness of the isogeometric function $\varphi$ along $\bE_\ell$ is equivalent to~\eqref{eq:C0-condition-interface}. Along that interface, recall the identity~\eqref{eq:def-direction-d-ell}, which is by~\eqref{eq-def-tk-nk} equivalent to
\begin{equation}\label{eq:direction-d-ell-with-Jacobian}
 (J \bfm F^{(k_1)}\circ \Qkl{k_1}{\ell}) \ \dkl{k_1}{\ell} \equiv (J \bfm F^{(k_2)}\circ \Qkl{k_2}{\ell})\ \dkl{k_2}{\ell} \equiv \bfm{d}_\ell.
\end{equation}

In order to analyze the $\C{1}$-smoothness, we extend the boundary lines $\{\Qkl{k}{\ell}(\xi) \, : \, \xi\in [0,1]\}$, $k \in \set{k_1,k_2}$, of the parametric domains $\parDom^{(k)}$ to three-dimensional regions 
\[
\widetilde{\parDom}^{(k)}_\ell = \defset{\xiBo\in \RR^2}{\|\xiBo-\Qkl{k}{\ell}(\xi)\| \leq \epsilon, \;\xi \in [0,1]} \times[-\epsilon,\epsilon], 
\]
and the geometry mappings $\bfm{F}^{(k)}$, $k \in \set{k_1,k_2}$,  
to 
\begin{equation} \label{eq-ext-of-F}
\widetilde{\bfm{F}}^{(k)}: \widetilde{\parDom}^{(k)}_\ell \to \widetilde{\bE}^{(k)}_{\ell} \subset\RR^3, \quad \widetilde{\bfm{F}}^{(k)}(\widetilde{\xiBo}) = \bfm{F}^{(k)}(\xiBo) + 
\xic  \bfm{N}^{(k)}(\xiBo)/\norm{\bfm{N}^{(k)}(\xiBo)},
\end{equation}
for $\epsilon$ small enough so that $\widetilde{\bfm{F}}^{(k)}$ is bijective and regular, where $\widetilde{\xiBo} = (\xiBo,\xic)=(\xia,\xib,\xic)$. 
Moreover, the isogeometric function $\varphi$ and the corresponding polynomials $f^{(k)} = \varphi\circ \bfm{F}^{(k)}$, $k\in \set{k_1,k_2}$, are extended as 
\[
\widetilde{f}^{(k)}(\widetilde{\xiBo}) := f^{(k)}(\xiBo) \quad \mbox{and} \quad 
\widetilde{\varphi}^{(k)}|_{\widetilde{\bE}^{(k)}_{\ell}} := \widetilde{f}^{(k)}\circ {\left(\widetilde{\bfm{F}}^{(k)}\right)^{-1}}|_{\widetilde{\bE}^{(k)}_{\ell}}.
\]
Clearly, it holds true that $\widetilde{\varphi}^{(k_1)}|_{\widetilde{\bE}_\ell} = \widetilde{\varphi}^{(k_2)}|_{\widetilde{\bE}_\ell}$, where $\widetilde{\bE}_\ell := \widetilde{\bE}^{(k_1)}_\ell\cap \widetilde{\bE}^{(k_2)}_\ell$.

Let $\xBo\in\bE_\ell$, with $\xBo=\bfm{E}_\ell(\xi)$, for some $\xi\in[0,1]$. 
The isogeometric function $\varphi$ is $\C{1}$-smooth at $\xBo$ if and only if the condition 
\begin{equation} \label{eq-cond-ext-C1}
\Dv_{\bfm{d}_\ell(\xi)}\widetilde{\varphi}^{(k_1)}(\xBo) = \Dv_{\bfm{d}_\ell(\xi)}\widetilde{\varphi}^{(k_2)}(\xBo)
\end{equation}
is satisfied, where $\bfm{d}_\ell(\xi)$ is as defined in~\eqref{eq:def-direction-d-ell} and $D_{\bfm{d}_\ell(\xi)}$ denotes the directional derivative in direction $\bfm{d}_\ell(\xi)$. Let $\nabla_{\widetilde{\xiBo}}$ denote the gradient operator with respect to the parameters $\widetilde{\xiBo}$, considered as a row vector. From the equality $\widetilde{f}^{(k)}(\widetilde{\xiBo}) = \widetilde{\varphi}^{(k)}\left(\widetilde{\bfm{F}}^{(k)}(\widetilde{\xiBo})\right)$
it is straightforward to compute that
$$
\nabla_{\widetilde{\xiBo}} \widetilde{f}^{(k)}(\widetilde{\xiBo}) = \nabla_{\xBo} \widetilde{\varphi}^{(k)}(\xBo) J \widetilde{\bfm{F}}^{(k)}(\widetilde{\xiBo}), 
$$
where $\xBo = \widetilde{\bfm{F}}^{(k)}(\widetilde{\xiBo})$, with $\widetilde{\xiBo} = (\Qkl{k}{\ell}(\xi),0)^T$, and $J \widetilde{\bfm{F}}^{(k)}$ is the Jacobian matrix of $\widetilde{\bfm{F}}^{(k)}$. Thus, we have
\[
\Dv_{\bfm{d}_\ell(\xi)}\widetilde{\varphi}^{(k)}(\xBo) = \nabla_{\xBo} \widetilde{\varphi}^{(k)}(\xBo)\bfm{d}_\ell(\xi) 
= \nabla_{\widetilde{\xiBo}} \widetilde{f}^{(k)}(\widetilde{\xiBo}) (J\widetilde{\bfm{F}}^{(k)}(\widetilde{\xiBo}))^{-1}\bfm{d}_\ell(\xi).
\]
From~\eqref{eq:direction-d-ell-with-Jacobian}, we further get
\begin{equation}\label{eq:directional-derivative-proof}
\Dv_{\bfm{d}_\ell(\xi)}\widetilde{\varphi}^{(k)}(\xBo) = \nabla_{\widetilde{\xiBo}} \widetilde{f}^{(k)}(\widetilde{\xiBo}) (J\widetilde{\bfm{F}}^{(k)}(\widetilde{\xiBo}))^{-1}\bfm{d}_\ell(\xi) = \nabla_{\widetilde{\xiBo}} \widetilde{f}^{(k)}(\widetilde{\xiBo}) (J\widetilde{\bfm{F}}^{(k)}(\widetilde{\xiBo}))^{-1}  J \bfm F^{(k)}( \Qkl{k}{\ell}(\xi)) \dkl{k}{\ell},
\end{equation}
for $k\in\{k_1,k_2\}$. Moreover, $\nabla_{\widetilde{\xiBo}} \widetilde{f}^{(k)}(\Qkl{k}{\ell}(\xi),0) = \left(\Da f^{(k)}(\Qkl{k}{\ell}(\xi)), \Db f^{(k)}(\Qkl{k}{\ell}(\xi)),0\right)$, and
\[
 J \bfm F^{(k)} (\Qkl{k}{\ell}(\xi)) \dkl{k}{\ell} = J \widetilde{\bfm{F}}^{(k)} (\widetilde{\xiBo}) \left(\begin{array}{c} \dkl{k}{\ell} \\ 0 \end{array}\right).
\]
Thus, \eqref{eq:directional-derivative-proof} simplifies to
\[
    \Dv_{\bfm{d}_\ell(\xi)}\widetilde{\varphi}^{(k)}(\xBo) = \left(\Da f^{(k)}(\Qkl{k}{\ell}(\xi)), \Db f^{(k)}(\Qkl{k}{\ell}(\xi)),0\right) \left(\begin{array}{c} \dkl{k}{\ell} \\ 0 \end{array}\right) = (\nabla f^{(k)} \circ \Qkl{k}{\ell})(\xi) \ \dkl{k}{\ell}.
\]
Consequently, \eqref{eq-cond-ext-C1} is equivalent to \eqref{eq:C1-condition-interface}, which completes the proof.
\end{proof}
In the following we analyze the $C^1$-smooth isogeometric space described in Theorem~\ref{thm-C1ForTwoElements} and provide a basis construction.

\section{A $\C{1}$ super-smooth isogeometric spline space over a surface domain}\label{sec:SplineSpace}

In this section we analyze the $C^1$-conditions across a single interface, as they are described in Theorem~\ref{thm-C1ForTwoElements}. We then develop a Hermite-like basis for all non-vanishing functions at an interface. This representation is then used to define a global $C^1$-smooth basis, which is additionally $C^2$ at all vertices.

\subsection{$\C{1}$ continuity conditions over an interface} \label{SubSec:C1ContinuityConditionsOverTheEdge}

For a single interface we have the following.
\begin{theorem} \label{thm-degree-of-omega}
Consider an interface $\bE_\ell = \Omega^{(k_1)}\cap\Omega^{(k_2)}$ and let $f^{(k)} = \varphi\circ \bfm{F}^{(k)}$, for $k\in\{k_1,k_2\}$, be the pull-backs of the isogeometric function $\varphi\in \mathcal{V}_\pd(\Omega)$. We define the \emph{trace function} $\theta_\ell:[0,1]\rightarrow \RR$ and the \emph{crossing derivative function} $\w_\ell:[0,1]\rightarrow \RR$ as
\begin{equation}\label{eq:def-theta-l}
    \theta_\ell(\xi) := f^{(k_1)} \circ \Qkl{k_1}{\ell} (\xi) \equiv f^{(k_2)} \circ \Qkl{k_2}{\ell} (\xi)
\end{equation}
and
\begin{equation}\label{eq:def-w-l}
    \w_\ell(\xi) := (\nabla f^{(k_1)} \circ \Qkl{k_1}{\ell})(\xi) \ \dkl{k_1}{\ell} \equiv (\nabla f^{(k_2)} \circ \Qkl{k_2}{\ell})(\xi) \ \dkl{k_2}{\ell},
\end{equation}
respectively, where
\begin{align}
& \w_\ell(\xi) = 
\frac{1}{{\alpha}_\ell^{(k)}(\xi)} 
\left(\beta_\ell(\xi) \,\vartheta^{(k)}_\ell(\xi) - \beta_\ell^{(k)}(\xi)\, \theta'_\ell(\xi) \right), \quad 
\quad \vartheta^{(k)}_\ell(\xi) := \nabla f^{(k)}\left(\Qkl{k}{\ell}(\xi)\right) 
{\nkl{k}{\ell}}, 
\label{eq-vartheta}
\end{align}
for $k\in\{k_1,k_2\}$, with gluing functions defined in \eqref{eq-def-of-gluing-functions} 
and \eqref{eq-def-of-gluing-functions-beta}. The trace function $\theta_\ell$ is a polynomial of degree at most $d$ and the crossing derivative function {$\w_\ell$} is a polynomial of degree at most {$\pd+2\delta-2$}.
\end{theorem}
\begin{proof}
Relation \eqref{eq-vartheta} is a direct consequence of Theorem~\ref{thm-C1ForTwoElements}, in particular of~\eqref{eq:C1-condition-interface} and the definition of the directions $\dkl{k}{\ell}$ in~\eqref{eq:def-d-l-k}. 
It follows directly from~\eqref{eq:def-theta-l} that, since $f^{(k)} \in \poly^2_{\pd,\pd}$ or $f^{(k)} \in \poly^2_{\pd}$, the function $\theta_\ell$ is a polynomial in $\poly_{d}$.
With the same argument, $\vartheta^{(k)}_\ell$ is also a polynomial {of degree at most $d$}. Inserting~\eqref{eq-vartheta} into~\eqref{eq:def-w-l} and recalling $\gcd\left({\alpha}_\ell^{(k_1)}, {\alpha}_\ell^{(k_2)}\right)=1$, 
it is clear from the polynomial equality  
\[
{\alpha}_\ell^{(k_2)} 
\left(\beta_\ell \,\vartheta^{(k_1)}_\ell - \beta_\ell^{(k_1)}\, \theta'_\ell \right) = {\alpha}_\ell^{(k_1)} 
\left(\beta_\ell \,\vartheta^{(k_2)}_\ell - \beta_\ell^{(k_2)}\, \theta'_\ell \right),
\]
that ${\alpha}_\ell^{(k)}$ must divide $\beta_\ell \,\vartheta^{(k)}_\ell - \beta_\ell^{(k)}\, \theta'_\ell$ for $k\in \{k_1, k_2\}$, so $\w_\ell$ must be a polynomial. Moreover, it is straightforward to see that  
\begin{equation} \label{eq-degrees-1}
\deg{(\bfm{t}_\ell)} \leq \delta-1, \quad \deg{(\bfm{n}^{(k)}_\ell)} \leq \delta-1 + \varsigma^{(k)}, \quad
\deg(\beta_\ell) \leq 2 \delta -2, \quad \deg(\beta^{(k)}_\ell) \leq 2 \delta -2 + \varsigma^{(k)},
\end{equation}
where $\varsigma^{(k)}$ is defined by \eqref{def-varsigma}. 
Thus, the maximal possible degree of $\w_\ell$ is 
{$\pd+2\delta-2$, which concludes the proof.}
\end{proof}
The following remark establishes a connection between the results of Theorem~\ref{thm-C1ForTwoElements} and those in \cite{GrPe15}, where smoothness is studied via the geometric continuity of graphs of isogeometric functions in $\RR^4$. In contrast to \cite{GrPe15}, Theorem~\ref{thm-C1ForTwoElements} provides additional insight by explicitly characterizing the directional derivatives along the interface. This explicit representation is particularly useful for further analysis and for the construction of spline spaces with prescribed smoothness properties.
\begin{remark}
The graph $\Phi \subset \Omega \times \RR\subset\RR^4$ 
of an isogeometric function $\varphi\in \mathcal{V}_\pd(\Omega)$ consists of patches in $\RR^4$, given by parameterizations
\begin{equation*}
\Phi^{(k)} = \begin{bmatrix}
\bfm{F}^{(k)}\\
f^{(k)}
\end{bmatrix}: \parDom^{(k)} \to \RR^4.
\end{equation*}
Theorem~\ref{thm-C1ForTwoElements} is closely related to the result proven in \cite{GrPe15}, that along all interfaces the function $\varphi$ is $\C{1}$-smooth if and only if its graph $\Phi$ is $\G{1}$-smooth. 
Namely, by \eqref{eq-C0-for-Qmega}, \eqref{eq-G1-for-Omega}, \eqref{eq-def-of-gluing-functions-conn}, \eqref{eq:C0-condition-interface} and \eqref{eq:C1-condition-interface}, it is straightforward to verify that for any interface $\bE_\ell = \Omega^{(k_1)}\cap\Omega^{(k_2)}$, we have 
\begin{align*}
\Phi^{(k_1)}\circ \Qkl{k_1}{\ell} = \Phi^{(k_2)}\circ \Qkl{k_2}{\ell}\quad \mbox{and} \quad \gamma_\ell\, \alpha_\ell^{(k_1)} \,\widehat{\bfm{n}}^{(k_2)}_\ell - \gamma_\ell \, \alpha_\ell^{(k_2)} \, \widehat{\bfm{n}}^{(k_1)}_\ell + \alpha_\ell \,\widehat{\bfm{t}}_\ell \equiv 0,
\end{align*}
where $\widehat{\bfm{t}}_\ell$, $\widehat{\bfm{n}}^{(k)}_\ell$, $k\in\{k_1,k_2\}$, are defined as in Lemma~\ref{lemma-1} with $\bfm{F}^{(k)}$ replaced by $\Phi^{(k)}$.
\end{remark}
The results presented in the following theorem play a crucial role in the construction of isogeometric spline spaces.

\begin{theorem} \label{thm-3}
Suppose that $\bVm \in {\Omega}^{(k)}$, $\bVm = \bfm{F}^{(k)}(\xiBo_m)$ for some $\xiBo_m\in \set{(0,0),(1,0),(0,1),(1,1)}$. Let 
$\bfm{E}_{m,1}=\left({E}_{m,1,1},{E}_{m,1,2},{E}_{m,1,3}\right)^T$ and $\bfm{E}_{m,2}=\left({E}_{m,2,1},{E}_{m,2,2},{E}_{m,2,3}\right)^T$ be two orthonormal vectors that span the tangent space $T_{\bVm}\Omega$. Let $\phi:\RR^3\rightarrow \RR$ be an arbitrary $C^2$-smooth function.
Given the values 
\begin{equation}\label{int-data-at-V}
\Dv_{\bfm{E}_{m,2}}^{\eta_2} \Dv_{\bfm{E}_{m,1}}^{\eta_1}\phi(\bVm) = c_{m,\eta_1,\eta_2}, \quad 0\leq \eta_1+\eta_2\leq 2,
\end{equation}
the values $\Db^{\eta_2} \Da^{\eta_1} f^{(k)}(\xiBo_m)$, $0\leq \eta_1+\eta_2\leq 2$, are uniquely determined by interpolating
\begin{equation*}
    \lim_{\Omega^{(k)}\ni\xBo \rightarrow \bVm}\Dv_{\bfm{E}_{m,2}}^{\eta_2} \Dv_{\bfm{E}_{m,1}}^{\eta_1}f^{(k)}\left((\bfm{F}^{(k)})^{-1}(\xBo)\right) = c_{m,\eta_1,\eta_2}.
\end{equation*}
\end{theorem}
\begin{proof}
Clearly 
\begin{equation} \label{eq-int-cond-Vm-0}
f^{(k)}(\xiBo_m) = \varphi\big(\bfm{F}^{(k)}(\xiBo_m)\big) = \varphi\big(\bVm\big) = c_{m,0,0}.    
\end{equation}
By considering in the neighborhood of the vertex $\bVm$ the extensions $\widetilde{\varphi}^{(k)}$,
$\widetilde{f}^{(k)}$ and $\widetilde{\bfm{F}}^{(k)} = (\widetilde{F}_1^{(k)},
\widetilde{F}_2^{(k)},\widetilde{F}_3^{(k)})^T$, as defined in the proof of Theorem~\ref{thm-C1ForTwoElements}, 
it is straightforward to compute from $\widetilde{f}^{(k)} = \widetilde{\varphi}^{(k)} \circ \widetilde{\bfm{F}}^{(k)}$ that
\begin{align*}
&
\nabla_{\widetilde{\xiBo}} \widetilde{f}^{(k)}(\widetilde{\xiBo}_{m}) = \nabla_{\bfm{x}}  \widetilde{\varphi}^{(k)}(\bVm) J_m, 
\\  
& H\widetilde{f}^{(k)}(\widetilde{\xiBo}_{m}) =J_m^T H_m \,J_m + 
\widetilde{\varphi}_{x_1}^{(k)}(\bVm) H \widetilde{F}^{(k)}_1(\widetilde{\xiBo}_{m}) + 
\widetilde{\varphi}_{x_2}^{(k)}(\bVm) H \widetilde{F}^{(k)}_2(\widetilde{\xiBo}_{m}) +
\widetilde{\varphi}_{x_3}^{(k)}(\bVm) H \widetilde{F}^{(k)}_3(\widetilde{\xiBo}_{m}),
\end{align*}
where
$\widetilde{\xiBo}_{m}= (\xiBo_m,0)=(\xi_{m,1},\xi_{m,2},0)$, 
$$J_m = J \widetilde{\bfm{F}}^{(k)}(\widetilde{\xiBo}_{m}) =\left(
\Da \bfm{F}^{(k)}(\xiBo_m), \Db \bfm{F}^{(k)}(\xiBo_m), \bfm{N}_{m}
\right), \quad H_m = H \widetilde{\varphi}^{(k)}(\bVm) \in \RR^{3\times 3},
$$
and $\bfm{N}_{m}=\bfm{N}(\bfm{V}^{(m)})$ denotes the unit normal of the surface domain at the vertex $\bVm$.
Let 
$$
{Q=\left(\bfm{E}_{m,1}, \bfm{E}_{m,2}, \bfm{N}_m\right)} \in  SO(3), \quad
E= \begin{pmatrix}
\sprodS{\bfm{E}_{m,1}}{\Da \bfm{F}^{(k)}(\xiBo_m)}& \sprodS{\bfm{E}_{m,1}}{\Db \bfm{F}^{(k)}(\xiBo_m)}\\  
\sprodS{\bfm{E}_{m,2}}{\Da \bfm{F}^{(k)}(\xiBo_m)}& \sprodS{ \bfm{E}_{m,2}}{\Db \bfm{F}^{(k)}(\xiBo_m)}
\end{pmatrix},$$
and note that $$
Q^T J_m = 
\left(
\begin{array}{cc}
\multirow{2}{*}{E}  & 0 \\
                    & 0 \\
0 \quad 0 & 1 \\
\end{array}
\right).
$$
Since the directional derivatives (in the tangent plane) \eqref{int-data-at-V} are independent of the extension of ${\varphi}^{(k)}$ to $\widetilde{\varphi}^{(k)}$, we have that $\Dv^{\eta_2}_{\bfm{E}_{m,2}} \Dv^{\eta_1}_{\bfm{E}_{m,1}}\widetilde{\varphi}(\bVm) = \Dv^{\eta_2}_{\bfm{E}_{m,2}} \Dv^{\eta_1}_{\bfm{E}_{m,1}}\varphi(\bVm)$ for $0\leq \eta_1+\eta_2\leq 2$. Thus
\begin{align}
\sprod{\nabla  \widetilde{\varphi}^{(k)}(\bVm)}{\bfm{E}_{m,1}} = c_{m,1,0},\quad 
\sprod{\nabla  \widetilde{\varphi}^{(k)}(\bVm)}{\bfm{E}_{m,2}} = c_{m,0,1}, \label{eq-tmp-1}
\end{align}
and from    
\begin{align*}
\nabla \widetilde{f}^{(k)}(\widetilde{\xiBo}_{m}) & = (\Da f^{(k)}(\xiBo_m), \Db f^{(k)}(\xiBo_m), 0) = 
\nabla  \widetilde{\varphi}^{(k)}(\bVm) J_m = \left(\nabla  \widetilde{\varphi}^{(k)}(\bVm) \, Q\right) \left(Q^T J_m\right) \\
& = \left(c_{m,1,0}, c_{m,0,1}, \sprodS{\nabla  \widetilde{\varphi}^{(k)}(\bVm)}{\bfm{N}_{m}}\right) \left(Q^T J_m\right),
\end{align*}
we obtain that
\begin{equation} \label{eq-int-cond-Vm-1}
(\Da f^{(k)}(\xiBo_m), \Db f^{(k)}(\xiBo_m)) = (c_{m,1,0}, c_{m,0,1})\, E \quad {\rm and}\quad \sprodS{\nabla  \widetilde{\varphi}^{(k)}(\bVm)}{\bfm{N}_{m}}=0.    
\end{equation}
Let $M(1:2,1:2)$ denote a submatrix of a matrix $M$ with first two rows and first two columns. From \eqref{int-data-at-V} it follows
$$
\left(Q^T H_m \, Q\right)(1:2,1:2) = \begin{pmatrix}
c_{m,2,0}& c_{m,1,1} \\
c_{m,1,1} & c_{m,0,2} \end{pmatrix} = C, 
$$
and $J_m^T H_m  \, J_m = (Q^T J_m)^T (Q^T H_m \, Q) (Q^T J_m)$
gives $\left(J_m^T H_m J_m\right)(1:2,1:2) = E^T C E$. 
Moreover (see \eqref{eq-ext-of-F}), 
$$\left(H \widetilde{F}^{(k)}_i(\widetilde{\xiBo}_{m})\right)(1:2,1:2) = 
\begin{pmatrix}
\Da^2 F^{(k)}_i(\xiBo_m) & \Da \Db F^{(k)}_i(\xiBo_m)\\
\Da \Db F^{(k)}_i(\xiBo_m)& \Db^2 F^{(k)}_i(\xiBo_m)
\end{pmatrix} =: H_{m,i} \in \RR^{2\times 2}, \quad i =1,2,3, $$ 
and one can check using \eqref{eq-tmp-1} and $\sprodS{\nabla  \widetilde{\varphi}^{(k)}(\bVm)}{\bfm{N}_{m}}=0$ that
\begin{align*}
& \left(\Dx \widetilde{\varphi}^{(k)}(\bVm) H \widetilde{F}^{(k)}_1(\widetilde{\xiBo}_{m}) + 
\Dy \widetilde{\varphi}^{(k)}(\bVm) H \widetilde{F}^{(k)}_2(\widetilde{\xiBo}_{m}) +
\Dz \widetilde{\varphi}^{(k)}(\bVm) H \widetilde{F}^{(k)}_3(\widetilde{\xiBo}_{m})\right)(1:2,1:2) =\\ 
& \sum_{\ell=1}^3 \left(c_{m,1,0} {E}_{m,1,\ell} + c_{m,0,1} {E}_{m,2, \ell}\right) H_{m,\ell},  
\end{align*}
Therefore
\begin{equation} \label{eq-int-cond-Vm-2}
Hf^{(k)}(\xiBo_{m}) =  E^T C E + \sum_{i=1}^3 \left(c_{m,1,0} {E}_{m,1, i} + c_{m,0,1} {E}_{m,2,i}\right) H_{m,i},   
\end{equation} 
which completes the proof.
\end{proof}

From Theorem~\ref{thm-C1ForTwoElements} it follows that the $\C{1}$ continuity conditions
along the interface $\bE_\ell$ are fulfilled if and only if
\begin{subequations} \label{eq-cond-remAndHighOrdT}
\begin{align}
f^{(k)}({\Qkl{k}{\ell}}(\xi)) & = \theta_\ell(\xi), \nonumber\\
\beta_\ell(\xi) \vartheta^{(k)}_\ell(\xi) & = r_\ell^{(k)}(\xi), \quad 
r_\ell^{(k)}(\xi):= \alpha_\ell^{(k)}(\xi) \w_\ell(\xi) + \beta_\ell^{(k)}(\xi)\, \theta'_\ell(\xi), \quad \xi \in [0,1], \quad k \in \set{k_1,k_2},  \label{def-rik}
\end{align}
for appropriately chosen trace function $\theta_\ell \in \poly_\pd$ and crossing derivative function $\w_\ell$, which by Theorem~\ref{thm-degree-of-omega} is in $\poly_{\pd + 2\delta-2}$.
The main difficulty arises from the fact that, when 
$\beta_\ell$ is not a constant polynomial, 
the trace function $\theta_\ell$ and the crossing derivative function $\w_\ell$ can not, in general, be chosen independently, since $r_\ell^{(k)}$ must be divisible by $\beta_\ell$, i.e., the following remainders must be identically zero:
\begin{equation} \label{eq-cond-remainders}
\widehat{r}^{(k)}_\ell(\xi) :={\rm rem}\left(r^{(k)}_\ell(\xi),\beta_\ell(\xi)\right) \equiv 0,  \quad k \in \set{k_1,k_2}.  
\end{equation}
Furthermore, the degree of $r^{(k)}_\ell$ is limited by 
\begin{equation} \label{eq-cond-higherOrderTerms}
\deg\left( r_\ell^{(k)}\right) \leq d_{k,\ell}, \quad 
d_{k,\ell} := \deg\left(\beta_\ell \vartheta^{(k)}_\ell\right) =\deg(\beta_\ell) + \pd-1+ \varsigma^{(k)}, 
\quad k \in \set{k_1,k_2}.
\end{equation}
\end{subequations}
In \cite[Theorem~1]{GrKaKnTaVi2024} it is proven that conditions \eqref{eq-cond-remAndHighOrdT} are sufficient and necessary for the isogeometric function $\varphi$ to be $C^1$-continuous. 
However, these conditions do not directly indicate how to construct $\theta_\ell$ and $\w_\ell$, nor how to determine the number of degrees of freedom for the general degree $\delta$. This issue is resolved in \cite{GrKaKnTaVi2024} for the special case of meshes satisfying $\deg(\beta_\ell) \leq 2$, $\deg(\beta_\ell^{(k)}) \leq 2 + \varsigma^{(k)}$, and $\deg(\alpha_\ell^{(k)}) \leq 2 + \varsigma^{(k)}$, which, in particular, includes meshes with quadratic or biquadratic geometry mappings ($\delta = 2$). Under these assumptions, it is shown that the number of degrees of freedom associated with the interface is at least $2\pd - 1$.

\subsection{Hermite-like basis construction for an interface}\label{sec:Hermite-like-basis}

We consider a single interface with index $\ell$ between two elements with indices $k_1$ and $k_2$. To simplify the notation in this section, we drop the index $\ell$ and replace $k_1$ by $1$ and $k_2$ by $2$. Thus, we have for $k\in\{1,2\}$ gluing functions $\alpha^{(k)}$, $\beta^{(k)}$, and $\beta$, which determine the relations between functions $\theta$, $\w$, and $\vartheta^{(k)}$, describing the trace, crossing derivative, and parametric normal derivative on patch $k$, respectively. The degree of $\w$ is set to its maximal possible value $d+2\delta-2$.

We present a numerical algorithm that provides a constructive procedure for defining $\theta$ and $\w$ for geometry mappings of arbitrary degree $\delta$, such that they satisfy the required conditions \eqref{eq-cond-remAndHighOrdT} and, in addition, attain prescribed values and derivatives at the boundary parameters $\xi = 0$ and $\xi = 1$, up to order two and one, 
respectively:
\begin{equation*}
\Dti{\xi}{i} \tr(\xi)|_{\xi=\eta} =a_{\eta,i}, \quad \eta\in\set{0,1}, \; i \in\set{0,1,2}, \quad
\Dti{\xi}{i} \w(\xi)|_{\xi=\eta} =b_{\eta,i}, \quad \eta\in\set{0,1}, \; i \in\set{0,1}.
\end{equation*}
These conditions are automatically satisfied by expressing
\begin{align}
\tr(\xi)&  =   a_{0,0} \bern{5}{0}(\xi) + 
\left(a_{0,0} + \tfrac{1}{5} a_{0,1}\right) \bern{5}{1}(\xi) +
\left(a_{0,0} + \tfrac{2}{5} a_{0,1}+\tfrac{1}{20} a_{0,2}\right) \bern{5}{2}(\xi)\nonumber\\
& + \left(a_{1,0} - \tfrac{2}{5} a_{1,1} + \tfrac{1}{20} a_{1,2}\right) \bern{5}{3}(\xi) +
\left(a_{1,0} - \tfrac{1}{5} a_{1,1}\right) \bern{5}{4}(\xi)+a_{1,0} \bern{5}{5}(\xi) 
+ \xi^3(1-\xi)^3 \sum_{j=0}^{\pd-6}c_{0,j} \bern{\pd-6}{j}(\xi), \label{eq-def-tr-w-1}\\
\w(\xi) & =   b_{0,0} \bern{3}{0}(\xi) + 
\left(b_{0,0} + \tfrac{1}{3} b_{0,1}\right) \bern{3}{1}(\xi) +
\left(b_{1,0} - \tfrac{1}{3} b_{1,1}\right) \bern{3}{2}(\xi)+b_{1,0} \bern{3}{3}(\xi) \nonumber\\
& + \xi^2(1-\xi)^2 \sum_{j=0}^{{\pd + 2\delta -6}}c_{1,j} \bern{{\pd + 2\delta -6}}{j}(\xi),\nonumber
\end{align}
where clearly the degree $\pd$ must satisfy $\pd \geq 5$ (for $\pd = 5$, the sum appearing in $\tr$ is empty). The parameters 
\begin{equation} \label{eq-def-x}
{\bfm{y}} = \begin{pmatrix}
{\bfm{y}_1}\\
{\bfm{y}_2}\\
{\bfm{y}_3}
\end{pmatrix} \in \RR^{\overline{n}}, \quad
\begin{array}{lll}
{\bfm{y}_1} = \left((a_{0,\eta})_{\eta=0}^2, (a_{1,\eta})_{\eta=0}^2, (b_{0,\eta})_{\eta=0}^1,(b_{1,\eta})_{\eta=0}^1\right)^T\in \RR^{n_{1}}, & \quad n_{1} = 10,\\[1mm]
{\bfm{y}_2} = (c_{0,j})_{j=0}^{\pd-6} \in \RR^{n_{2}},  & \quad n_{2} = \pd-5,   \\[1mm] 
{\bfm{y}_3} = (c_{1,j})_{j=0}^{{\pd + 2\delta -6}} \in \RR^{n_{3}}  & \quad n_{3} = {\pd + 2\delta -5},
\end{array}
\end{equation}
where $\overline{n} = n_{1}+n_{2}+n_{3} ={2\pd+ 2\delta}$,
are linked through the linear constraints
\begin{equation} \label{eq-def-eq-Mk}
\Dti{\xi}{i}\widehat{r}^{(k)}(\xi)\big|_{\xi = 0} = 0, \; \;
i=0,1,\dots, \deg(\beta)-1, \quad 
\Dti{\xi}{i}r^{(k)}(\xi)\big|_{\xi = 0} = 0, \; \;
i=\pd+\deg(\beta) + \varsigma^{(k)}, \dots, \deg(r^{(k)}), 
\end{equation}
where 
\begin{equation} \label{def-rk-for-alg}
r^{(k)} =\alpha^{(k)} \, \w + \beta^{(k)} \, \tr' \quad \mbox{and} \quad \widehat{r}^{(k)}={\rm rem}\left(r^{(k)},\beta\right)    
\end{equation}
for $k \in \set{1, 2}$, that follow directly from \eqref{eq-cond-remAndHighOrdT}. 
Note that, while the computation of the polynomial remainder in~\eqref{def-rk-for-alg} is a symbolic operation, it can be performed in a numerically stable way if a tolerance is set for checking whether the leading term of a polynomial is zero.
The conditions above give rise to a linear system of the form 
$\bfm{M} \bfm{y} = 0$, so $\bfm{y}$ must lie in the kernel of $\bfm{M}$. 
We compute the singular value decomposition $\bfm{M} = \bfm{U}\bfm{\Sigma}\bfm{V}^T$ and determine the rank $R = \mathrm{rank}(\bfm{M})$, given by the number of nonzero singular values. In the following, we adopt Matlab notation for denoting the submatrices. The columns of the matrix $\bfm{V}_1 = \bfm{V}(\,:\,, (R+1):\overline{n})$ form a basis of the kernel of $\bfm{M}$, and its dimension $\chi:= \overline{n}-R$ reveals the number of degrees of freedom for the interface. If $\chi$ is less than $\chi_{1}:=n_{1}$, the degree $\pd$ must be increased. Otherwise, we introduce new  parameters ${\bfm{z}} \in \mathbb{R}^{\chi}$ and express
${\bfm{y}}=\bfm{V}_1 {\bfm{z}}$. To distinguish the free parameters in ${\bfm{z}}$ corresponding to the interpolation conditions from the remaining ones, we proceed as follows. We partition $\bfm{V}_1$ into three submatrices $\bfm{V}_{1,1}$, $\bfm{V}_{1,2}$, and $\bfm{V}_{1,3}$, containing the first $n_{1}$, the next $n_{2}$, and the remaining $n_{3}$ rows of $\bfm{V}_{1}$, respectively. If ${\rm rank}(\bfm{V}_{1,1})$ is less than $n_{1}$, the degree $\pd$ must be increased. Else, from the singular value decomposition $\bfm{V}_{1,1} = \widetilde{\bfm{U}}\widetilde{\bfm{\Sigma}}\widetilde{\bfm{V}}^T$, $\widetilde{\bfm{\Sigma}}=(\widetilde{\bfm{S}}, \bfm{0})$, where $\widetilde{\bfm{S}}$ is a non-singular diagonal $n_{1}\times n_{1}$ matrix, 
we express ${\bfm{z}}$ with new parameters ${\widetilde{\bfm{z}}}$ as ${\bfm{z}} = \widetilde{\bfm{V}} {\widetilde{\bfm{z}}}$. Choosing first 
$n_{1}$ components of ${\widetilde{\bfm{z}}}$ as ${\widetilde{\bfm{z}}}(1:n_{1}) = \widetilde{\bfm{S}}^{-1}\widetilde{\bfm{U}}^T {\bfm{y}_1}$ implies that $\bfm{V}_{1,1} {\bfm{z}} = {\bfm{y}_1}$, so such a choice preserves ${\bfm{y}_1}$ as parameters of freedom. 
Further, we compute the singular value decomposition of 
$\bfm{B}=\bfm{V}_{1,2} \widetilde{\bfm{V}}(\,:\, , (n_{1}+1):\chi)$, $\bfm{B}=\widehat{\bfm{U}} \widehat{\Sigma}\widehat{\bfm{V}}^T$, and select 
\begin{align*}
& {\widetilde{\bfm{z}}}(n_{1}+1:\chi) = \widehat{\bfm{V}} \bfm{w}, \quad 
 \bfm{w} = \begin{pmatrix}
\bfm{w}_{2}\\\bfm{w}_{3}
\end{pmatrix}, \quad \bfm{w}_{2} := (\zeta_j)_{j=1}^{\chi_{2}} \in \RR^{\chi_{2}}, 
\quad \bfm{w}_{3} := (\gamma_j)_{j=1}^{\chi_{3}} \in \RR^{\chi_{3}}, \\
& \chi_{2}={\rm rank} (\bfm{B}), \quad 
\chi_{3} = \chi- n_{1} - \chi_{2},
\end{align*}
so the components of $\bfm{w}_{2}$ and $\bfm{w}_{3}$ represent degrees of freedom associated with $\bfm{y}_2$ and $\bfm{y}_3$, respectively. 
To summarize, the parameters $\bfm{y}$ are expressed with new ones $
\bfm{w}_{1}:= {\bfm{y}}_1$, $\bfm{w}_{2}$ and $\bfm{w}_{3}$ as
\begin{equation} \label{eq-def-x*}
{\bfm{y}} = \bfm{V}_1 {\widetilde{\bfm{V}}} 
\begin{pmatrix}
\widetilde{\bfm{S}}^{-1}\widetilde{\bfm{U}}^T \bfm{w}_{1} \\
{\widehat{\bfm{V}}} \bfm{w}
\end{pmatrix}.  
\end{equation}
Algorithm~\ref{alg:DegreesOfFreedom} summarizes the presented construction.

\begin{algorithm}
\setstretch{1.2}
\SetKwInput{KwData}{Input}
\SetKwInput{KwResult}{Output}
\SetAlgoLined
\KwData{Gluing functions $\beta$, ${\alpha}^{(k)}$, $\beta^{(k)}$,  $k \in \set{1,2}$
}
{
1. Express $\tr$ and $\w$ via \eqref{eq-def-tr-w-1} using ${\bfm{y}}$ from \eqref{eq-def-x}, and fix the values $n_{1}$, $n_{2}$, $n_{3}$, and $\overline{n} = n_{1}+n_{2}+n_{3}$;\\
2. Define $r^{(k)}$, $\widehat{r}^{(k)}$ for $k \in \set{1,2}$ as in 
\eqref{def-rk-for-alg};\\
3. Compute $\bfm{M}$ from equations \eqref{eq-def-eq-Mk}; \\
4. Compute SVD decomposition $\bfm{M} = \bfm{U}\bfm{\Sigma}\bfm{V}^T$ and set $R={\rm rank}(\bfm{M})$, $\chi = \overline{n} - R$;\\
5. {\bf if}\, $\chi< n_{1}$\, {\bf then}\, {{\tt exit}}\\
6. Define $\bfm{V}_1 = \bfm{V}(\,:\,, (R+1):\overline{n})$, $\bfm{V}_{1,1} = \bfm{V}_{1}(1:n_{1},\,:\,)$; \\
7. Define $\bfm{V}_{1,2} = \bfm{V}_{1}((n_{1}+1):(n_{1}+n_{2}),\,:\,)$,
$\bfm{V}_{1,3} = \bfm{V}_{1}((n_{1}+n_{2}+1):\overline{n},\,:\,)$;\\
8. {\bf if}\, ${\rm rank}(\bfm{V}_{1,1})< n_{1}$\, {\bf then}\, {{\tt exit}}\\
9. Compute SVD decomposition $\bfm{V}_{1,1} = \widetilde{\bfm{U}}(\widetilde{\bfm{S}}, \bfm{0})\widetilde{\bfm{V}}^T$;\\
10. Define $\bfm{B}=\bfm{V}_{1,2} \widetilde{\bfm{V}}(\,:\, , (n_{1}+1):\chi)$ and compute SVD decomposition $\bfm{B}=\widehat{\bfm{U}} \widehat{\Sigma}\widehat{\bfm{V}}^T$;\\
11. Set $\chi_{1}=n_{1}$, $\chi_{2}={\rm rank} (\bfm{B})$, $\chi_{3} = \chi- \chi_{1} - \chi_{2}$;\\
12. Define new free parameters $\bfm{w}_{1} :={\bfm{y}}(1:\chi_{1})$, $\bfm{w}_{2} := (\zeta_j)_{j=1}^{\chi_{2}}$ and  
$\bfm{w}_{3} := (\gamma_j)_{j=1}^{\chi_{3}}$;\\ 
13. Express ${\bfm{y}}$ according to \eqref{eq-def-x*};\\
14. 
Using ${\bfm{y}}$, represent $\tr$ and $\w$ in terms of $\bfm{w}_{1}$, $\bfm{w}_{2}$, and $\bfm{w}_{3}$;\\
}
\KwResult{(trace $\tr$; directional derivative $\w$; free parameters $\bfm{w}_{1}$, $\bfm{w}_{2}$ and $\bfm{w}_{3}$; dimension $\chi$)}
\caption{{\tt Trace $\tr$ and directional derivative $\w$}}
\label{alg:DegreesOfFreedom}
\end{algorithm}

For the isogeometric functions over more than two elements, clearly some of the free parameters associated with edges do not represent the degrees of freedom for the spline as at the vertices of the partition one must also provide the $\C{1}$ continuity. To ensure that the dimension is independent of the geometry of the mesh at the
vertices, one can either require higher order of smoothness (super-smoothness) at vertices or employ 
macro-element splitting techniques. In this paper, we apply the first approach. In particular, we additionally impose the $C^2$ continuity conditions at the vertices of the partition. 
This leads to the $C^1$ super-smooth isogeometric spline space
$$\mathcal{V}^{1,2}_\pd(\Omega) := \left\{
\varphi \in \mathcal{V}^{1}_\pd(\Omega), \quad \varphi \in \C{2}(\bVm), \; \forall \,
m\in \indV
\right\},
$$
which is analyzed in detail in the next subsections.

\subsection{Bernstein representation over a mixed surface mesh} \label{subsec:BernRep}

The Bernstein representation of polynomials in the parametric domain $\parDomT$ and $\parDomQ$ allows us to define for the isogeometric spline $\varphi\in \mathcal{V}_\pd(\Omega)$ the control structure, given by domain points in ${\Omega}$ and the corresponding Bernstein coefficients. In more detail, for $k\in\indT$, i.e.,
$\Omega^{(k)} = \bfm{F}^{(k)}(\parDomT)$ and $k\in\indQ$, i.e., $\Omega^{(k)} = \bfm{F}^{(k)}(\parDomQ)$, the domain points are defined as
\begin{equation*}
\bfm{X}^{(k)}_{\mathbf i} = \bfm{F}^{(k)}\left(\tfrac{i_1}{d},\tfrac{i_2}{d}\right), \quad 
\mathbf i \in \indSetT{d}, \quad
\bfm{X}^{(k)}_{\mathbf i} = \bfm{F}^{(k)}\left(\tfrac{i_1}{d},\tfrac{i_2}{d}\right), \quad \mathbf i \in \indSetQ{d},
\end{equation*}
respectively. 
Choosing the coefficients $\mathcal{G}^{(k)} :=\set{b_{\mathbf i}^{(k)}; \mathbf i \in \indSetT{d}}$ for $k\in\indT$ and
$\mathcal{G}^{(k)} :=\set{b_{\mathbf i}^{(k)}; \mathbf i \in \indSetQ{d}}$ for $k\in\indQ$
one uniquely determines $\varphi^{(k)}$ through the  polynomials 
$$f^{(k)}(\xiBo) = \sum_{\mathbf i \in \indSetT{d}} b_{\mathbf i}^{(k)} \bernT{d}{\mathbf i}(\xiBo), \mbox{ for }k\in \indT, \quad
f^{(k)}(\xiBo) = \sum_{\mathbf i \in \indSetQ{d}} b_{\mathbf i}^{(k)} \bernB{d}{\mathbf i}(\xiBo), \mbox{ for }k\in \indQ,
$$
in the parametric space.
\begin{figure}[hbt]
\begin{center}
\begin{minipage}{0.48\textwidth}
\centering
\includegraphics[width=1\textwidth]{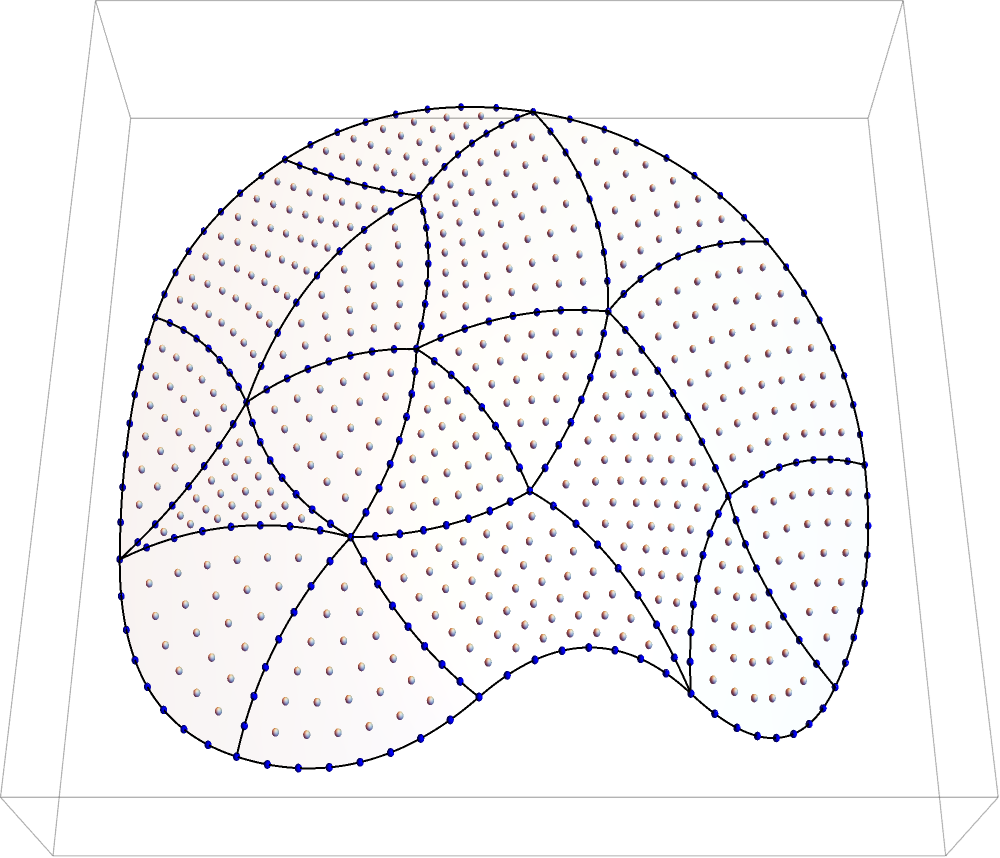}
\end{minipage}
\begin{minipage}{0.48\textwidth}
\centering
\includegraphics[width=1\textwidth]{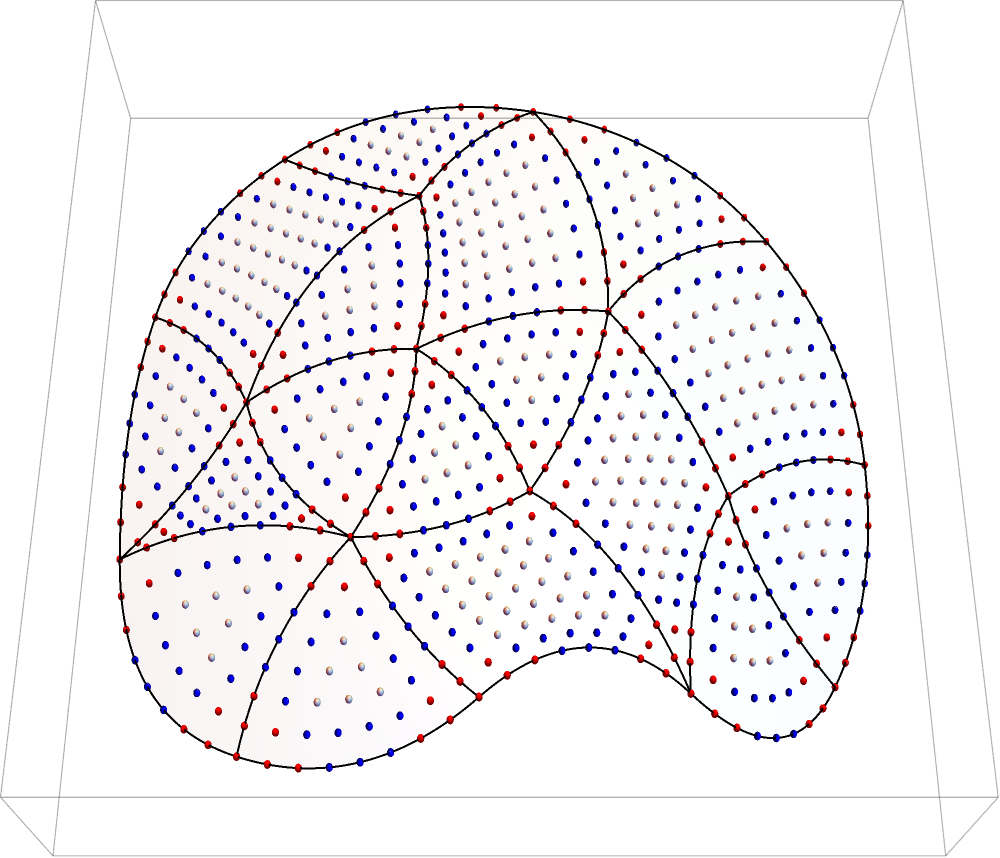}
\end{minipage}
\end{center}
\caption{Left: Domain points defining a $C^0$-continuous isogeometric function. The domain points on common element interfaces (shown in blue) are shared by the adjacent elements and therefore represent the same degrees of freedom.
Right: Partition of the domain points into vertex (red), edge (blue), and triangle/quadrilateral (gray) groups used in the construction of $C^1$-continuous isogeometric functions.}
\label{fig:DomainPoints-C1}
\end{figure}

Clearly, to achieve $\C{1}$ continuity, not all of these coefficients can be chosen freely. 
For a $\C{0}$ continuous isogeometric function $\varphi\in \mathcal{V}_d(\Omega)$, 
the coefficients associated with domain points on each edge shared by two elements must coincide (see Figure~\ref{fig:DomainPoints-C1}, left). 
Consequently,
$$\dim \left(\mathcal{V}_d(\Omega) \cap C^0(\Omega) \right) = |\indV| +(d-1) |\indE| + \binom{d-1}{2} |\indT| + (d-1)^2 |\indQ|. $$
For a $\C{1}$ continuous isogeometric function, we organize the coefficients of each element 
into four groups: {\it vertex}, {\it edge}, {\it triangle} and {\it quadrilateral} group (see Figure~\ref{fig:DomainPoints-C1}, right). In more detail, for the triangular element $\Omega^{(k)}$, 
$k \in \indT$, the coefficients are grouped into  
\begin{itemize}
\item Vertex group: $\mathcal{G}^{(k)}_{\bV} = \set{b^{(k)}_{\mathbf i}, b^{(k)}_{(d-i_1 -i_2,i_2)}, b^{(k)}_{(i_1,d-i_1-i_2)}; 0\leq i_1+i_2\leq2}$
\item Triangle group: $\mathcal{G}^{(k)}_{\triangle} = \set{b^{(k)}_{\mathbf i+\mathbf 2}; \mathbf i\in \indSetT{d-6}}$
\item Edge group: $\mathcal{G}^{(k)}_{\bE} = \mathcal{G}^{(k)} \setminus (\mathcal{G}^{k)}_{\bV}\cup \mathcal{G}^{(k)}_{\triangle}) $
\end{itemize}
and for the quadrilateral element $\Omega^{(k)}$, $k \in \indQ$, we define
\begin{itemize}
\item Vertex group: $\mathcal{G}^{(k)}_{\bV} =  \set{b^{(k)}_{\mathbf i}, b^{(k)}_{(d-i_1,i_2)}, b^{(k)}_{(d-i_1,d-i_2)}, b^{(k)}_{(i_1,d-i_2)}; 0\leq i_1+i_2\leq2}$
\item Quadrilateral group: $\mathcal{G}^{(k)}_{{\tiny \Box}} = \set{b^{(k)}_{\mathbf i+\mathbf 2}; \mathbf i\in \indSetQ{d-4}}$
\item Edge group: $\mathcal{G}^{(k)}_{\bE} = \mathcal{G}^{(k)} \setminus (\mathcal{G}^{(k)}_{\bV}\cup \mathcal{G}^{(k)}_{{\tiny \Box}}) $
\end{itemize}
Here, $\mathbf i+\mathbf 2 = (i_1+2,i_2+2)$.

\subsection{Structure of the $C^1$ super-smooth isogeometric spline space}

To construct an isogeometric function $\varphi \in \mathcal{V}_d^{1,2}(\Omega)$, the coefficients in the triangle and quadrilateral group can be chosen freely, whereas the remaining ones require special attention.
For each interface $\bE_{\ell}$, we first apply Algorithm~\ref{alg:DegreesOfFreedom}
and express the trace function $\theta_{\ell}$ and the crossing derivative function $\w_\ell$ in terms of free parameters $\bfm{w}_{\ell,1}$, $\bfm{w}_{\ell,2}$, $\bfm{w}_{\ell,3}$, whose total number is $\chi_\ell$. The parameters $\bfm{w}_{\ell,2}$ and $\bfm{w}_{\ell,3}$ can be chosen freely, while $\bfm{w}_{\ell,1}$ must be determined  from the $C^2$ interpolation data at the vertices of the edge. 
More precisely, since each coefficient in the vertex group is determined by the trace and the crossing derivative functions of two different edges (that share the vertex), consistency must be ensured. This can be achieved by prescribing at every vertex $\bV_m$ the $\C{2}$ interpolation conditions, specified in \eqref{int-data-at-V} with parameters 
\begin{equation} \label{eq-int-data-for-V}
\bfm{c}_m := \left(c_{m,0,0}, c_{m,1,0}, c_{m,0,1}, c_{m,2,0}, c_{m,1,1}, c_{m,0,2} \right), 
\quad m\in \indV,   
\end{equation}
that $\varphi$ must fulfill. 
As shown in Theorem~\ref{thm-3}, these interpolation conditions uniquely determine, through \eqref{eq-int-cond-Vm-0}, \eqref{eq-int-cond-Vm-1} and \eqref{eq-int-cond-Vm-2}, for each $k$, the values as well as the first and second-order derivatives of $f^{(k)}$ at all boundary parameters.  
Furthermore, from \eqref{eq:def-theta-l}--\eqref{eq-vartheta} one can see that these quantities determine for every edge $\bE_\ell$ the following values: $\tr_\ell(\xi), \tr_\ell'(\xi), \tr_\ell''(\xi)$ and $\w_\ell(\xi), \w_\ell'(\xi)$ for $\xi=0$ and $\xi=1$. These values coincide precisely with the components of the vector $\bfm{w}_{\ell,1}$. 

For a boundary edge $\bE_{\ell} \subset \Omega^{(k)}$, all coefficients in the edge group associated with $\bE_{\ell}$ can be chosen freely. To unify the notation and the counting of degrees of freedom, we denote the free coefficients within this edge group corresponding to the domain points lying on the edge by $\bfm{w}_{\ell,2}$ and those corresponding to the domain points in the first row adjacent to the edge by $\bfm{w}_{\ell,3}$. Accordingly, we set
$\chi_{\ell,2}=d-5$, $\chi_{\ell,3}=d-4+\varsigma^{(k)}$, where $\varsigma^{(k)}$ is defined in \eqref{def-varsigma}.

To summarize, determining $\bfm{w}_{\ell,1}$ from the interpolation data \eqref{eq-int-data-for-V}, together with a free choice of parameters in $\bfm{w}_{\ell,2}$ and $\bfm{w}_{\ell,3}$ (comprising $\chi_{\ell,2} + \chi_{\ell,3} = \chi_\ell - 10$ parameters), uniquely determine $\tr_\ell$ and $\w_\ell$, and consequently the coefficients in the vertex and edge groups. Thus, the dimension of 
$\mathcal{V}^{1,2}_\pd(\Omega)$ equals 
\begin{equation} \label{eq-dimension}
\dim \mathcal{V}^{1,2}_\pd(\Omega)  = 6|\indV| +
\sum_{\ell\in \indE} (\chi_\ell-10) +  \binom{d-4}{2} |\indT| + (d-3)^2 |\indQ|, 
\end{equation}
and, based on the above results, it is straightforward to determine a basis of the space $\mathcal{V}^{1,2}_\pd(\Omega)$
as explained in the following subsection.

\subsection{Basis for the $\C{1}$ super-smooth isogeometric spline space}\label{subsec:basis}

In accordance with the grouping of coefficients, we introduce four types of basis functions.
\begin{itemize}
\item
{Vertex basis functions}: $\Bv{m,i}$, $i = 1, 2, \dots, 6$, for every vertex $\bV_m$, $m \in \indV$;
\item
{Edge basis functions}:  $\Betr{\ell, i}$, $i = 1, 2, \dots, \chi_{\ell,2}$, 
and $\Bew{\ell, i}$, $i = 1, 2, \dots, \chi_{\ell,3}$, for every edge $\bE_{\ell}$, $\ell \in \indE$;
\item
{Triangle basis functions}:
$\Bt{k,\mathbf i}$, $\mathbf i\in \indSetT{d-6}$, for every $\Omega^{(k)}$, $k \in \indT$;
\item
{Quadrilateral basis functions}:
$\Bq{k,\mathbf i}$, $\mathbf i\in \indSetQ{d-4}$, for every $\Omega^{(k)}$, $k \in \indQ$.
\end{itemize}
To define these basis functions, we need to prescribe the interpolation values $\bfm{c}_m$, $\forall m\in \indV$,
values of free parameters $\bfm{w}_{\ell,2}$, $\bfm{w}_{\ell,3}$, $\forall \ell \in \indE$, and values of coefficients in the triangle group $\mathcal{G}^{(k)}_{\triangle}$, $\forall k \in \indT$, and quadrilateral group $\mathcal{G}^{(k)}_{{\tiny \Box}}$, $\forall k \in \indQ$. To shorten the notation, let  
$$ 
\mathcal{G}^{0}_{\triangle} := \set{0, \dots, 0}, \quad 
\mathcal{G}^{0}_{{\tiny \Box}} := \set{0, \dots, 0},
$$
and let $\bfm{0}$ and $\bfm{e}_r$ denote the zero and $r$-th unit vectors of appropriate dimension.

\begin{itemize}
\item
For the vertex $\bV_{m^*}$ and $i^* \in \set{1,2,\dots, 6}$, the function $\Bv{m^*,i^*}$ is defined by choosing
$$
\bfm{c}_{m} = \begin{cases}
\bfm{e}_{i^*}, & m = m^*\\
\bfm{0}, &  m\not = m^*
\end{cases}, \quad \bfm{w}_{\ell,2}=\bfm{0}, \quad \bfm{w}_{\ell,3}= \bfm{0}, \quad
\mathcal{G}^{(k)}_{\triangle}=\mathcal{G}^{0}_{\triangle}, \quad \mathcal{G}^{(k)}_{{\tiny \Box}} =\mathcal{G}^{0}_{{\tiny \Box}}. 
$$
\item
For the edge $\bE_{\ell^*}$, we define $\Betr{\ell^*,i^*}$ for $i^* \in \set{1, 2 \dots, \chi_{\ell^*,2}}$ by choosing
$$
\bfm{c}_{m} = \bfm{0}, \quad 
\bfm{w}_{\ell,2} = 
\begin{cases}
\bfm{e}_{i^*}, & \ell = \ell^*\\
\bfm{0}, &  \ell \not = \ell^*
\end{cases}, \quad \quad \bfm{w}_{\ell,3}= \bfm{0}, \quad
\mathcal{G}^{(k)}_{\triangle}=\mathcal{G}^{0}_{\triangle}, \quad \mathcal{G}^{(k)}_{{\tiny \Box}} =\mathcal{G}^{0}_{{\tiny \Box}}, 
$$
and $\Bew{\ell^*,i^*}$ for $i^* \in \set{1, 2, \dots, \chi_{\ell^*,3}}$ by 
$$
\bfm{c}_{m} = \bfm{0}, \quad 
\bfm{w}_{\ell,2} = \bfm{0}, \quad  
\bfm{w}_{\ell,3} = \begin{cases}
\bfm{e}_{i^*}, & \ell = \ell^*\\
\bfm{0}, &  \ell\not = \ell^*
\end{cases}, \quad
\mathcal{G}^{(k)}_{\triangle}=\mathcal{G}^{0}_{\triangle}, \quad \mathcal{G}^{(k)}_{{\tiny \Box}} =\mathcal{G}^{0}_{{\tiny \Box}}. 
$$
\item
The triangle basis function $\Bt{k^*,{\mathbf i}^*}$ results from    
$\bfm{c}_{m} = \bfm{0}$, $\bfm{w}_{\ell,2} = \bfm{0}$, $\bfm{w}_{\ell,3} = \bfm{0}$,
$\mathcal{G}^{(k)}_{{\tiny \Box}} =\mathcal{G}^{0}_{{\tiny \Box}}$ and
$$
\mathcal{G}^{(k)}_{\triangle} = \set{b^{(k)}_{\mathbf i+\mathbf 2}; \mathbf i\in \indSetT{d-6}} \quad {\rm with}\quad
b^{(k)}_{\mathbf i+\mathbf 2} = \begin{cases}
1, & k=k^*, \,\mathbf i={\mathbf i}^*\\
0, & {\rm otherwise}
\end{cases}.
$$
\item
The quadrilateral basis function $\Bq{k^*,{\mathbf i}^*}$ results from    
$\bfm{c}_{m} = \bfm{0}$, $\bfm{w}_{\ell,2} = \bfm{0}$, $\bfm{w}_{\ell,3} = \bfm{0}$,
$\mathcal{G}^{(k)}_{\triangle} =\mathcal{G}^{0}_{\triangle}$ and
$$
\mathcal{G}^{(k)}_{{\tiny \Box}} := \set{b^{(k)}_{\mathbf i+\mathbf 2}; \mathbf i\in \indSetQ{d-4}} \quad {\rm with}\quad
b^{(k)}_{\mathbf i+\mathbf 2} = \begin{cases}
1, & k=k^*, \,\mathbf i={\mathbf i}^*\\
0, & {\rm otherwise}
\end{cases}.
$$
\end{itemize}

After constructing the basis functions as described above, we normalize each of them with respect to its B\'ezier coefficients. More precisely, for each basis function, we determine the largest absolute value among all its B\'ezier coefficients and divide the function by this value. As a result, the largest absolute B\'ezier coefficient of every normalized basis function is equal to one.
This normalization is motivated by the fact that the magnitudes of the B\'ezier coefficients can vary substantially between different basis functions. Numerical examples show that, without normalization, these values may differ by several orders of magnitude. The normalization therefore provides a more uniform scaling of the basis functions.

\subsection{Interpolation degrees of freedom for the $C^1$ super-smooth isogeometric spline space}\label{subsec:interpolation}

In this subsection, we establish a connection between the isogeometric spline space $\mathcal{V}^{1,2}_\pd(\Omega)$ and a well-posed interpolation problem. For this purpose, we introduce the following linear functionals, where $\epsilon_{\bV}$ denotes the linear functional that assigns to each (tri-variate) function its value at $\bV$.

\begin{subequations} \label{subeq-interpolFunctionals}
\begin{itemize}
\item
For each vertex $\bV_m$, $m\in \indV$, we define
\begin{equation} \label{subeq-interpolFunctionals-1}
\lfunvert{m,\eta_1, \eta_2} := \lim_{\Omega\ni\xBo \rightarrow \bVm} \epsilon_{\xBo}\left(D^{\eta_2}_{\bfm{E}_{m,2}} D^{\eta_1}_{\bfm{E}_{m,1}}\right), 
\quad 0\leq \eta_1+\eta_2\leq 2, 
\end{equation}
where $\bfm{E}_{m,1}$, $\bfm{E}_{m,2}$ are two orthonormal vectors that span the tangent plane of $\Omega$ at $\bV_m$. 
\item 
For each interface $\bE_\ell$, $\ell\in\indE$, we select two sets of distinct parameter values,
$$
\parV^{0}_{\ell,j}\in(0,1), \quad j=1,\dots,\chi_{\ell,2}, \quad {\rm and} \quad 
\parV^{1}_{\ell,j}\in(0,1), \quad j=1,\dots,\chi_{\ell,3},
$$
and compute points 
$$\bfm{R}^{0}_{\ell,j} = \bfm{F}^{(k_1)} \bigl(q^{(k_1)}_\ell(\parV^{0}_{\ell,j})\bigr),
\quad j=1,\dots,\chi_{\ell,2}, \quad {\rm and} \quad 
\bfm{R}^{1}_{\ell,j} = \bfm{F}^{(k_1)}\bigl(q^{(k_1)}_\ell(\parV^{1}_{\ell,j})\bigr),
\quad j=1,\dots,\chi_{\ell,3},
$$
where $k_1\in \indOmega$ is such that $\bE_\ell \subset \Omega^{(k_1)}$. Then, we define
\begin{equation} \label{subeq-interpolFunctionals-2}
\lfunedgetr{\ell,j}:= \epsilon_{\bfm{R}^{0}_{\ell,j}}, \quad j=1,\dots, \chi_{\ell,2}, \qquad 
\lfunedgew{\ell,j}:= \epsilon_{\bfm{R}^{1}_{\ell,j}}
\left({\Dv_{\bfm{d}_\ell(\parV^{1}_{\ell,j})}} \right), \quad j=1,\dots, \chi_{\ell,3},
\end{equation}
{where $\bfm{d}_\ell$ is given by~\eqref{eq:def-direction-d-ell}}.
\item
For each element $\Omega^{(k)}$, we define
\begin{equation} \label{subeq-interpolFunctionals-3}
\lfunTri{k,\mathbf i} = \epsilon_{\bfm{X}^{(k)}_{\mathbf i+\mathbf 2}}, \quad  \mathbf i\in \indSetT{\pd-6}, \; k \in \indT, \qquad 
\lfunQ{k,\mathbf i} = \epsilon_{\bfm{X}^{(k)}_{\mathbf i+\mathbf 2}}, \quad  \mathbf i\in \indSetQ{\pd-4}, \; k \in \indQ.
\end{equation}
\end{itemize}
\end{subequations}
Interpolation functionals \eqref{subeq-interpolFunctionals} form a nodal minimal determining set for the isogeometric spline space $\mathcal{V}^{1,2}_\pd(\Omega)$. 
Given the data values, the interpolation problem is to find the isogeometric function $\varphi \in \mathcal{V}^{1,2}_\pd({\Omega})$ satisfying the conditions
\begin{align}
& \lfunvert{m,\mathbf i} \varphi = v_{m,\mathbf i}, \quad  0\leq i_1+i_2\leq 2, \quad m \in \indV, \nonumber\\
& \lfunedgetr{\ell,j} \varphi = v^{e,0}_{\ell,j},  \quad j=1,\dots, \chi_{\ell,2}, \quad 
\lfunedgew{\ell,j} \varphi = v^{e,1}_{\ell,j}, \quad j=1,\dots, \chi_{\ell,3}, \quad \ell \in \indE, \\ 
& \lfunTri{k,\mathbf i} \varphi = v^{t}_{k,\mathbf i}, \quad  \mathbf i\in \indSetT{\pd-6}, \; k \in \indT, \quad  
\lfunQ{k,\mathbf i} \varphi = v^{q}_{k,\mathbf i}, \quad  \mathbf i\in \indSetQ{\pd-4}, \; k \in \indQ. \nonumber
\end{align}
The problem admits a unique solution and can be used as an alternative approach for constructing basis functions, where each basis function is associated with exactly one interpolation condition equal to $1$, while all others are equal to $0$.

\section{Examples}
\label{sec:Examples}

\begin{figure}[hbt]
\begin{center}
\begin{minipage}{0.32\textwidth}
\centering
\includegraphics[width=1\textwidth]{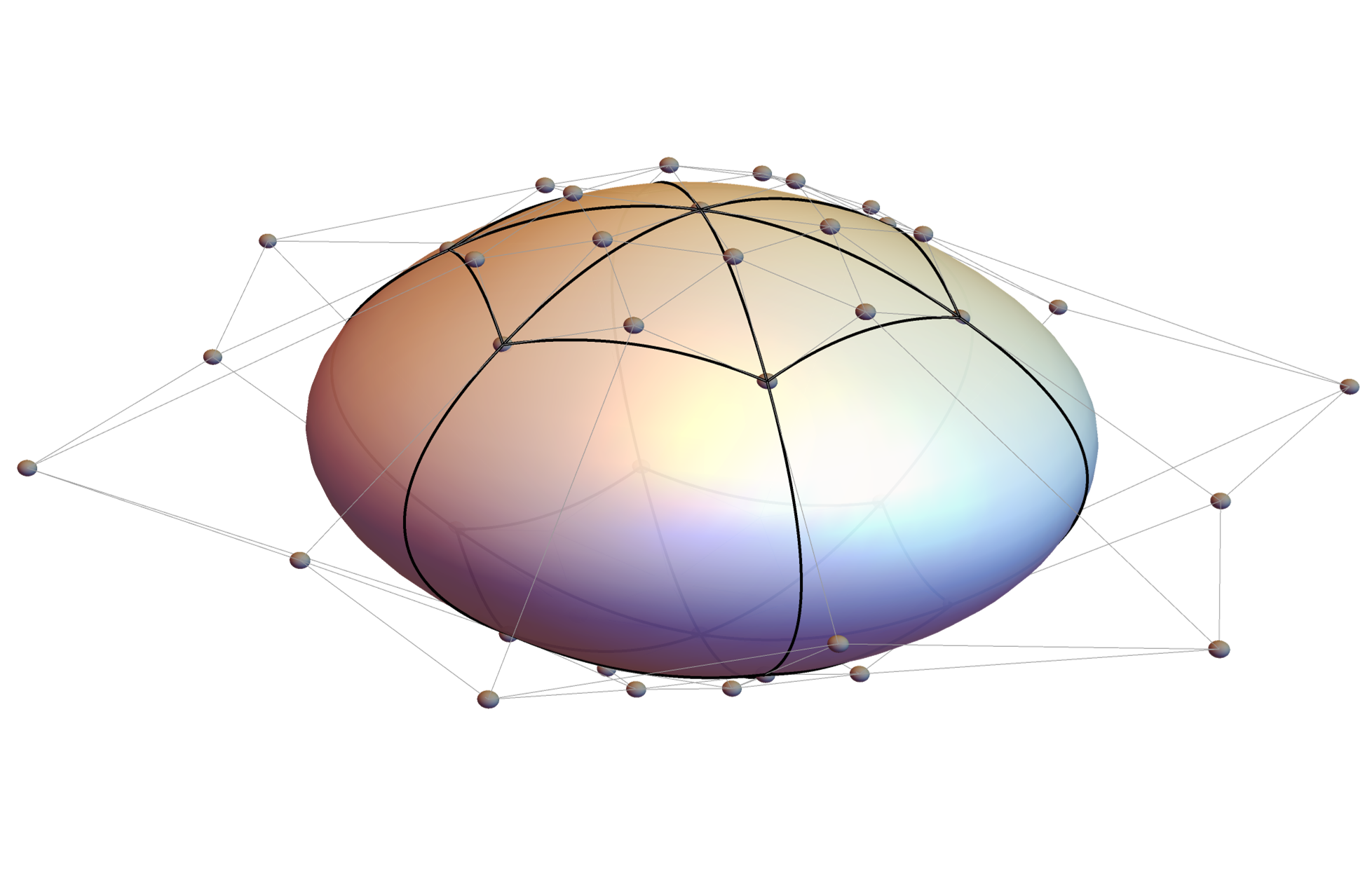}
\end{minipage}
\begin{minipage}{0.32\textwidth}
\centering
\includegraphics[width=1\textwidth]{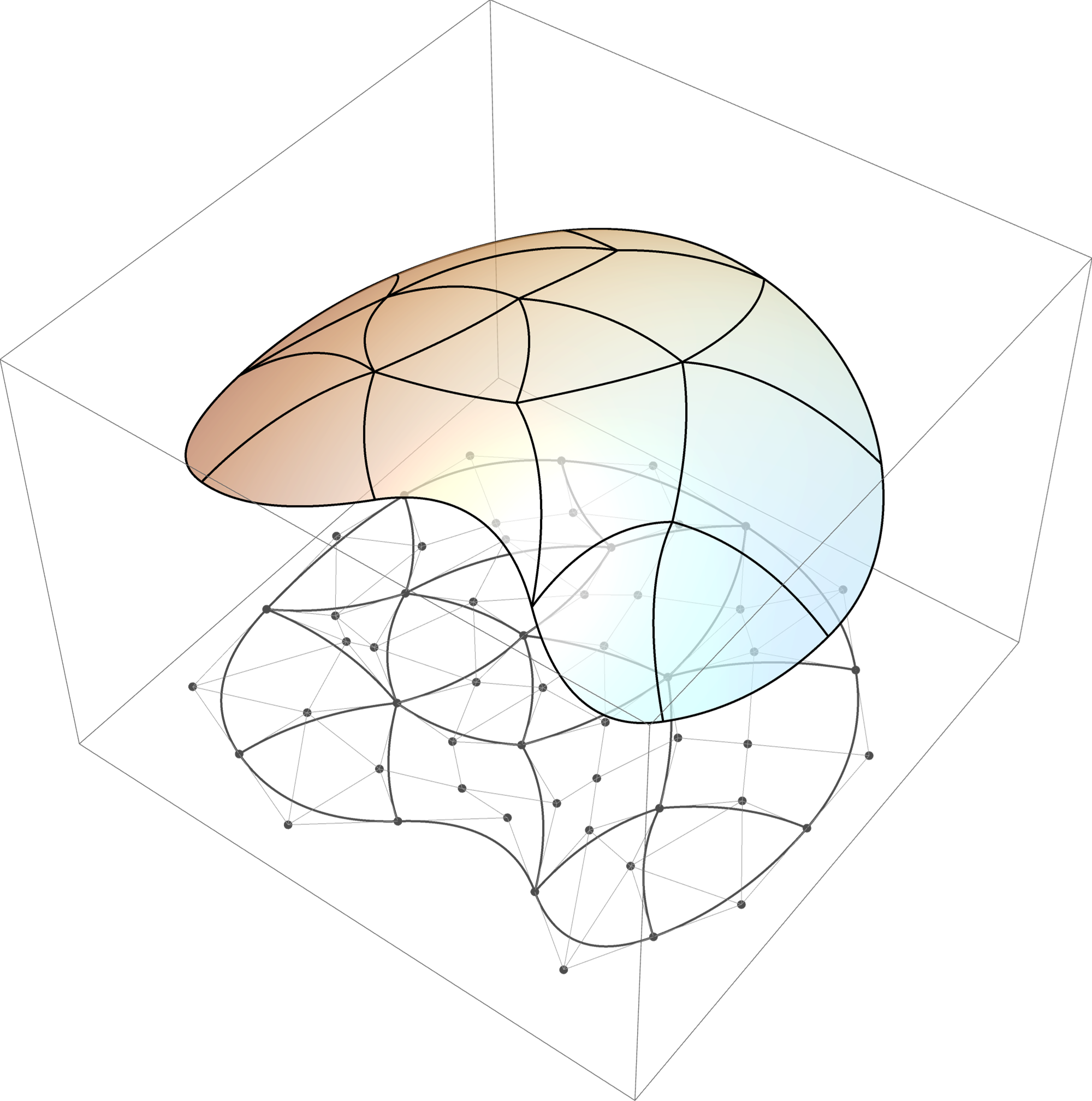}
\end{minipage}
\begin{minipage}{0.32\textwidth}
\centering
\includegraphics[width=1\textwidth]{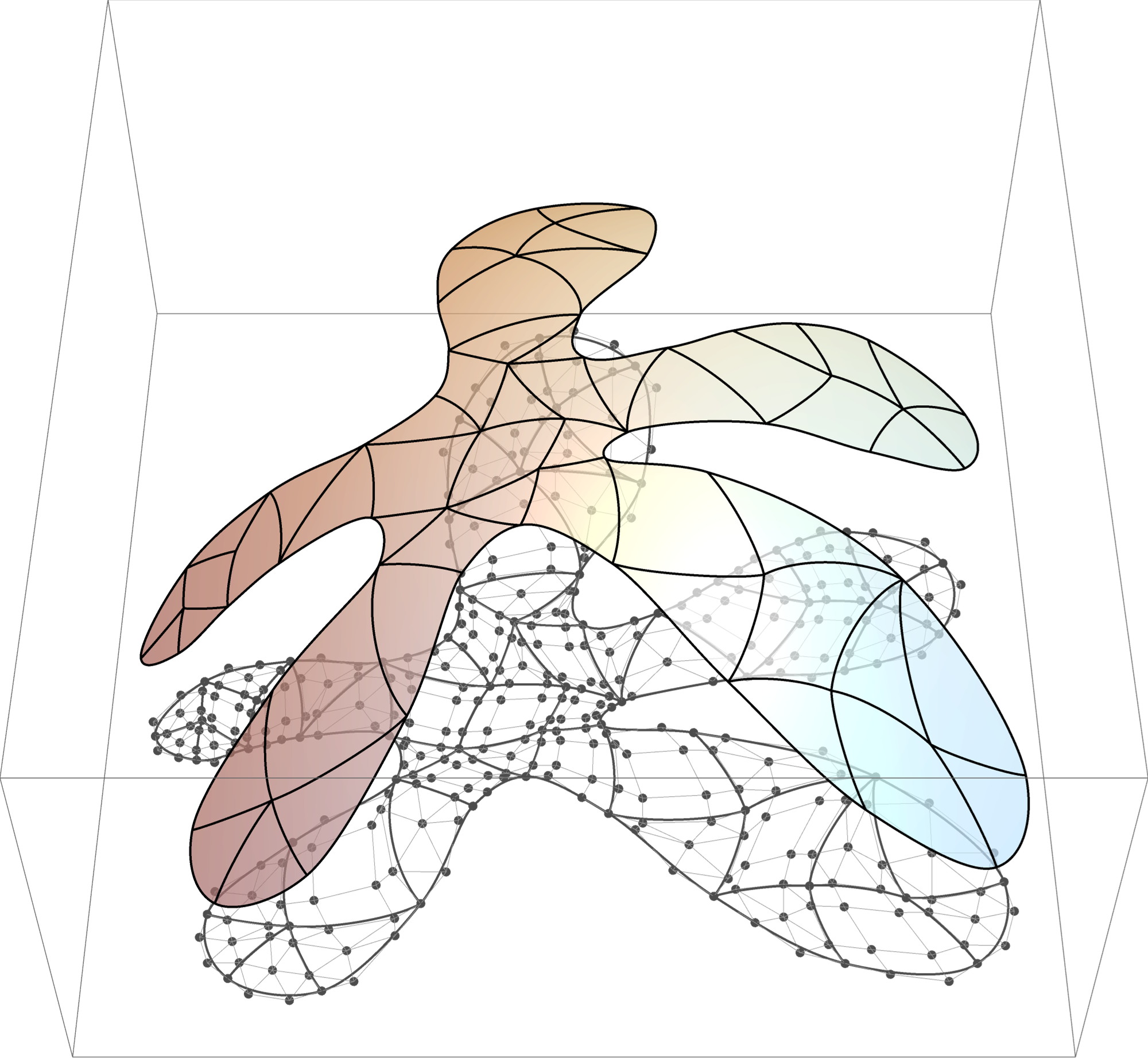}
\end{minipage}\\
\begin{minipage}{0.32\textwidth}
\centering
Mesh 1
\end{minipage}
\begin{minipage}{0.32\textwidth}
\centering
Mesh 2
\end{minipage}
\begin{minipage}{0.32\textwidth}
\centering
Mesh 3
\end{minipage}
\end{center}
\caption{
Control meshes of the geometry mappings associated with the mixed surface meshes shown in Figure~\ref{fig:MixedSurfaceMesh}.
For better visualization the planar geometry mappings (for Mesh~2 and Mesh~3) are slightly moved down. 
}
\label{fig:MixedSurfaceMeshKP}
\end{figure}

In this section, we illustrate the proposed spline construction and its properties through several numerical examples involving basis construction, interpolation, and least-squares approximation. We consider the three surface meshes shown in Figure~\ref{fig:MixedSurfaceMeshKP},
together with the control meshes of the corresponding geometry mappings. The basic characteristics of the meshes are summarized below:
$$
\renewcommand{\arraystretch}{1.1}
\begin{array}{|c|c|c|c|c|c|}
\hline
 & |\indV| & |\indE| & |\indT| & |\indQ| & \text{bounding box of control points} \\
\hline
{\tt Mesh}\,1 & 14 & 30 & 12 & 6 & [-2.1,2.1]\times [-1.8,1.8] \times [-0.9,0.9] \\
{\tt Mesh}\,2 & 17 & 33  & 12 & 5 & [1.1,9.4]\times [1.3,8.4]\times [-1.2,0.1] \\
{\tt Mesh}\,3 & 57 & 108 & 34 & 18 & [0.3,8.9]\times [0.2,9.1] \times [-1.4,0.1] \\
\hline
\end{array}
$$

Mesh~1 is represented directly by quadratic surface geometry mappings. The geometry is taken from~\cite{Cesek26}. Mesh~2 and Mesh~3 are constructed from planar domains represented by quadratic and cubic geometry mappings, respectively. The corresponding surface domains are obtained by mapping the planar domains onto the graph of a quadratic polynomial $p$ via $\Psi(x,y)=(x,y,p(x,y))^T$.
Thus, for Mesh~2 and Mesh~3, the surface geometry mappings are given by
$\bfm F^{(k)}=\Psi\circ\widehat{\bfm F}^{(k)}$,
where $\widehat{\bfm F}^{(k)}$ denotes the corresponding planar geometry mapping. 
Although the resulting surface mappings are of higher polynomial degree, their gluing data are closely related to those of the underlying quadratic and cubic planar geometries, respectively. For each edge, the gluing data determine a local reparameterization, such that the reparameterized geometry mappings join $C^1$. While the planar and corresponding surface domains possess different gluing data, the resulting local reparameterizations are the same and thus the gluing conditions yield the same numbers of free edge parameters in the spline construction.
Figure~\ref{fig:MixedSurfaceMeshKP} shows the surface meshes together with the control meshes of the quadratic  geometry mappings for Mesh~1 and of the underlying planar geometry mappings for Mesh~2 and Mesh~3. The explicit expressions of the geometry will be provided upon request and will be published in an open repository together with the final version of the paper.

\subsection{Basis functions}
\label{Examples-basis}

We first examine the dimensions and illustrate the construction of the basis functions for spline spaces of different degrees. For each of the meshes shown in Figure~\ref{fig:MixedSurfaceMeshKP}, we use Algorithm~\ref{alg:DegreesOfFreedom} to determine, for a chosen degree $d$, the number of free parameters associated with each edge.

For Mesh~1,  which is represented directly by quadratic surface geometry mappings, and Mesh~2, whose gluing data are inherited from quadratic planar geometry mappings, the algorithm shows that the lowest admissible spline degree is $d=6$. For Mesh~3, whose gluing data are inherited from cubic planar geometry mappings, the lowest admissible spline degree is $d=7$. 
The distinct pairs $(\chi_{\ell,2},\chi_{\ell,3})$ obtained for the three meshes, together with the numbers of edges for which they occur, are summarized in Table~\ref{tab:edge-dof}.

\begin {table}[h!]
\begin{align*}
\renewcommand{\arraystretch}{1.1}
\begin{array}{|c|c|c|c|}
\hline
{\rm Mesh} & d & (\chi_{\ell,2},\chi_{\ell,3})
& {\rm Number\ of\ edges} \\
\hline
{\tt Mesh}\,1 & d\geq 6 & (d-5,d-5) & 12 \\
               &          & (d-5,d-4) & 12 \\
               &          & (d-5,d-3) & 6  \\
\hline
{\tt Mesh}\,2 & d\geq 6 & (d-5,d-6) & 23 \\
               &          & (d-5,d-4) & 7  \\
               &          & (d-5,d-3) & 3  \\
\hline
{\tt Mesh}\,3 & 7 & (1,0) & 66 \\
               &  & (2,3) & 29 \\
               &  & (2,4) & 13 \\
\cline{2-4}
               & d\geq 8 & (d-5,d-8) & 66 \\
               &          & (d-5,d-4) & 29 \\
               &          & (d-5,d-3) & 13 \\
\hline
\end{array}
\end{align*}
\caption{Distinct pairs $(\chi_{\ell,2},\chi_{\ell,3})$ and their
multiplicities for the three test meshes.}
\label{tab:edge-dof}
\end{table}

Regarding the values reported in Table~\ref{tab:edge-dof}, note that, independently of the particular mesh and geometry mappings, the expressions
for the trace and the crossing derivative in \eqref{eq-def-tr-w-1} imply that increasing the spline degree $d$ by one increases the total number
$\chi_{\ell,2}+\chi_{\ell,3}$ of free edge parameters by two. In particular, for Mesh~1 and Mesh~2, increasing the spline degree by one increases both $\chi_{\ell,2}$ and $\chi_{\ell,3}$ by one. 
For Mesh~3, on the other hand, the behavior at the lowest admissible degree $d=7$ differs from that at higher degrees. Namely, in Algorithm~\ref{alg:DegreesOfFreedom}, the free parameters are distributed between the trace and the crossing derivative $\omega$ so that as many as
possible are associated with the trace, while the remaining ones are associated with $\omega$. Subject to the constraint
$\chi_{\ell,2}\leq d-5$, this leads, for some edges, to a different distribution at $d=7$ than at higher degrees.
This explains why the case $d=7$ is listed separately in Table~\ref{tab:edge-dof}.

Note also that, since the total number of free parameters associated with an edge is $10+\chi_{\ell,2}+\chi_{\ell,3}$, the values reported in
Table~\ref{tab:edge-dof} for Mesh~1 and Mesh~2 are consistent with the theoretical lower bound $2d-1$ for quadratic geometry mappings. More
precisely, this bound is attained on every interior edge of Mesh~2, whereas every edge of Mesh~1 has at least $2d$ free parameters. For Mesh~3, the
minimum number of free parameters associated with an edge is $2d-3$. 

\begin{figure}[h!]
\begin{center}
\begin{minipage}{0.3\textwidth}
\centering
\includegraphics[width=1\textwidth]{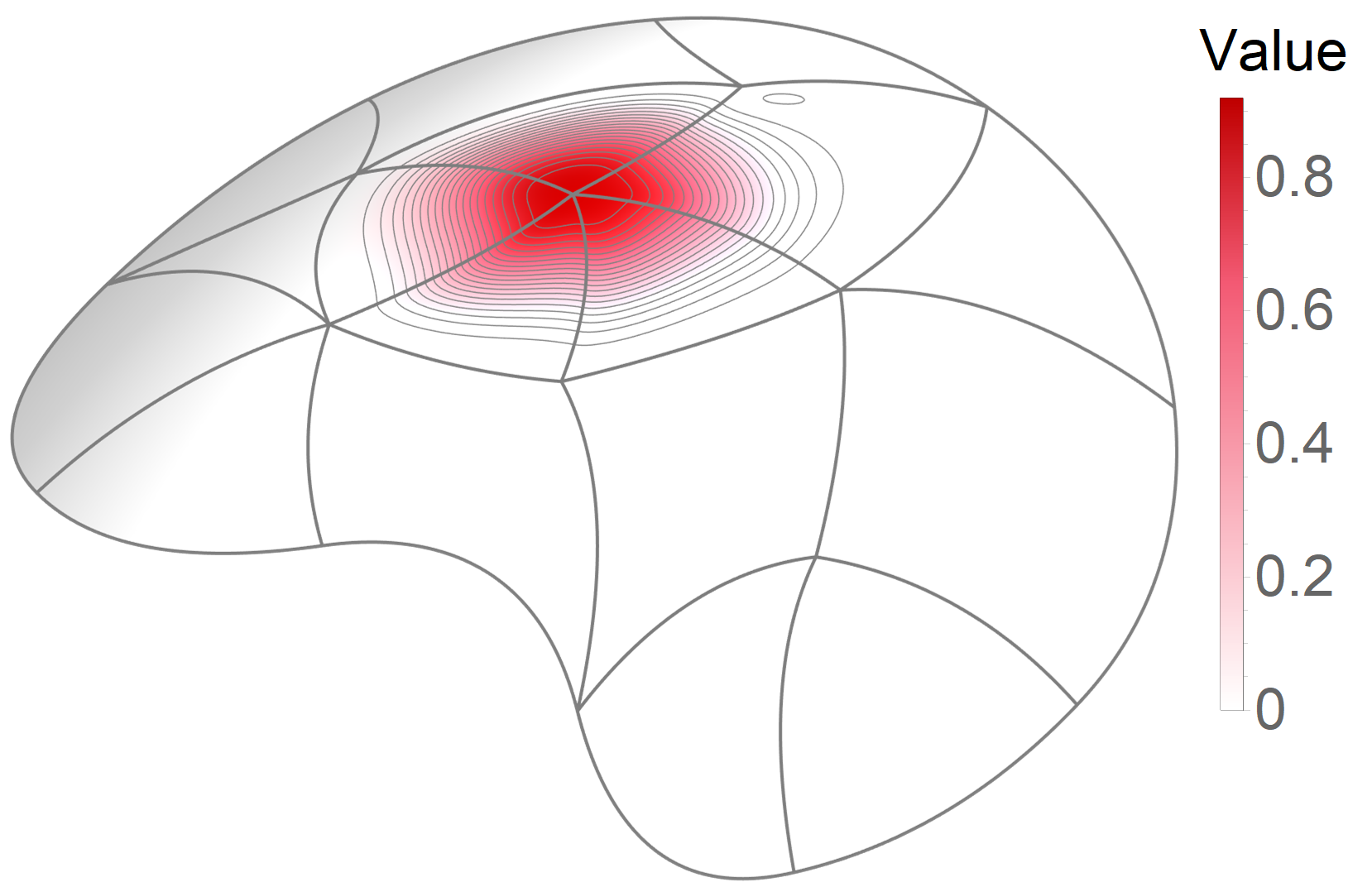}
\end{minipage}
\begin{minipage}{0.3\textwidth}
\centering
\includegraphics[width=1\textwidth]{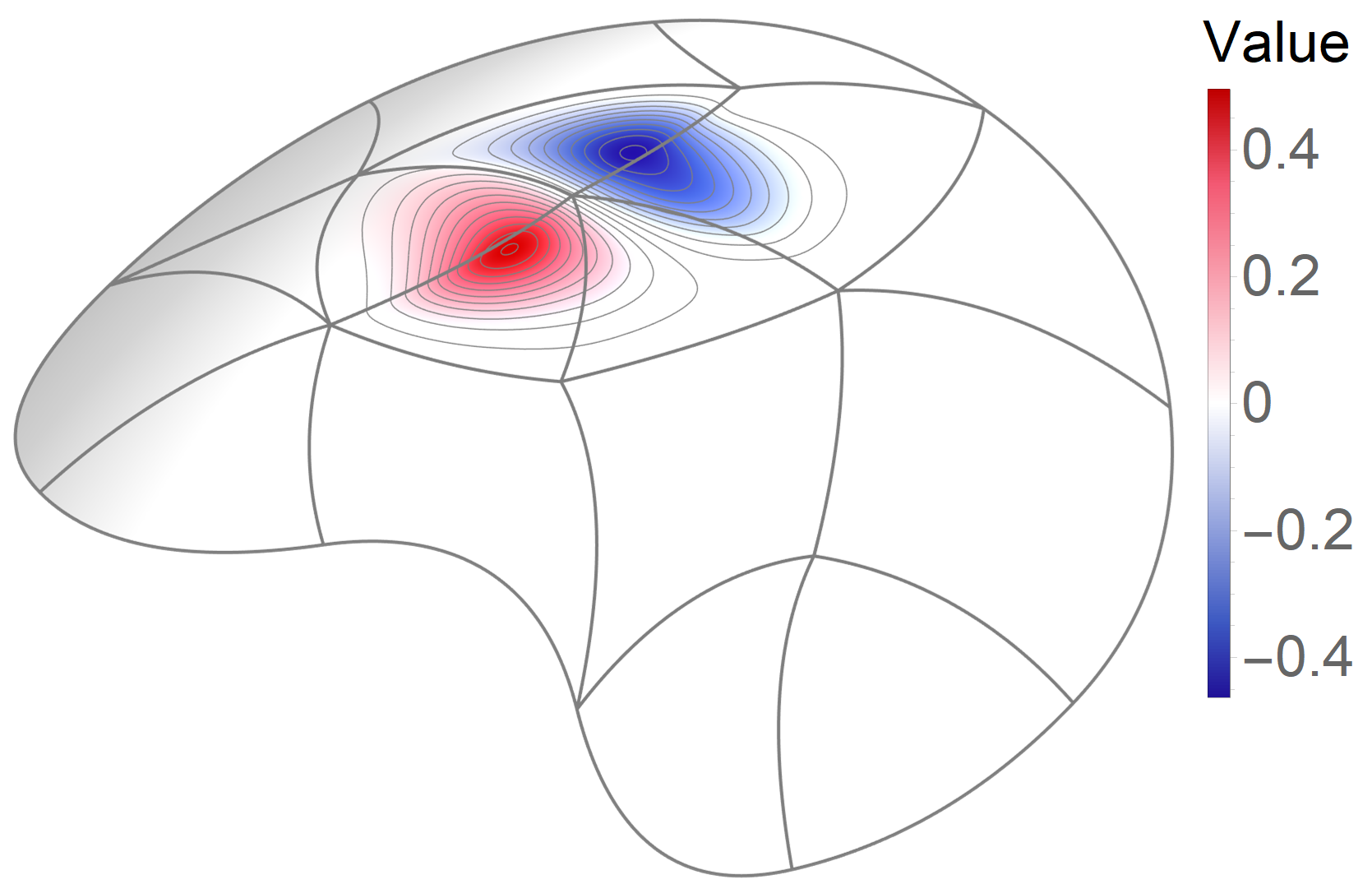}
\end{minipage}
\begin{minipage}{0.3\textwidth}
\centering
\includegraphics[width=1\textwidth]{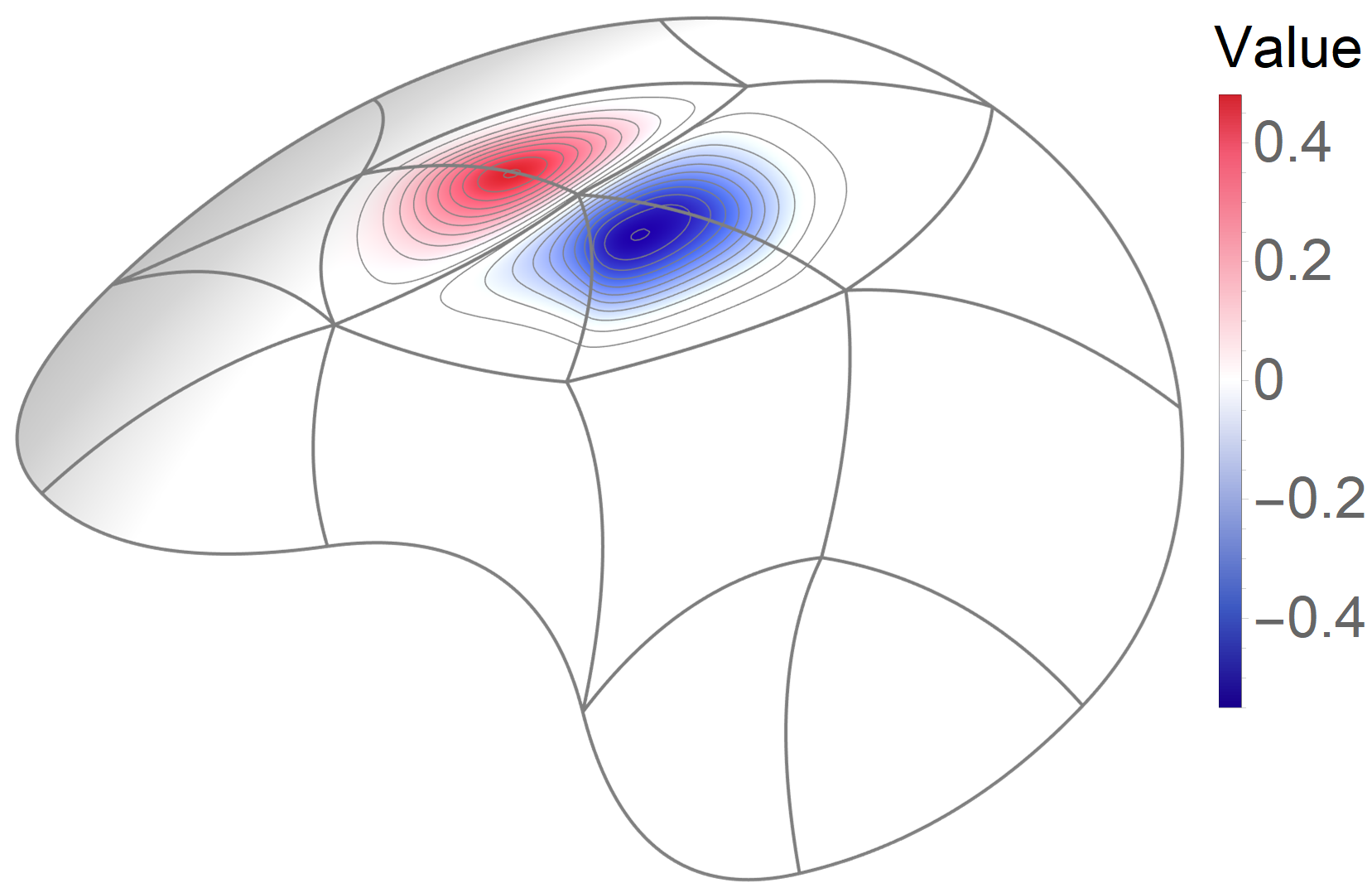}
\end{minipage}\\
\begin{minipage}{0.3\textwidth}
\centering
\includegraphics[width=1\textwidth]{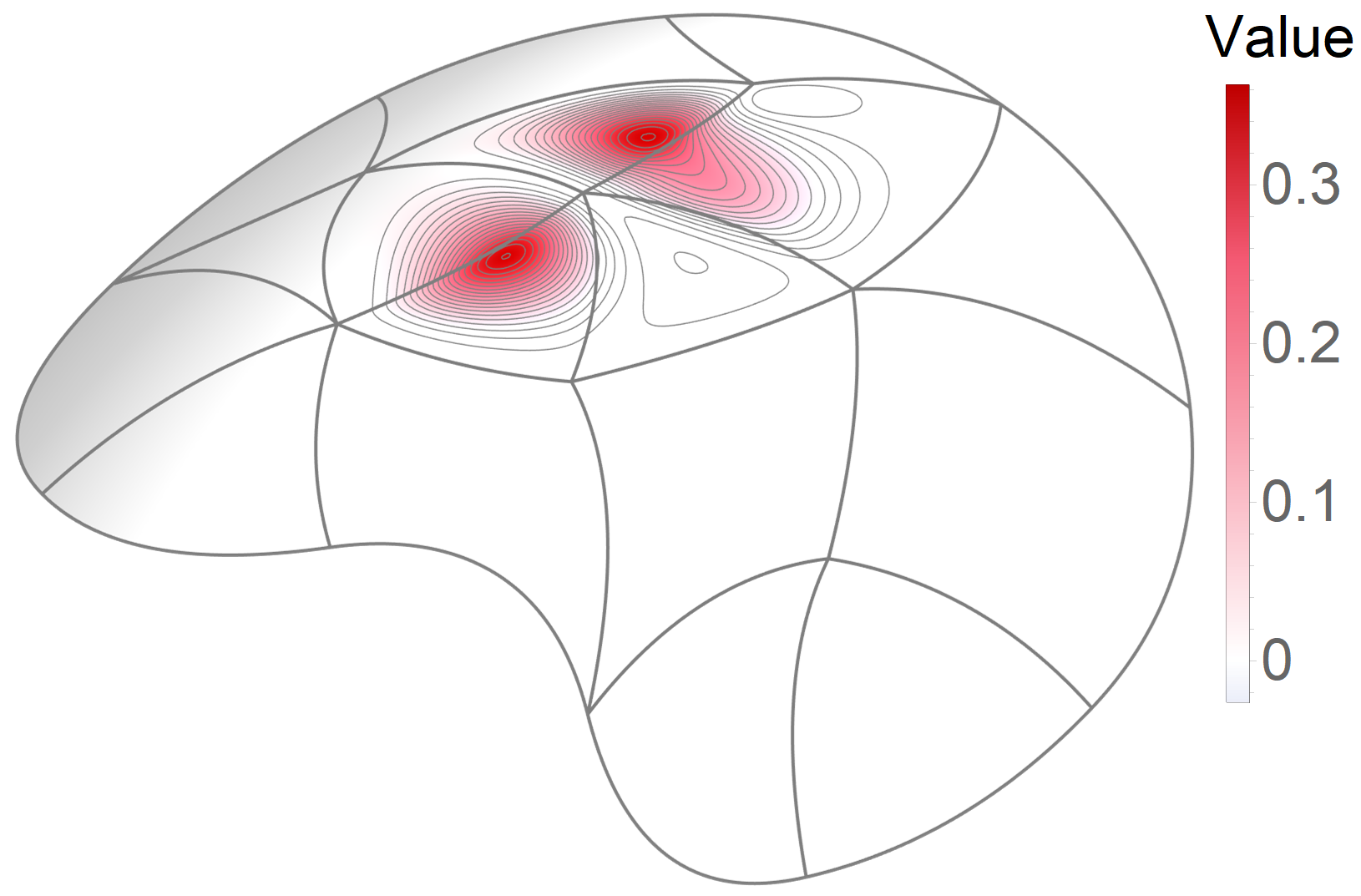}
\end{minipage}
\begin{minipage}{0.3\textwidth}
\centering
\includegraphics[width=1\textwidth]{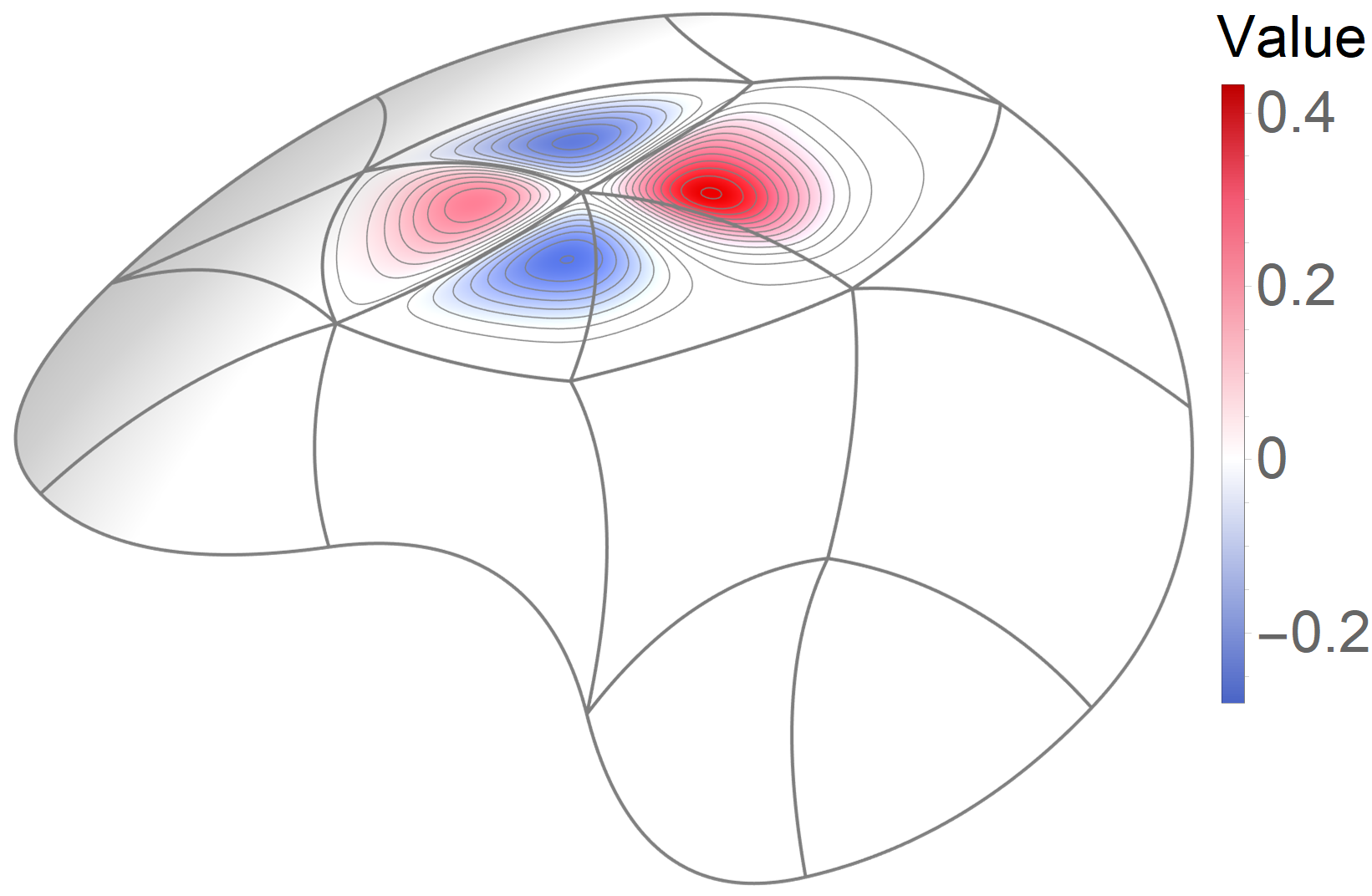}
\end{minipage}
\begin{minipage}{0.3\textwidth}
\centering
\includegraphics[width=1\textwidth]{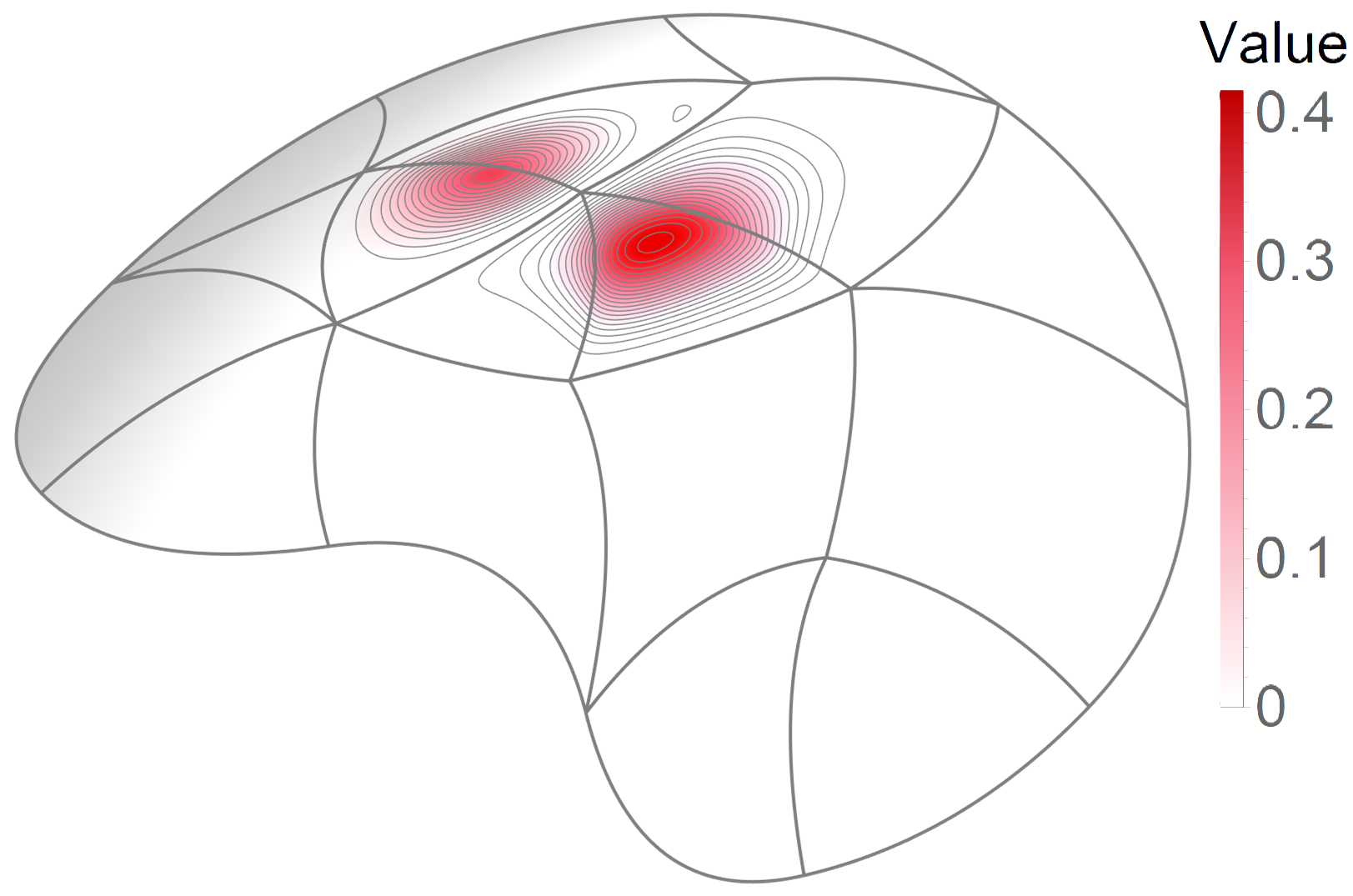}
\end{minipage}
\end{center}
\caption{Selected vertex basis functions of degree $6$ constructed over the domain partitioned by {\tt Mesh}~2.}
\label{fig:BasisFunctionsV}
\end{figure}

\begin{figure}[h!]
\begin{center}
\begin{minipage}{0.3\textwidth}
\centering
\includegraphics[width=1\textwidth]{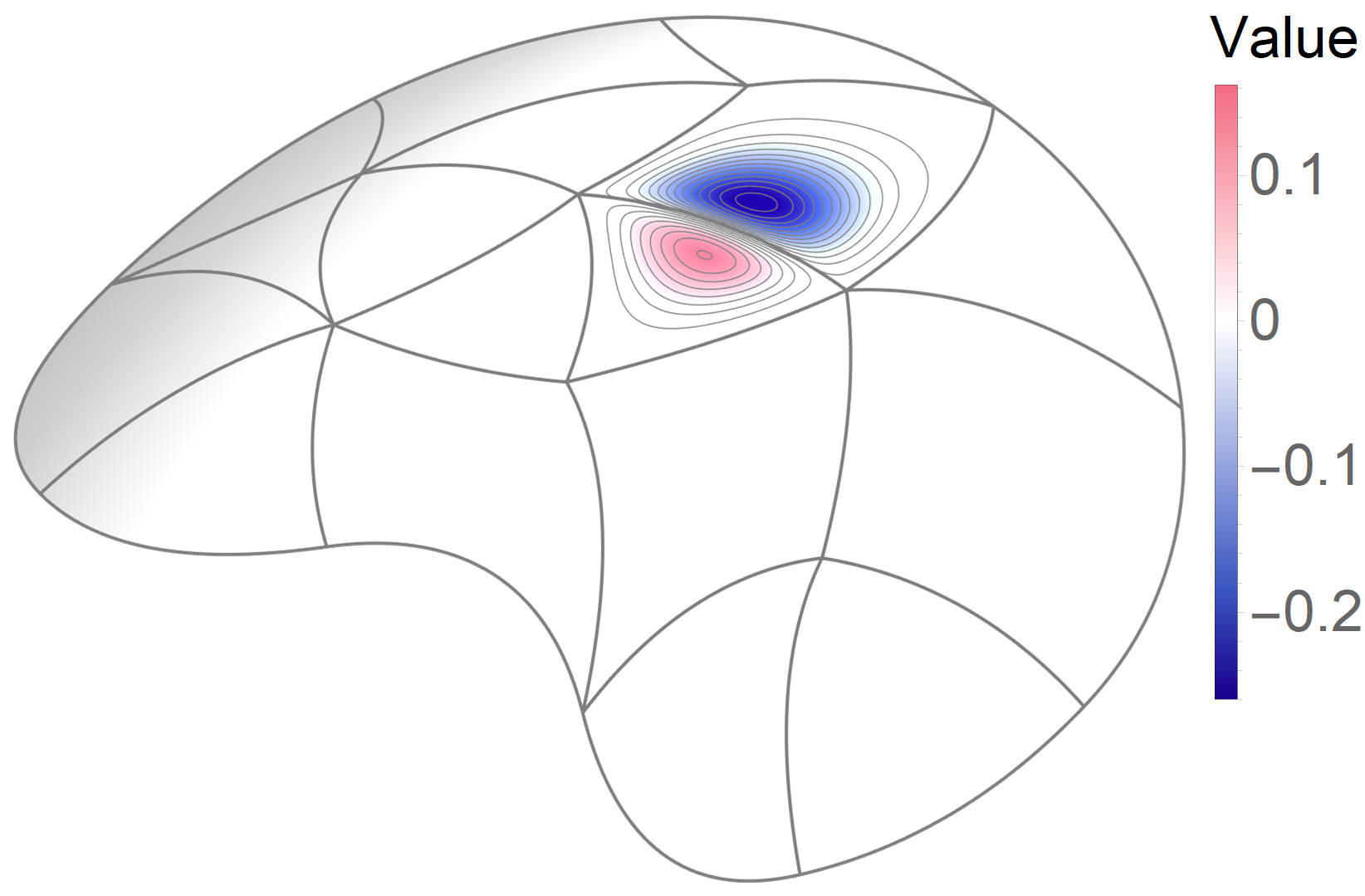}
\end{minipage}
\begin{minipage}{0.3\textwidth}
\centering
\includegraphics[width=1\textwidth]{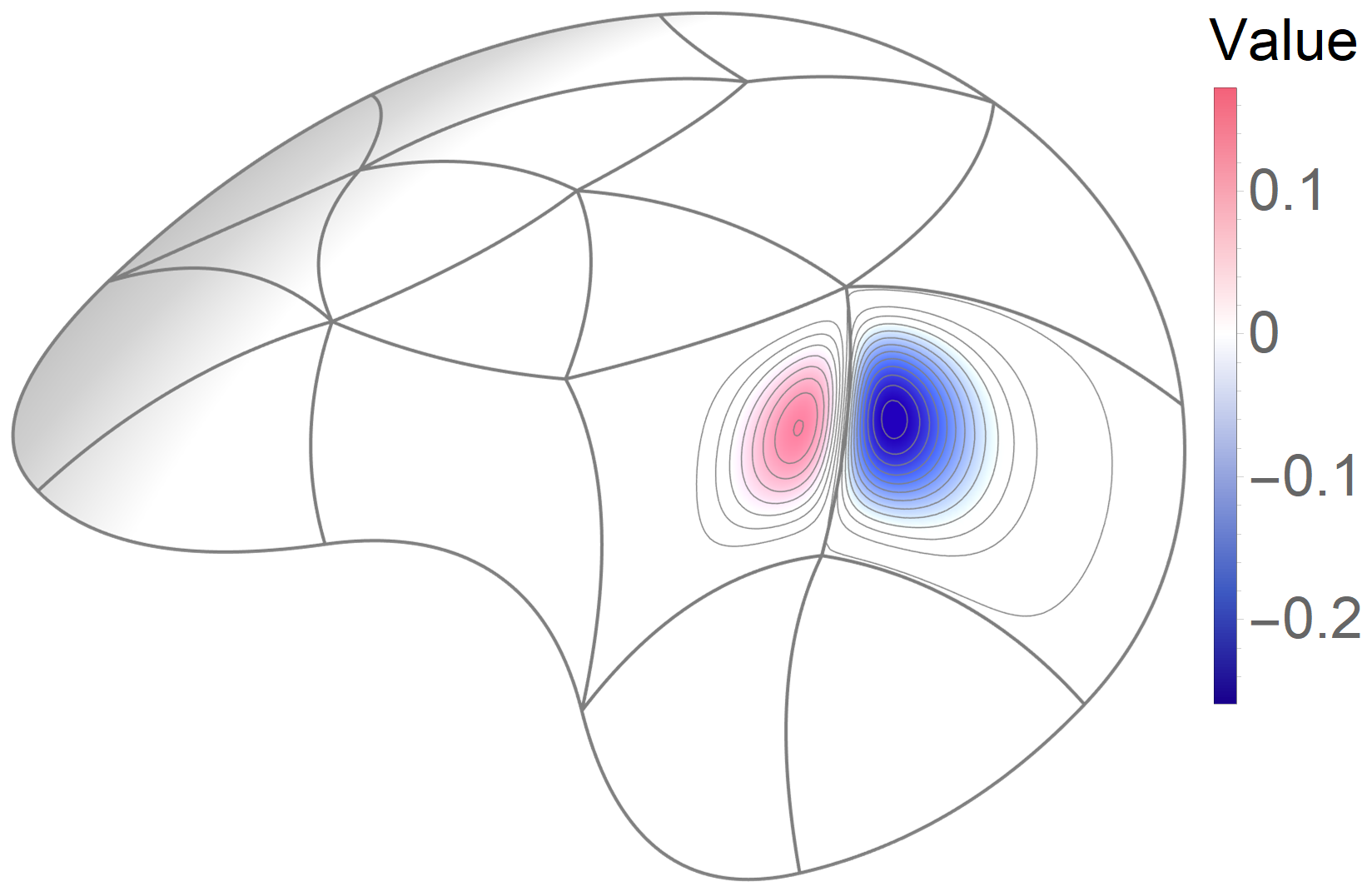}
\end{minipage}
\begin{minipage}{0.3\textwidth}
\centering
\includegraphics[width=1\textwidth]{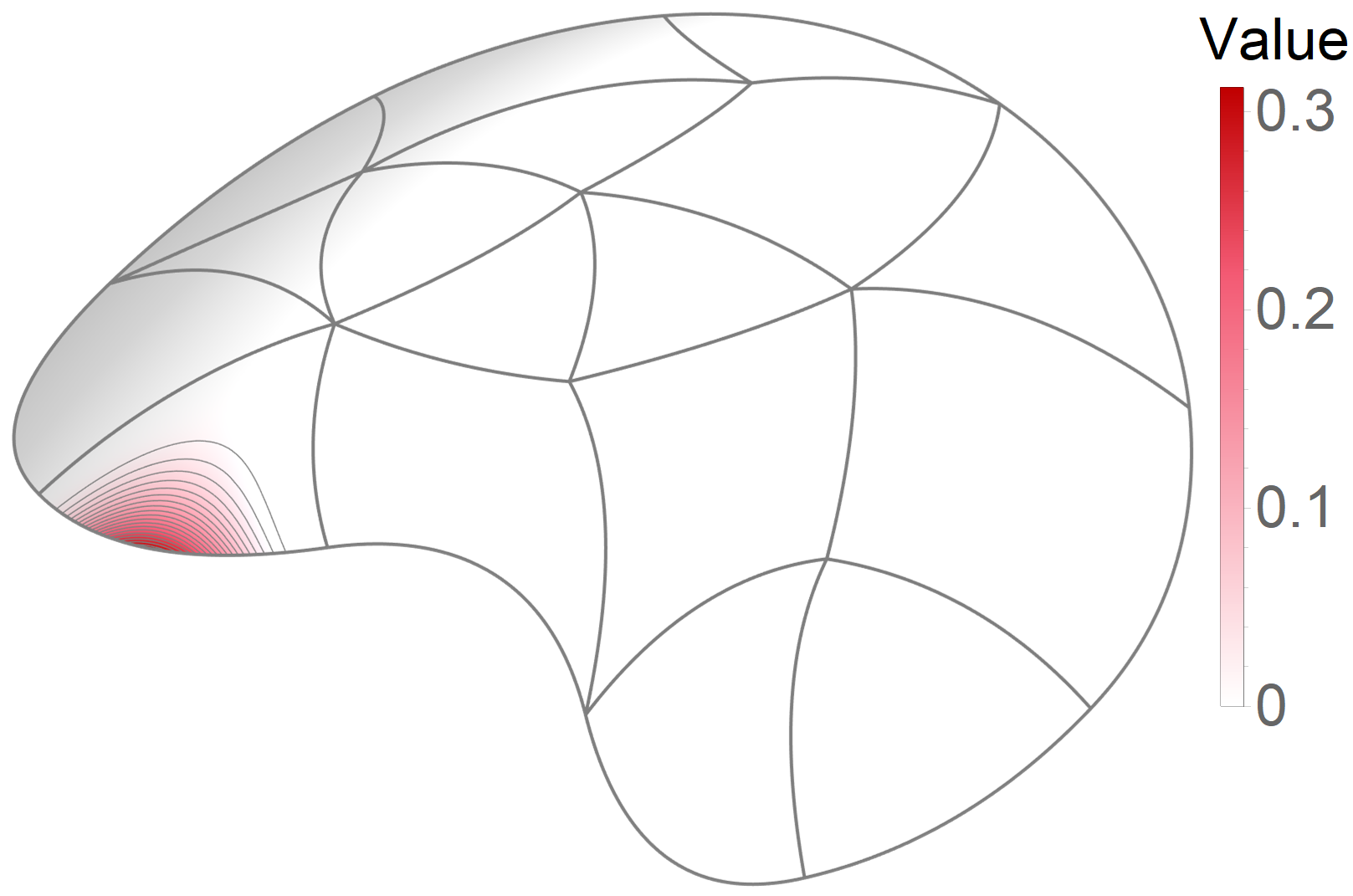}
\end{minipage}
\end{center}
\caption{Selected edge basis functions of degree $6$ constructed over the domain partitioned by {\tt Mesh}~2.}
\label{fig:BasisFunctionsE}
\end{figure}

\begin{figure}[h!]
\begin{center}
\begin{minipage}{0.3\textwidth}
\centering
\includegraphics[width=1\textwidth]{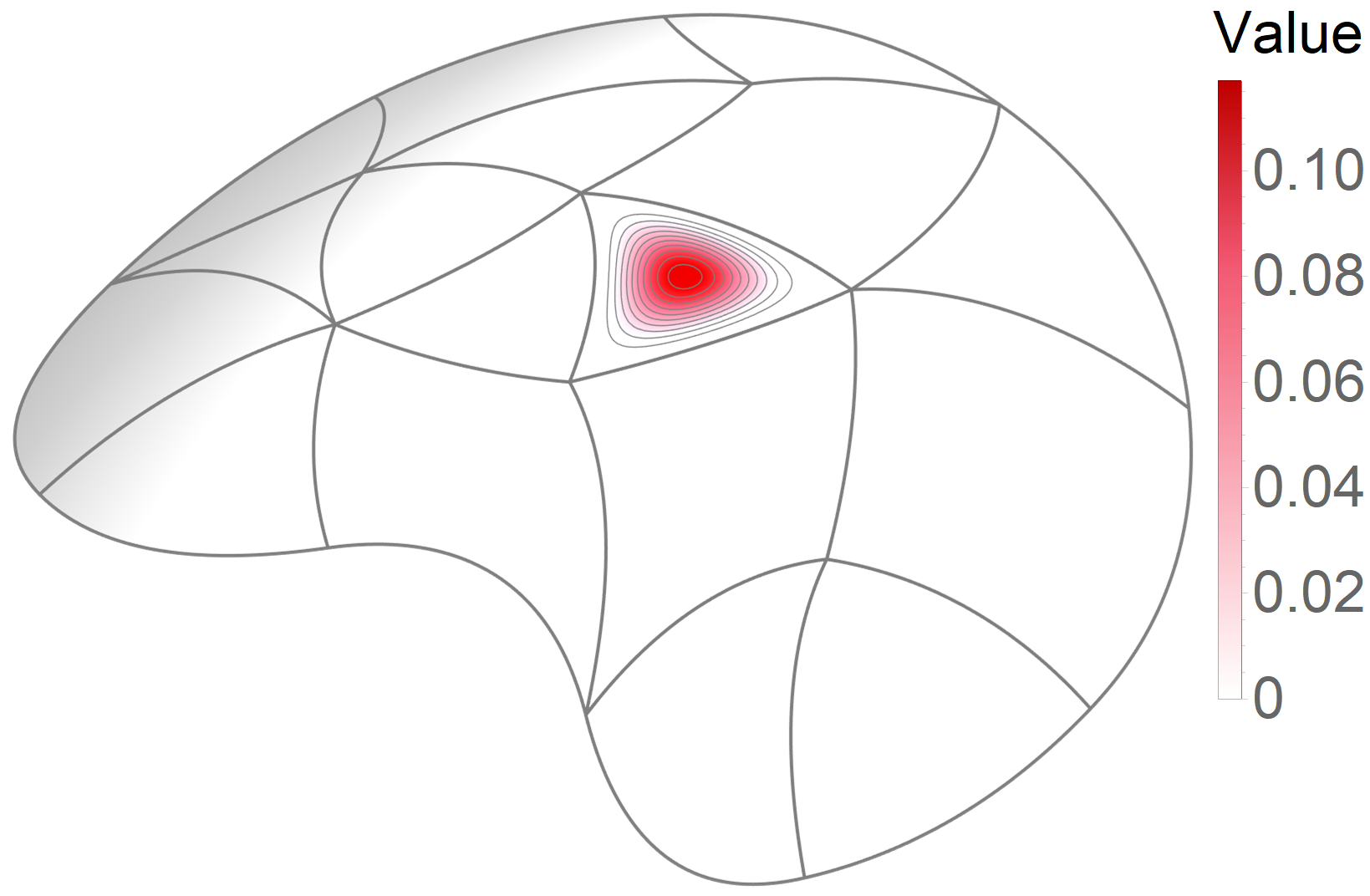}
\end{minipage}
\begin{minipage}{0.3\textwidth}
\centering
\includegraphics[width=1\textwidth]{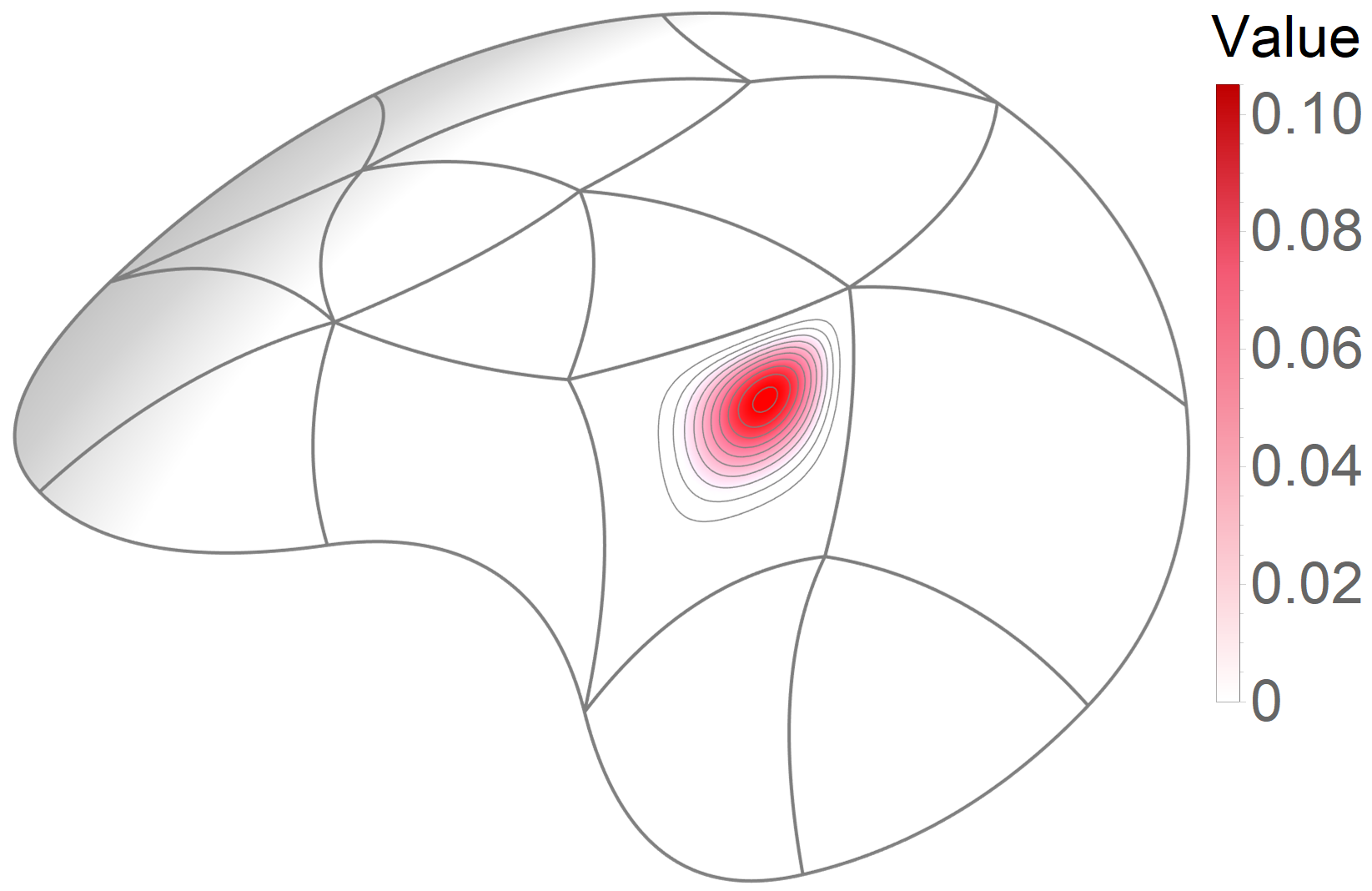}
\end{minipage}
\begin{minipage}{0.3\textwidth}
\centering
\includegraphics[width=1\textwidth]{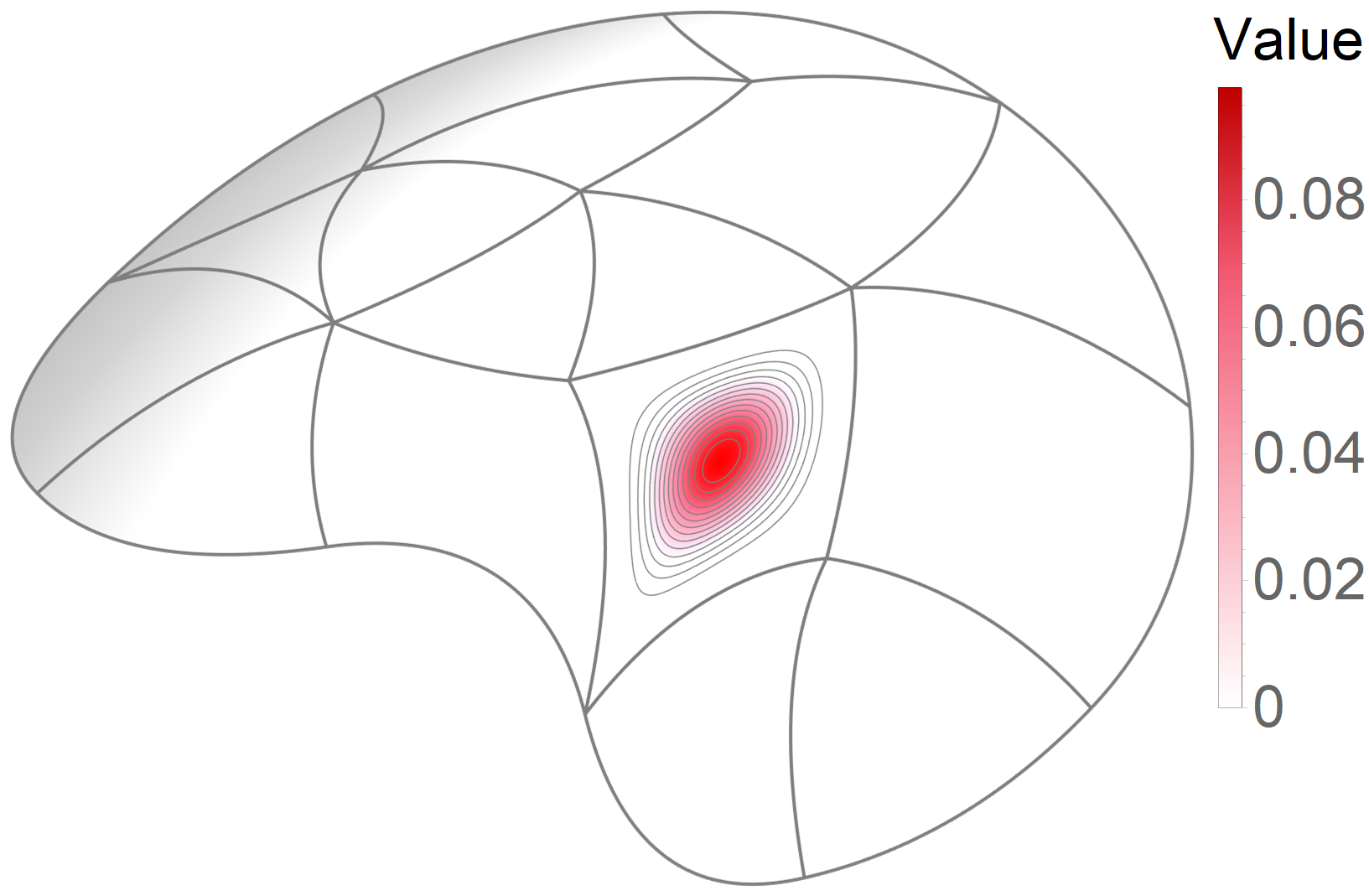}
\end{minipage}
\end{center}
\caption{{Selected triangle and quadrilateral basis functions of degree $6$ constructed over the domain partitioned by {\tt Mesh}~2.}} 
\label{fig:BasisFunctionsT}
\end{figure}

Using the values reported in Table~\ref{tab:edge-dof}, the dimension of the corresponding spline space can be computed directly from \eqref{eq-dimension}, and the basis functions can then be constructed as described in Subsection~\ref{subsec:basis}. Selected basis functions of degree $6$ over the domain partitioned by {\tt Mesh}~2 are illustrated in Figures~\ref{fig:BasisFunctionsV}, \ref{fig:BasisFunctionsE}, and
\ref{fig:BasisFunctionsT}. Similar graphs  are obtained for {\tt Mesh}~1 and {\tt Mesh}~3.

\subsection{Interpolation}
\label{Examples-int}

\begin{figure}[hbt]
\begin{center}
\begin{minipage}{0.32\textwidth}
\centering
\includegraphics[width=0.9\textwidth]{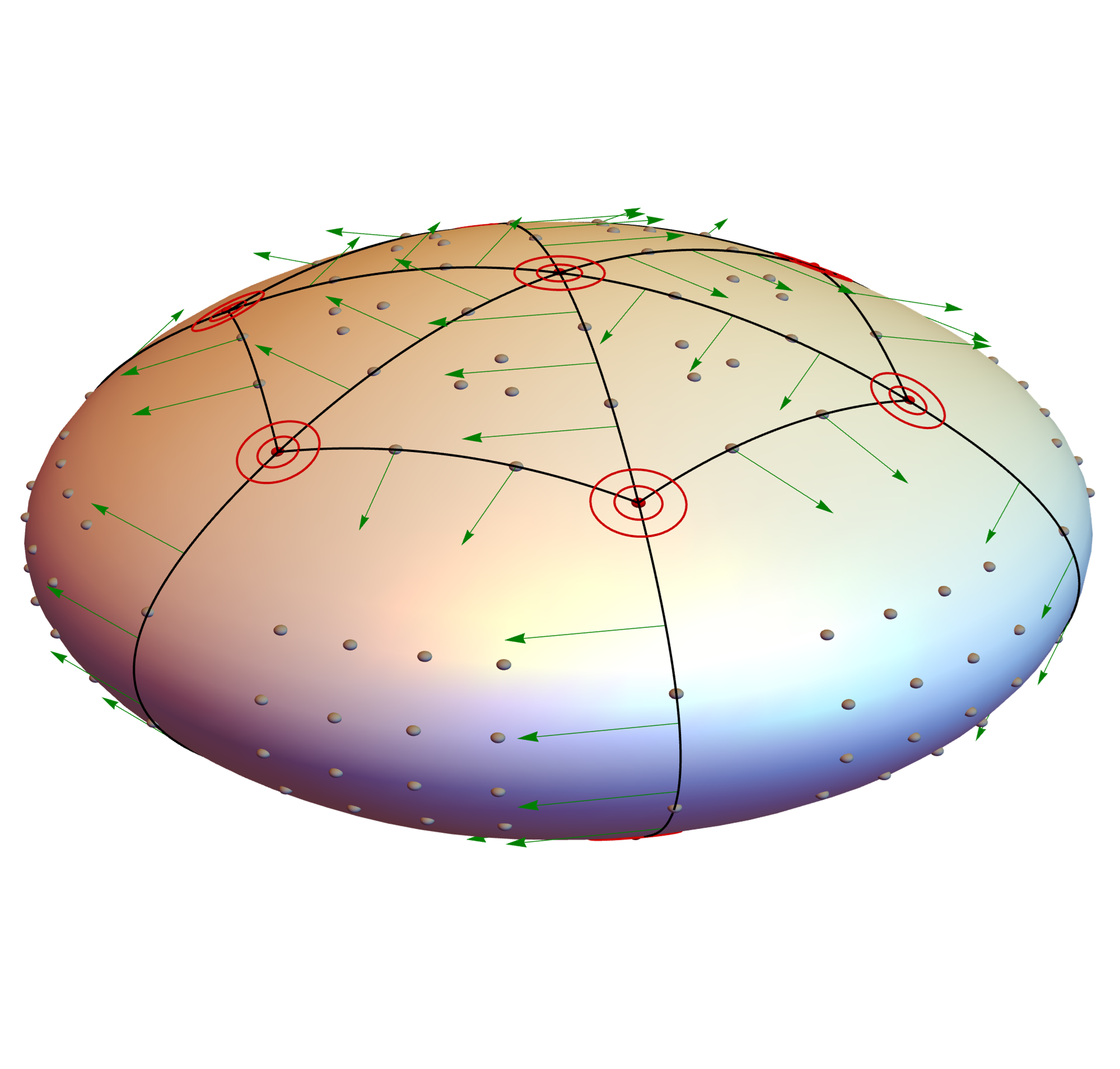}
\end{minipage}
\begin{minipage}{0.32\textwidth}
\centering
\includegraphics[width=1\textwidth]{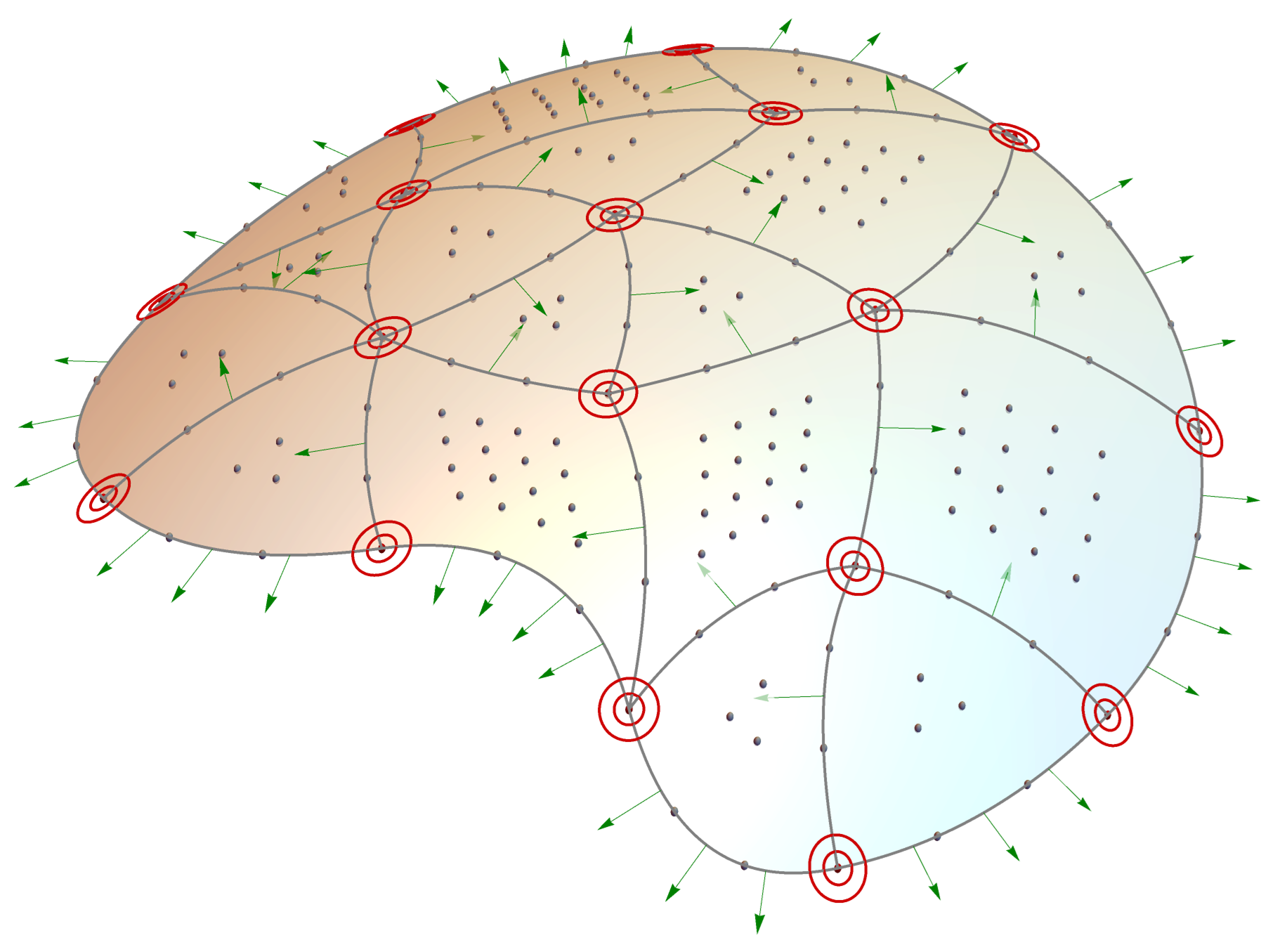}
\end{minipage}
\begin{minipage}{0.32\textwidth}
\centering
\includegraphics[width=1\textwidth]{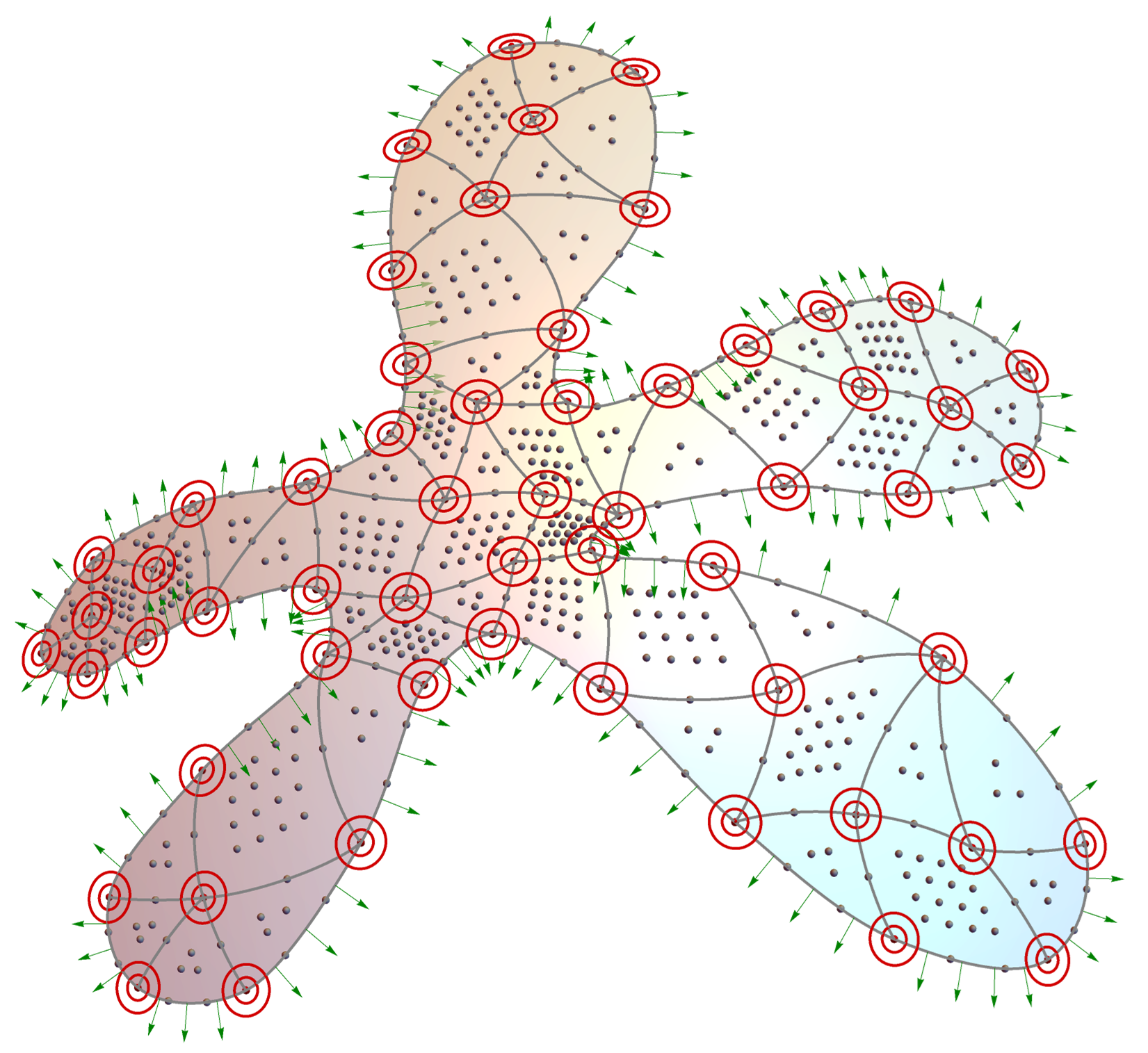}
\end{minipage}
\end{center}
\caption{Distribution of the interpolation conditions for spline degree $d=7$ on Meshes~1--3.}
\label{fig:IntCond-deg7}
\end{figure}

\begin{figure}[hbt]
\begin{center}
\begin{minipage}{0.32\textwidth}
\centering
\includegraphics[width=0.9\textwidth]{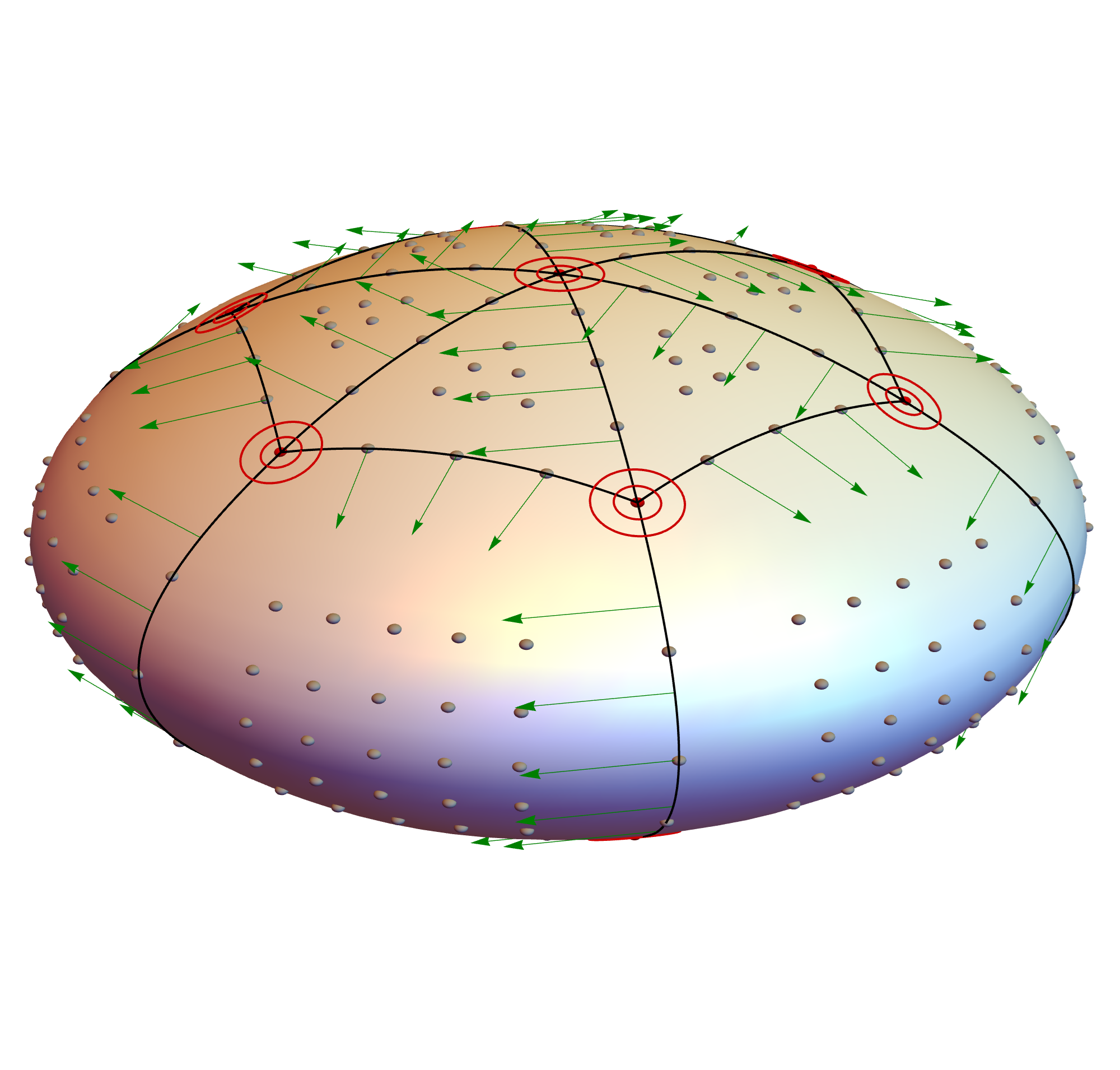}
\end{minipage}
\begin{minipage}{0.32\textwidth}
\centering
\includegraphics[width=1\textwidth]{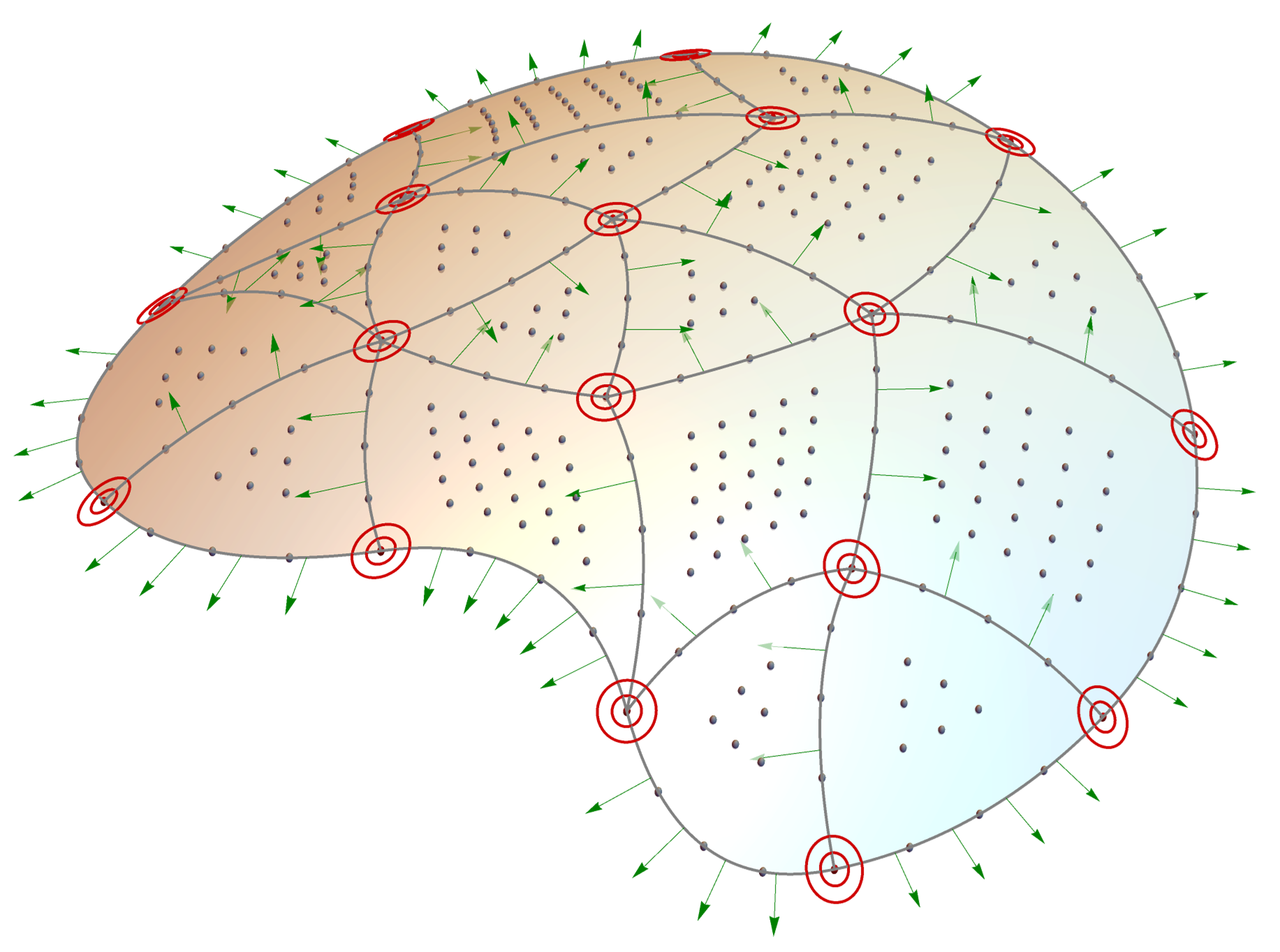}
\end{minipage}
\begin{minipage}{0.32\textwidth}
\centering
\includegraphics[width=1\textwidth]{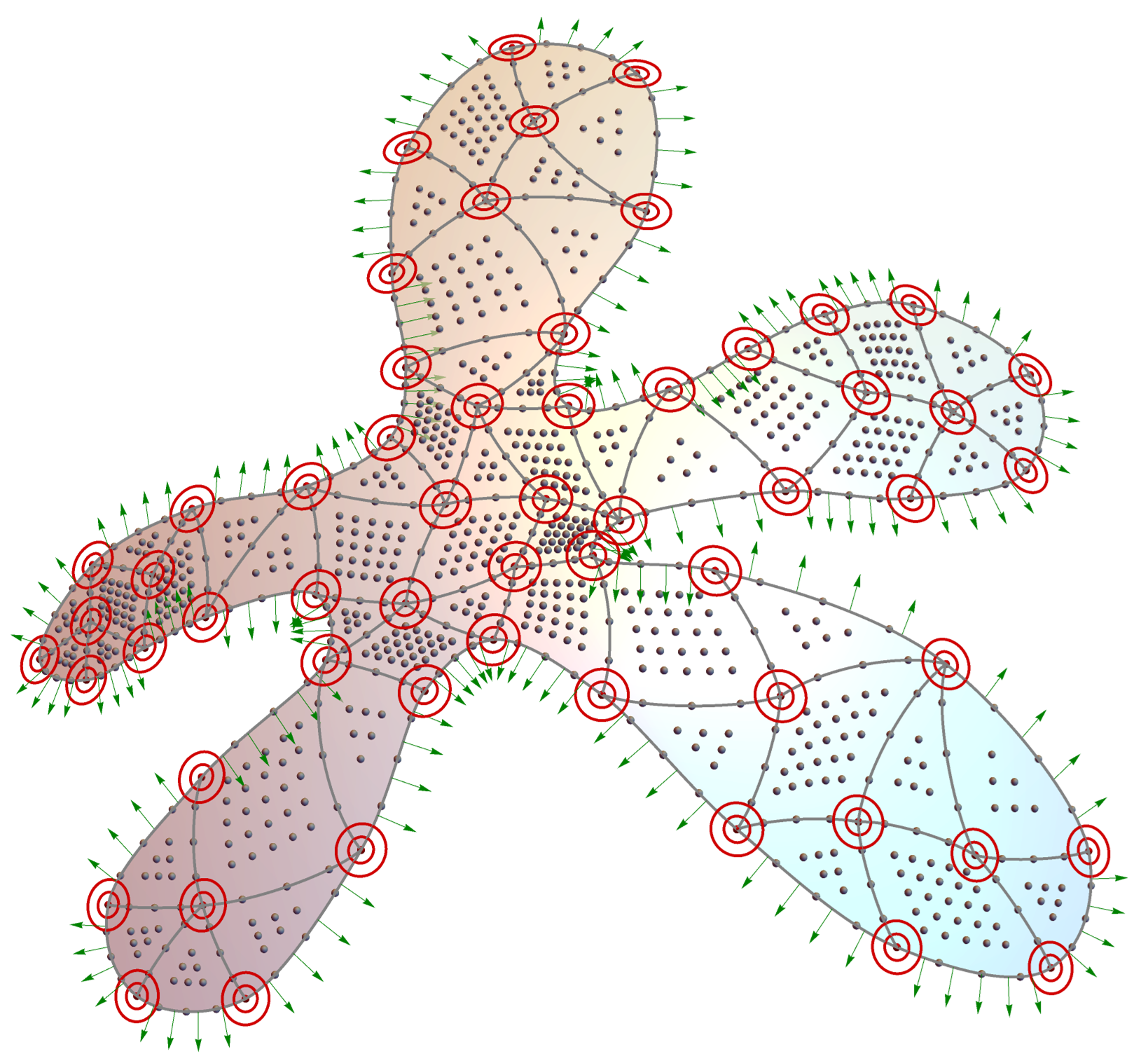}
\end{minipage}
\end{center}
\caption{Distribution of the interpolation conditions for spline degree $d=8$ on Meshes~1--3.}
\label{fig:IntCond-deg8}
\end{figure}

We first consider the interpolation scheme introduced in Subsection~\ref{subsec:interpolation}. 
For the edge interpolation parameters in \eqref{subeq-interpolFunctionals-2}, we choose the uniformly distributed
values 
$$\parV^{0}_{\ell,j} = \frac{j}{\chi_{\ell,2}+1}, \quad j=1,\dots,\chi_{\ell,2}, 
\quad {\rm and} \quad \parV^{1}_{\ell,j} = \frac{j}{\chi_{\ell,3}+1}, \quad j=1,\dots,\chi_{\ell,3}.
$$
The resulting interpolation conditions for spline degrees $d=7$ and $d=8$ for all three meshes are illustrated in Figures~\ref{fig:IntCond-deg7} and~\ref{fig:IntCond-deg8}. These conditions also reflect the edge degrees of freedom reported in Table~\ref{tab:edge-dof}.

The interpolation data are obtained from a sufficiently smooth trivariate function $g$ defined on a three-dimensional domain containing $\Omega$. The approximation error between $g$ and its interpolant $\varphi$ is estimated
by the discrete maximum norm 
\begin{equation} \label{eq-def-maxError}
\|g-\varphi\|_{\infty} \approx  {\rm err}_{\infty,N}(g, \varphi):=\|g-\varphi\|_{\infty,N} := \max_{k \in \indOmega} \max_{\bfm{\xi} \in \Xi_N} \big|g (\bfm{F}^{(k)}(\bfm{\xi} )) - f^{(k)}(\bfm{\xi} )\big|,
\end{equation}
where $f^{(k)}$ is the local representation of $\varphi$ over the surface
element $\Omega^{(k)}$, and $\Xi_N$ denotes a uniform set of evaluation
parameters in the corresponding parametric domain. More precisely,
$\Xi_N = \defset{\left(\tfrac{i_1}{N},\tfrac{i_2}{N}\right)}{\mathbf i \in \indSetT{N}}$ for $\parDomT$ and
$\Xi_N = \defset{\left(\tfrac{i_1}{N},\tfrac{i_2}{N}\right)}{\mathbf i \in \indSetQ{N}}$ for $\parDomQ$.

\begin{example} \label{example:interpolation1}
Consider the surfaces corresponding to Meshes~1--3 and let
\begin{equation} \label{def:function-g-1}
g(x,y,z)= 2 \, e^{-\frac{1}{2} ((x-4)^2+(y-5)^2+(z+1)^2)}+\sin \left(\frac{\pi  x}{4}\right) \sin \left(\frac{\pi  y}{4}\right) \cos \left(\frac{\pi  z}{8}\right)
\end{equation}
be the chosen test function. Contour plots of the restriction of $g$ to the three surface domains are shown in Figure~\ref{fig:IntFun1}.
\begin{figure}[hbt]
\begin{center}
\begin{minipage}{0.32\textwidth}
\centering
\includegraphics[width=1\textwidth]{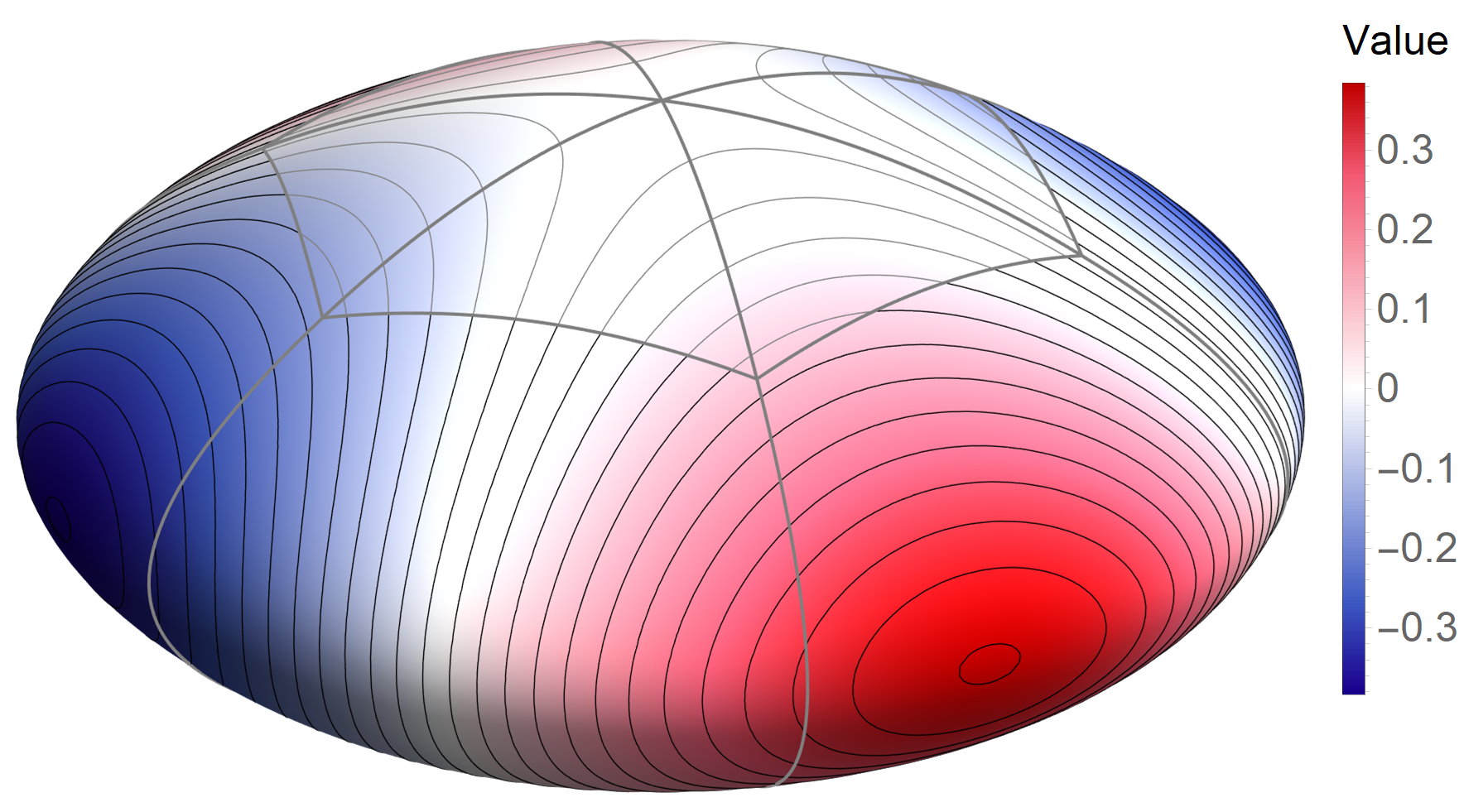}
\end{minipage}
\begin{minipage}{0.32\textwidth}
\centering
\includegraphics[width=0.95\textwidth]{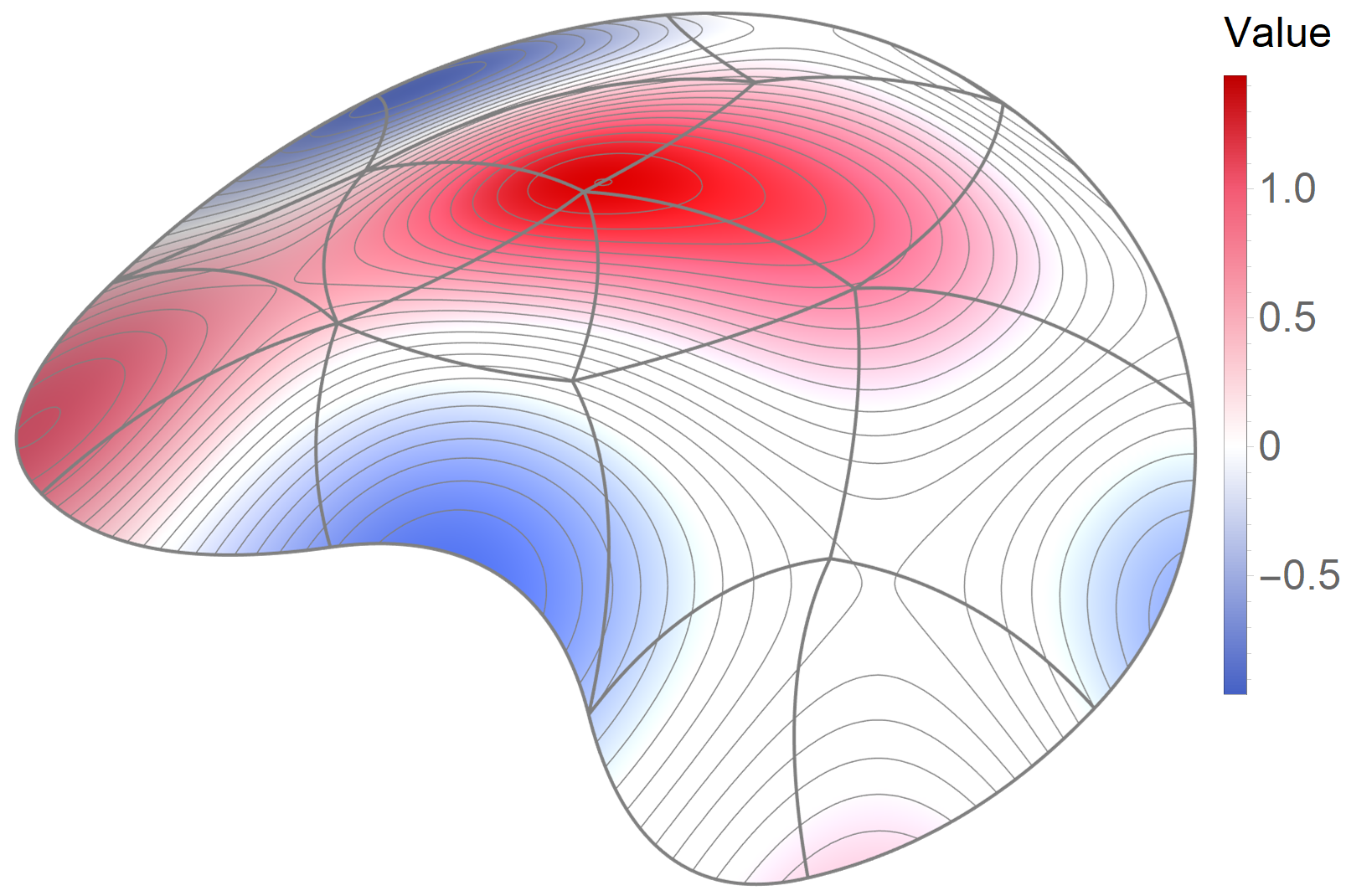}
\end{minipage}
\begin{minipage}{0.32\textwidth}
\centering
\includegraphics[width=0.9\textwidth]{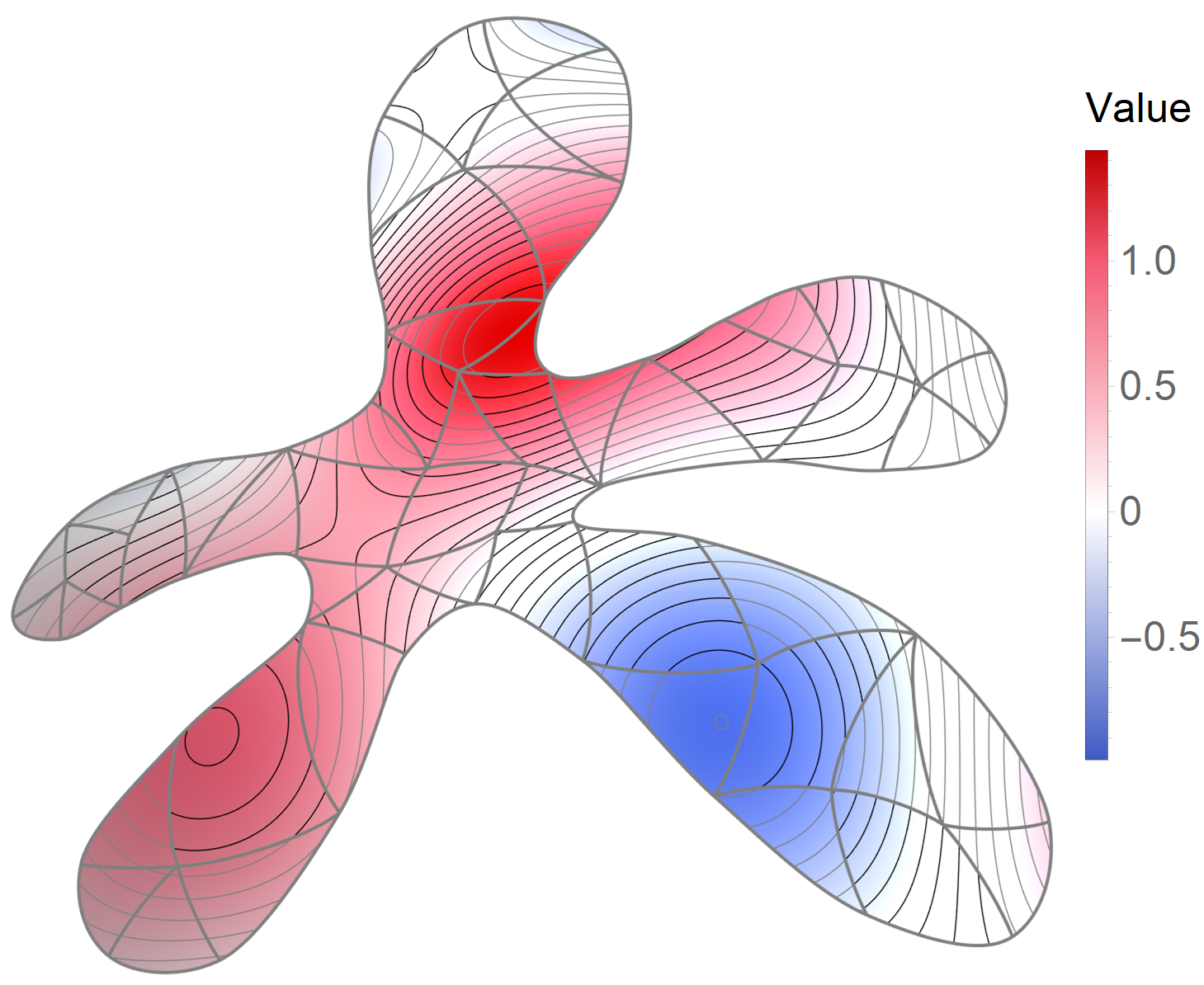}
\end{minipage}
\end{center}
\caption{Contour plots of the restriction of $g$, given by \eqref{def:function-g-1}, to the surface domains represented by Meshes~1--3.}
\label{fig:IntFun1}
\end{figure}

For Mesh~1 and Mesh~2, we compute the interpolating splines of degrees $d=6,7,8,9,10,11$, while for Mesh~3 we consider degrees
$d=7,8,9,10,11$. We evaluate the error \eqref{eq-def-maxError} for $N=50$. The errors are reported in Table~\ref{tab:err-exm-1}, together with the dimensions of the corresponding spline spaces. Figure~\ref{fig:err-exm-1} provides a graphical illustration of the decay of the interpolation error with respect to the spline degree and the square root of the number of degrees of freedom.
\begin{table}[h!]
\begin{align*}
\renewcommand{\arraystretch}{1.1}
\begin{array}{|c|c|c|}
\hline
 \multicolumn{3}{|c|}{{\tt Mesh}\, 1} \\
\hline
d & {\rm NDOF} & {\rm err}_{\infty,N} \\
\hline
6  & 234  & 1.77364\cdot 10^{-3} \\
7  & 360  & 1.48834\cdot 10^{-3} \\
8  & 510  & 5.55784\cdot 10^{-4} \\
9  & 684  & 4.37672\cdot 10^{-4} \\
10 & 882  & 2.37634\cdot 10^{-4} \\
11 & 1104 & 1.91113\cdot 10^{-4} \\
\hline
\end{array}
\qquad 
\begin{array}{|c|c|c|}
\hline
 \multicolumn{3}{|c|}{{\tt Mesh}\, 2} \\
\hline
d & {\rm NDOF} &  {\rm err}_{\infty,N} \\
\hline
6  & 215  & 4.45887\cdot 10^{-2} \\
7  & 340  & 6.72468\cdot 10^{-3} \\
8  & 487  & 4.98082\cdot 10^{-3} \\
9  & 656  & 1.95684\cdot 10^{-3} \\
10 & 847  & 2.00559\cdot 10^{-3} \\
11 & 1060 & 8.19947\cdot 10^{-4} \\
\hline
\end{array}
\qquad 
\begin{array}{|c|c|c|}
\hline
 \multicolumn{3}{|c|}{{\tt Mesh}\, 3} \\
\hline
d & {\rm NDOF} & {\rm err}_{\infty,N} \\
\hline
6  & -  & - \\
7  & 1021 & 1.20815\cdot 10^{-2} \\
8  & 1501 & 3.02609\cdot 10^{-2} \\
9  & 2051 & 2.49053\cdot 10^{-3} \\
10 & 2671 & 4.68947\cdot 10^{-3} \\
11 & 3361 & 8.08770\cdot 10^{-4} \\
\hline
\end{array}
\qquad
\end{align*}
\caption{
Example~\ref{example:interpolation1}: {Discrete maximum interpolation errors} with $N=50$ for spline degrees $d=6,\ldots,11$ for Mesh~1 and Mesh~2, and $d=7,\ldots,11$ for Mesh~3, together with the corresponding numbers of degrees of freedom.
}
\label{tab:err-exm-1}
\end{table}

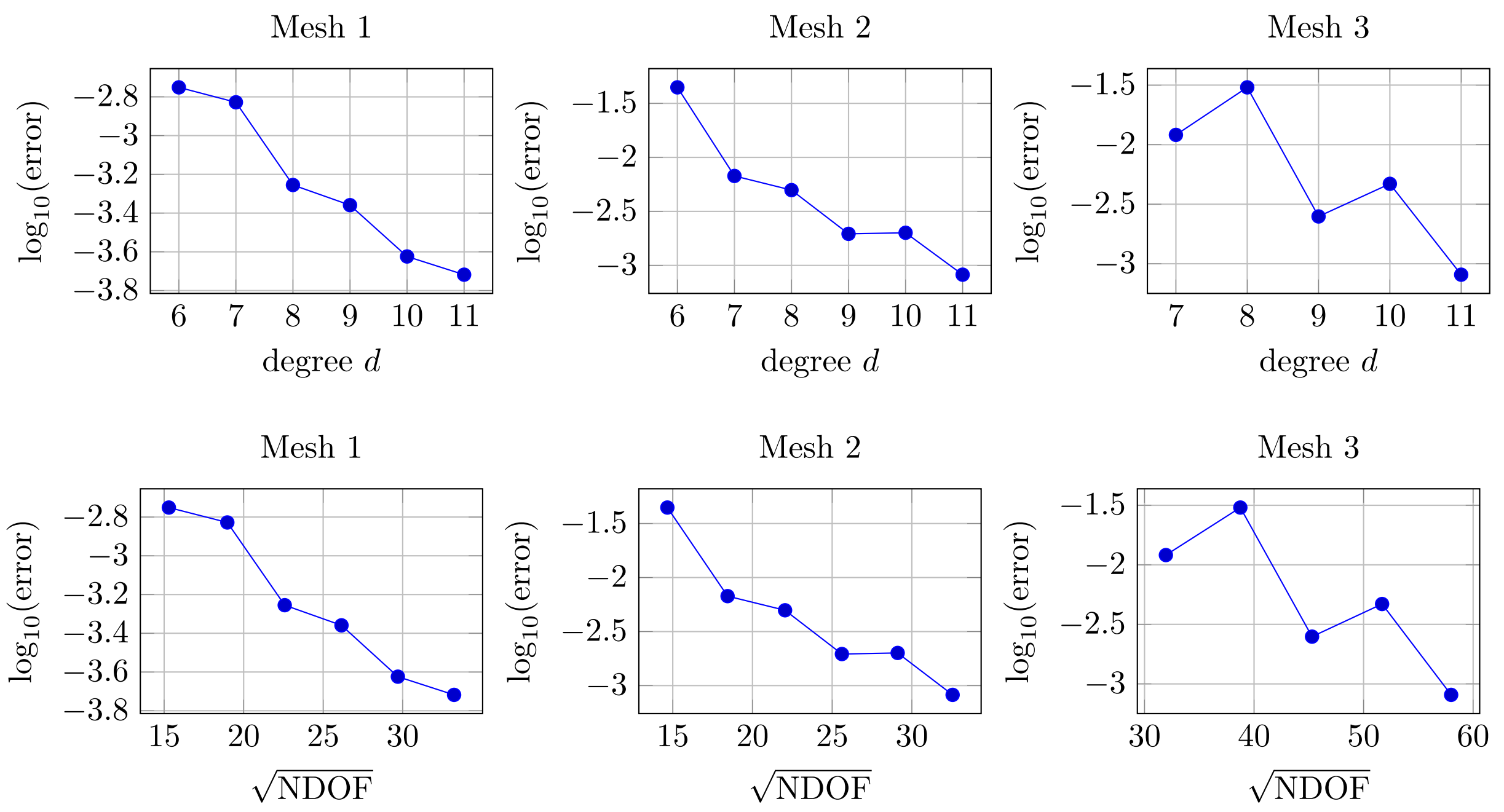
\begin{figure}[h!]
\begin{center}
\begin{tikzpicture}
\begin{axis}[
    width=.32\textwidth, height=4cm,
    xlabel={degree $d$},
    ylabel={$\log_{10}(\text{error})$},
    xtick={6, 7, 8, 9, 10, 11},
    ytick={-2.8,-3,-3.2,-3.4,-3.6,-3.8},
    grid=major,
    title={Mesh 1}
]
\addplot coordinates {
    (6, -2.751)
    (7, -2.828)
    (8, -3.255)
    (9, -3.359)
    (10, -3.624)
    (11, -3.718)
};
\end{axis}
\end{tikzpicture}
\begin{tikzpicture}
\begin{axis}[
    width=.32\textwidth, height=4cm,
    xlabel={degree $d$},
    ylabel={$\log_{10}(\text{error})$},
    xtick={6, 7, 8, 9, 10, 11},
    grid=major,
    title={Mesh 2}
]
\addplot coordinates {
    (6, -1.352)
    (7, -2.172)
    (8, -2.303)
    (9, -2.708)
    (10, -2.698)
    (11, -3.086)
};
\end{axis}
\end{tikzpicture}
\begin{tikzpicture}
\begin{axis}[
    width=.32\textwidth, height=4cm,
    xlabel={degree $d$},
    ylabel={$\log_{10}(\text{error})$},
    xtick={7, 8, 9, 10, 11},
    grid=major,
    title={Mesh 3}
]
\addplot coordinates {
    (7, -1.918)
    (8, -1.519)
    (9, -2.603)
    (10, -2.329)
    (11, -3.092)
};
\end{axis}
\end{tikzpicture}
\\[10pt]
\begin{tikzpicture}
\begin{axis}[
    width=.32\textwidth, height=4cm,
    xlabel={$\sqrt{\text{NDOF}}$},
    ylabel={$\log_{10}(\text{error})$},
    ytick={-2.8,-3,-3.2,-3.4,-3.6,-3.8},
    xtick={15,20,25,30},
    grid=major,
    title={Mesh 1}
]
\addplot coordinates {
    (15.297, -2.751)
    (18.974, -2.828)
    (22.583, -3.255)
    (26.153, -3.359)
    (29.698, -3.624)
    (33.233, -3.718)
};
\end{axis}
\end{tikzpicture}
\begin{tikzpicture}
\begin{axis}[
    width=.32\textwidth, height=4cm,
    xlabel={$\sqrt{\text{NDOF}}$},
    ylabel={$\log_{10}(\text{error})$},
    grid=major,
    title={Mesh 2}
]
\addplot coordinates {
    (14.662, -1.352)
    (18.439, -2.172)
    (22.045, -2.303)
    (25.612, -2.708)
    (29.110, -2.698)
    (32.558, -3.086)
};
\end{axis}
\end{tikzpicture}
\begin{tikzpicture}
\begin{axis}[
    width=.32\textwidth, height=4cm,
    xlabel={$\sqrt{\text{NDOF}}$},
    ylabel={$\log_{10}(\text{error})$},
    grid=major,
    title={Mesh 3}
]
\addplot coordinates {
    (31.953, -1.918)
    (38.748, -1.519)
    (45.293, -2.603)
    (51.688, -2.329)
    (57.993, -3.092)
};
\end{axis}
\end{tikzpicture}
\end{center}
\caption{Discrete maximum interpolation errors for Meshes~1--3 as functions of the spline degree $d$ (top) and of the square root of the number of degrees of freedom (bottom), shown on a logarithmic scale. 
}
\label{fig:err-exm-1}
\end{figure}

Although the interpolation error decreases overall as the spline degree and the number of degrees of freedom increase, the convergence with respect to the spline degree is not monotone. As often observed in polynomial interpolation, the error levels for even and odd degrees are different. In particular, for Mesh~3 a slight increase in the error is observed at even degrees.
\end{example}
\begin{example} \label{example:interpolation2}
In this example, we numerically verify the polynomial reproduction property of the proposed interpolation scheme. To this end, we consider two planar gear-shaped domains with a hole, represented by Mesh~4 and Mesh~5, shown in Figure~\ref{fig:exm-2}. 
\begin{figure}[h!]
\begin{center}
\begin{minipage}{0.45\textwidth}
\includegraphics[width=1.05\textwidth]{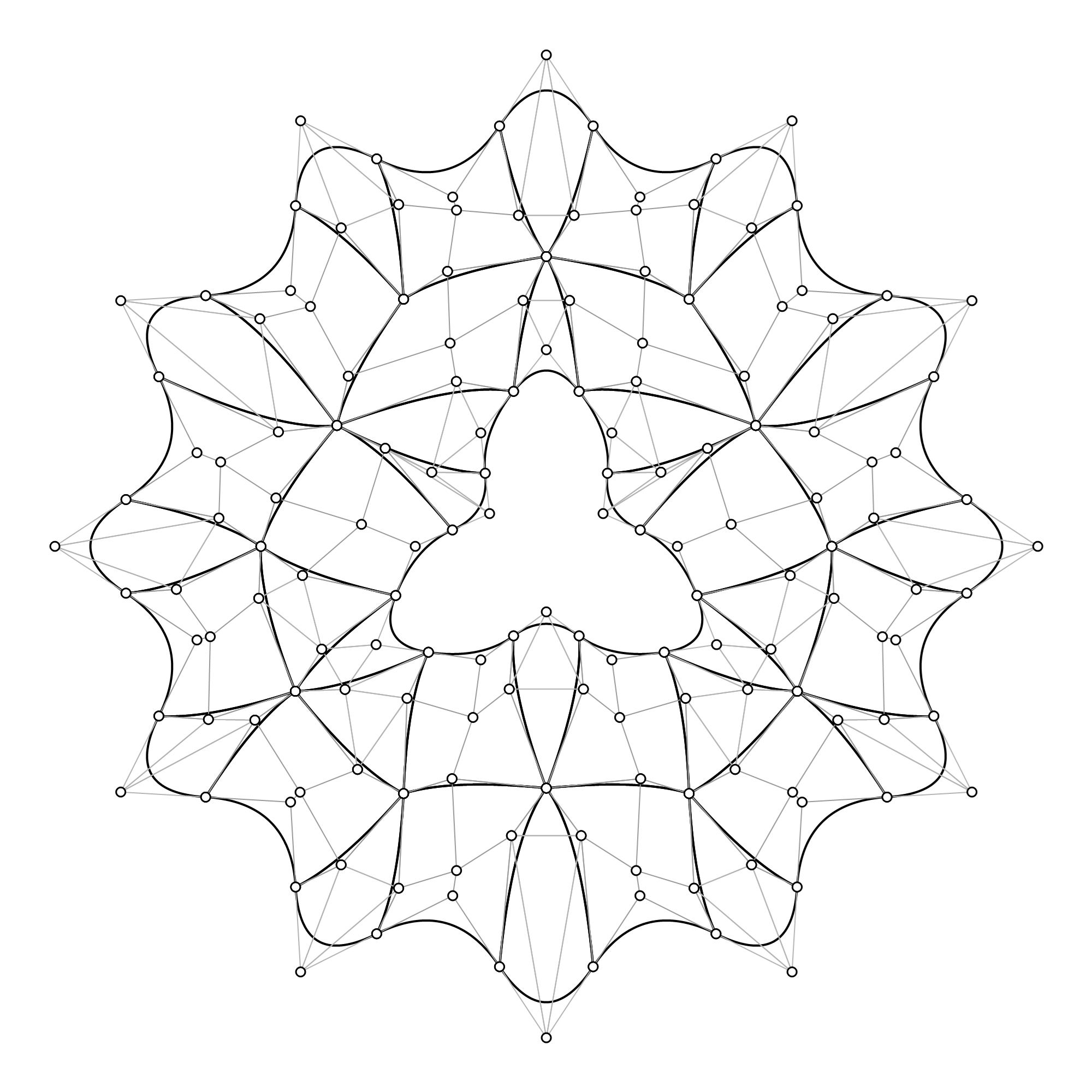}
\end{minipage}
\quad
\begin{minipage}{0.45\textwidth}
\includegraphics[width=1\textwidth]{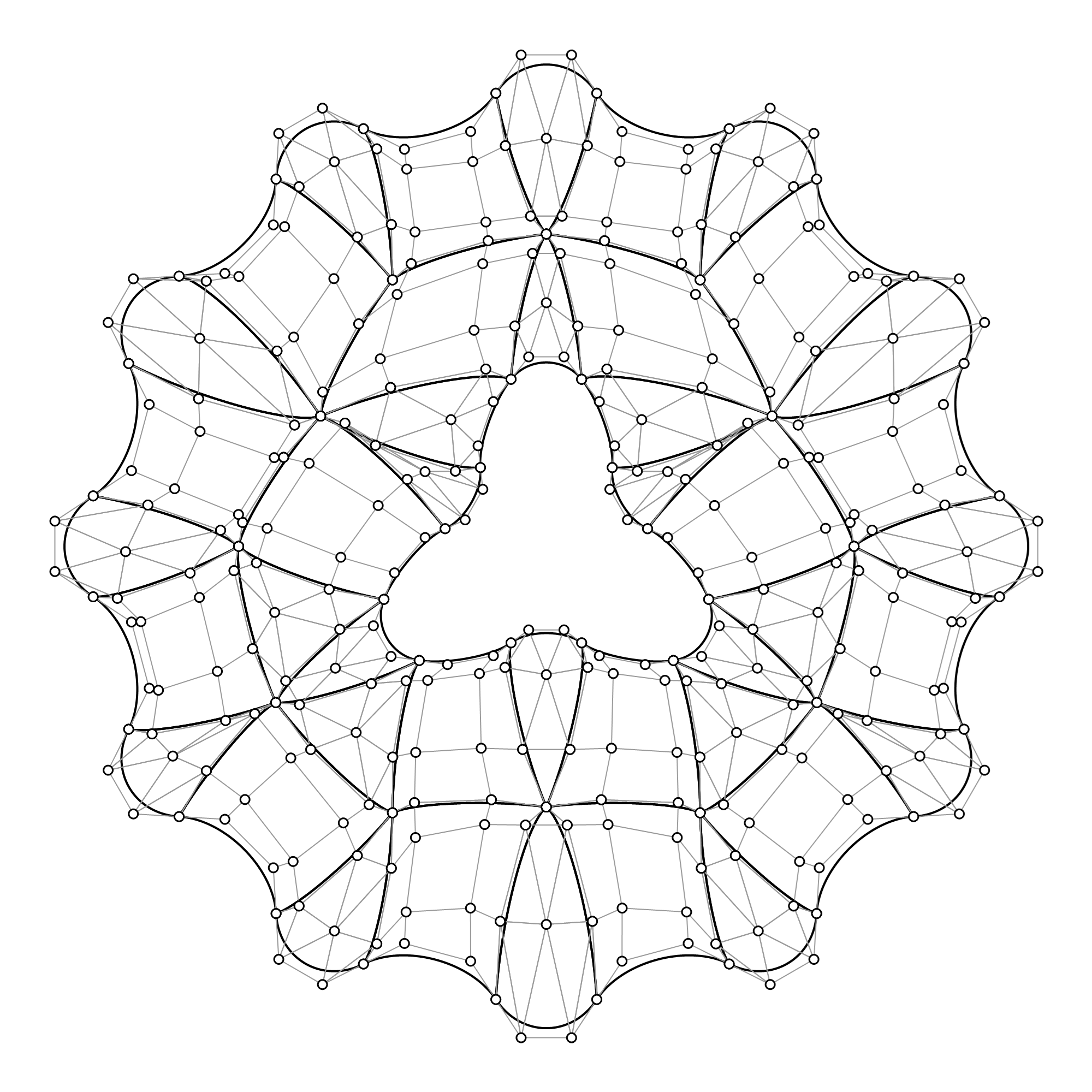}
\end{minipage}
\end{center}
\caption{Planar gear-shaped domains with a hole represented by Mesh~4 (left) and Mesh~5 (right), with the control meshes illustrating the corresponding geometry mappings. The bounding box for the control points is $[-1.25,1.25]^2$ for both meshes.}
\label{fig:exm-2}
\end{figure}
The geometry mappings defining Mesh~4 are quadratic. As in the previous example, the minimum number of degrees of freedom associated with an edge is $2d-1$, and hence the spline construction is applicable for $d\geq 6$. For Mesh~5, which is given by cubic geometry mappings, the computation of the edge degrees of freedom using Algorithm~\ref{alg:DegreesOfFreedom} shows that the minimum number of degrees of freedom per edge is $2d-3$, and that the spline construction requires $d\geq 7$. 
We consider the polynomial test functions
$$
g_k(x,y)=2^{-k-2}(x-5)(y-4)(x+y)^{k-2}, \quad k=2,3,\ldots,6,
$$
and compute their interpolating splines of degrees $d=6,7,\dots,11$ for Mesh~4 and $d=7,8\dots,11$ for Mesh~5. 

For Mesh~4, the numerical results show that polynomials of degree up to $\lfloor \tfrac{d}{2}\rfloor$ are reproduced exactly, up to numerical precision.
The errors for the remaining combinations of $k$ and $d$ are reported in
Table~\ref{tab:poly-reproduction-mesh4}.
For Mesh~5, the same experiment shows exact reproduction, up to numerical precision, for polynomial degrees up to $\lfloor \tfrac{d}{3}\rfloor$.
The corresponding errors for the remaining combinations are listed in Table~\ref{tab:poly-reproduction-mesh5}.
Interpolating splines of degree $9$ for the test functions $g_4$, $g_5$, and $g_6$ over Mesh~5 are shown in Figure~\ref{fig:exm-2-int}.
\begin{table}[h!]
\begin{align*}
\renewcommand{\arraystretch}{1.05}
\begin{array}{|c|c|c|c|c|c|c|}
\hline
 & d=6 &  d=7 & d=8 & d=9 & d=10 & d=11\\
\hline
{\rm NDOF} & 652 & 1006 & 1420 & 1894 & 2428 & 3023\\
\hline
k=2 & 0 & 0 & 0 & 0 & 0 & 0\\
k=3 & 0 & 0 & 0 & 0 & 0  & 0 \\
k=4 & 7.4360\cdot 10^{-6} & 2.6332\cdot 10^{-7} & 0 & 0 & 0 & 0 \\
k=5 & 4.8896\cdot 10^{-5} & 2.5057\cdot 10^{-6} & 9.9365\cdot 10^{-8} & 3.5682\cdot 10^{-9} & 0  & 0 \\
k=6 & 1.3379\cdot 10^{-4} & 1.2190\cdot 10^{-5} & 5.4457\cdot 10^{-7} & 3.8516\cdot 10^{-8} & 1.0894\cdot 10^{-9} & 4.6719\cdot 10^{-11}\\
\hline
\end{array}
\end{align*}
\caption{Polynomial reproduction for Mesh~4: Discrete maximum interpolation errors ${\rm err}_{\infty,50}$ between polynomials $g_k$ and interpolating splines of degree $d=6,7,\dots,11$. A zero indicates reproduction up
to numerical precision.}
\label{tab:poly-reproduction-mesh4}
\end{table}

\begin{table}[h!]
\begin{align*}
\renewcommand{\arraystretch}{1.05}
\begin{array}{|c|c|c|c|c|c|}
\hline
 & d=7 & d=8 & d=9 & d=10 & d=11 \\
\hline
{\rm NDOF} & 898 & 1312 & 1786 & 2320 & 2914 \\
\hline
k=2 & 0 & 0 & 0 & 0 & 0 \\
k=3 & 8.3622\cdot 10^{-6} & 2.5356\cdot 10^{-6} & 0 & 0 & 0 \\
k=4 &  6.2846\cdot 10^{-5} & 2.3294\cdot 10^{-5} & 7.2136\cdot 10^{-7} & 5.2029\cdot 10^{-8} & 6.3961\cdot 10^{-9} \\
k=5 & 2.8356\cdot 10^{-4} & 7.4922\cdot 10^{-5} & 4.4407\cdot 10^{-6} & 3.8197\cdot 10^{-7} & 5.0638\cdot 10^{-8} \\
k=6 & 7.3564\cdot 10^{-4} & 1.5120\cdot 10^{-4} & 2.0610\cdot 10^{-5} & 1.0754\cdot 10^{-6} & 2.2433\cdot 10^{-7} \\
\hline
\end{array}
\end{align*}
\caption{Polynomial reproduction for Mesh~5: Discrete maximum interpolation errors ${\rm err}_{\infty,50}$ between polynomials $g_k$ and interpolating splines of degree $d=7, 8, \dots, 11$. A zero indicates reproduction up
to numerical precision.}
\label{tab:poly-reproduction-mesh5}
\end{table}
\begin{figure}[hbt]
\centering
\includegraphics[width=1\textwidth]{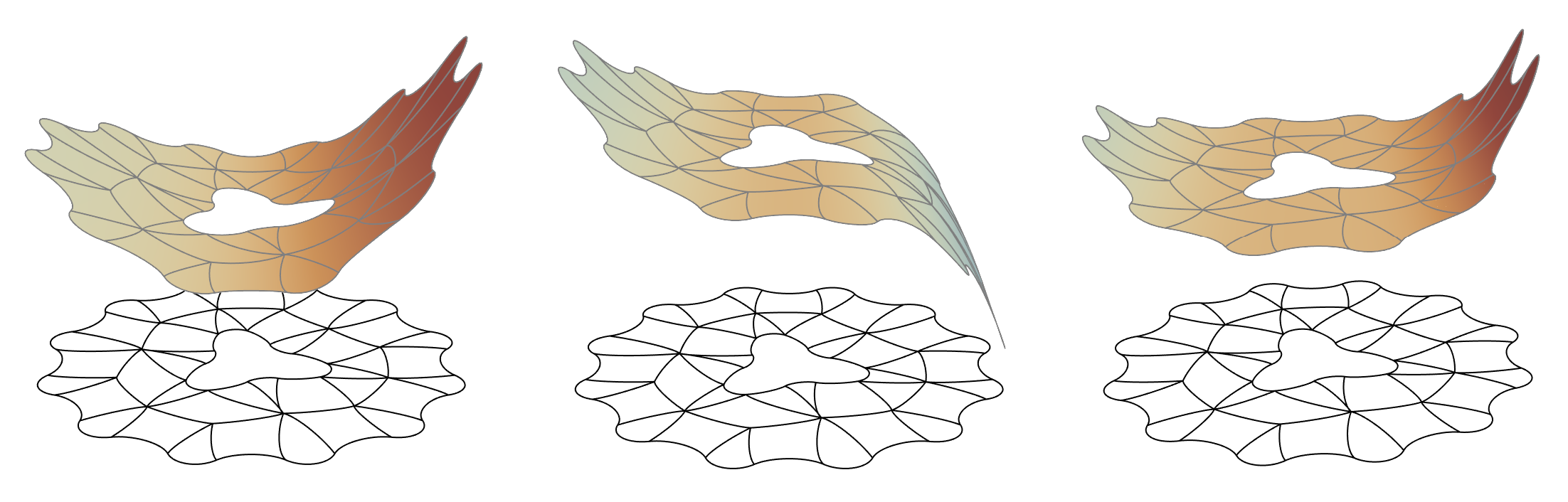}
\caption{Interpolating splines of degree $9$ over Mesh~5 for the polynomial test functions $g_4$ (left), $g_5$ (middle), and $g_6$ (right). }
\label{fig:exm-2-int}
\end{figure}
\end{example}

\subsection{Least-squares approximation}
\label{Examples-ls}

Having an explicit basis for the constructed spline space, it is straightforward to formulate and solve least-squares approximation problems. Let
$\{\bfm{x}_\ell\}_{\ell=1}^{M}\subset\Omega$ be a set of sampling points and
let the corresponding data be given by the values of a sufficiently smooth
trivariate function $g$. We seek a spline $\varphi$ in the considered spline
space that minimizes the discrete least-squares functional
$$
\sum_{\ell=1}^{M}
\bigl(\varphi(\bfm{x}_\ell)-g(\bfm{x}_\ell)\bigr)^2 = \sum_{\ell=1}^{M} \bigl(f^{(k_\ell)}(\bfm{\xi}_\ell)
      -g(\bfm{F}^{(k_\ell)}(\bfm{\xi}_\ell))\bigr)^2,
$$
where $\bfm{x}_\ell=\bfm{F}^{(k_\ell)}(\bfm{\xi}_\ell)$ and $f^{(k_\ell)}$ denotes the local representation of $\varphi$ on the surface element containing
$\bfm{x}_\ell$. 

To generate the sampling points, we fix an integer $K$ and choose uniformly distributed parameter values of degree $K$ in each parametric domain $\Omega^{(k)}$. These parameter values are then mapped to the surface by the corresponding geometry
mappings $\bfm{F}^{(k)}$. Equivalently, the resulting sampling points coincide with the surface domain points associated with a $C^0$-continuous isogeometric discretization of degree $K$.

Expressing the unknown spline in the constructed basis leads to an
overdetermined linear system for its coefficients, which is solved in the
least-squares sense. The system matrix is assembled from the values of the
basis functions at the sampling points. 

\begin{figure}[h!]
\begin{center}
\begin{tikzpicture}
\begin{axis}[
    width=.32\textwidth, height=5cm,
    xlabel={degree $d$},
    ylabel={$\log_{10}(\text{error})$},
    grid=major,
    title={Mesh 1},
    ymax=-2,
    xtick={6, 7, 8, 9, 10},
    ytick={-9, -8, -7, -6, -5, -4, -3}
]
\addplot coordinates {
    (6, -5.304)
    (7, -5.374)
    (8, -7.043)
    (9, -7.094)
    (10, -8.721)
};
\addlegendentry{$\log_{10}(\text{err}_{\rm LS})$}

\addplot coordinates {
    (6, -4.666)
    (7, -4.707)
    (8, -6.135)
    (9, -6.134)
    (10, -7.414)
};
\addlegendentry{$\log_{10}(\text{err}_{\infty,50})$}
\end{axis}
\end{tikzpicture}
\begin{tikzpicture}
\begin{axis}[
    width=.32\textwidth, height=5cm,
    xlabel={degree $d$},
    ylabel={$\log_{10}(\text{error})$},
    grid=major,
    title={Mesh 2},
    xtick={6, 7, 8, 9, 10},
    ytick={-7, -6, -5, -4, -3, -2}
]
\addplot coordinates {
    (6, -3.378)
    (7, -4.138)
    (8, -4.998)
    (9, -5.879)
    (10, -6.490)
};
\addplot coordinates {
    (6, -2.678)
    (7, -3.194)
    (8, -4.036)
    (9, -4.808)
    (10, -5.307)
};
\end{axis}
\end{tikzpicture}
\begin{tikzpicture}
\begin{axis}[
    width=.32\textwidth, height=5cm,
    xlabel={degree $d$},
    ylabel={$\log_{10}(\text{error})$},
    grid=major,
    title={Mesh 3},
    xtick={7, 8, 9, 10},
    ytick={-7, -6, -5, -4, -3}
]
\addplot coordinates {
    (7, -4.542)
    (8, -5.208)
    (9, -6.141)
    (10, -6.784)
};
\addplot coordinates {
    (7, -3.601)
    (8, -3.947)
    (9, -5.041)
    (10, -5.324)
};
\end{axis}
\end{tikzpicture}
\\[10pt]
\begin{tikzpicture}
\begin{axis}[
    width=.32\textwidth, height=5cm,
    xlabel={$\sqrt{\text{NDOF}}$},
    ylabel={$\log_{10}(\text{error})$},
    grid=major,
    title={Mesh 1},
    ymax=-2,
    xtick={15,20,25,30},
    ytick={-9, -8, -7, -6, -5, -4, -3}
]
\addplot coordinates {
    (15.297, -5.304)
    (18.974, -5.374)
    (22.583, -7.043)
    (26.153, -7.094)
    (29.698, -8.721)
};
\addlegendentry{$\log_{10}(\text{err}_{\rm LS})$}

\addplot coordinates {
    (15.297, -4.666)
    (18.974, -4.707)
    (22.583, -6.135)
    (26.153, -6.134)
    (29.698, -7.414)
};
\addlegendentry{$\log_{10}(\text{err}_{\infty,50})$}
\end{axis}
\end{tikzpicture}
\begin{tikzpicture}
\begin{axis}[
    width=.32\textwidth, height=5cm,
    xlabel={$\sqrt{\text{NDOF}}$},
    ylabel={$\log_{10}(\text{error})$},
    grid=major,
    title={Mesh 2},
    xtick={15,20,25,30},
    ytick={-7, -6, -5, -4, -3, -2}
]
\addplot coordinates {
    (14.662, -3.378)
    (18.439, -4.138)
    (22.045, -4.998)
    (25.612, -5.879)
    (29.110, -6.490)
};
\addplot coordinates {
    (14.662, -2.678)
    (18.439, -3.194)
    (22.045, -4.036)
    (25.612, -4.808)
    (29.110, -5.307)
};
\end{axis}
\end{tikzpicture}
\begin{tikzpicture}
\begin{axis}[
    width=.32\textwidth, height=5cm,
    xlabel={$\sqrt{\text{NDOF}}$},
    ylabel={$\log_{10}(\text{error})$},
    grid=major,
    title={Mesh 3},
    xtick={30,35,40,45,50},
    ytick={-7, -6, -5, -4, -3}
]
\addplot coordinates {
    (31.953, -4.542)
    (38.748, -5.208)
    (45.293, -6.141)
    (51.688, -6.784)
};
\addplot coordinates {
    (31.953, -3.601)
    (38.748, -3.947)
    (45.293, -5.041)
    (51.688, -5.324)
};
\end{axis}
\end{tikzpicture}
\end{center}
\caption{Discrete least-squares errors ${\rm err}_{\rm LS}$ and discrete maximum errors ${\rm err}_{\infty,50}$ as functions of the spline degree $d$ (top) and of the square root of the number of degrees of freedom (bottom), shown on a logarithmic scale. 
}
\label{fig:err-exm-3}
\end{figure}
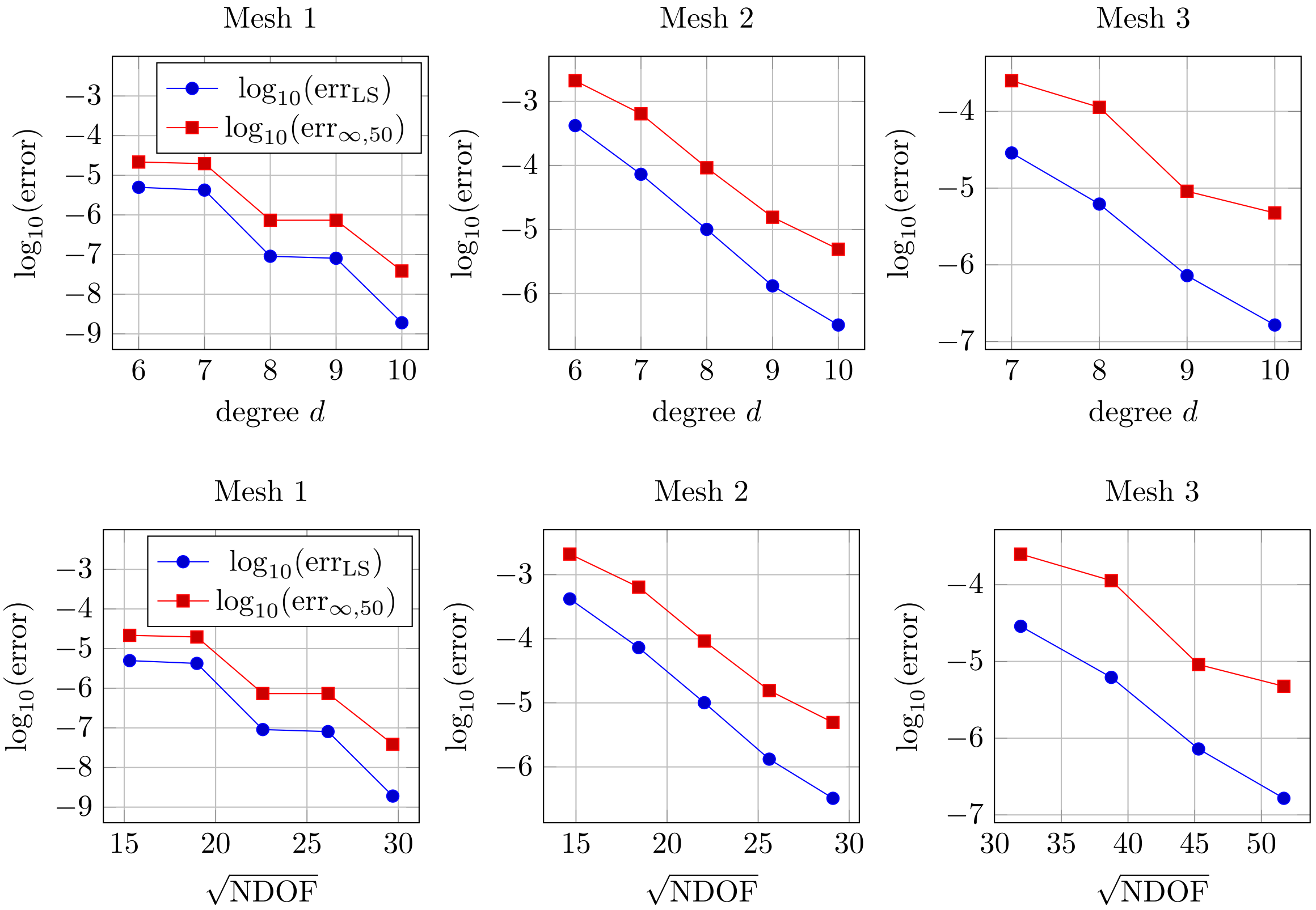

\begin{example} \label{example:LS-1}
We consider least-squares approximation on all three test surfaces corresponding to Meshes~1--3. For each surface, the sampling points are generated from uniformly distributed parameter values of (bi)-degree $K=15$ and mapped to the surface by the corresponding geometry mappings. The resulting sets of sampling points are shown in Figure~\ref{fig:LS-points}. They consist of $M=2702$, $M=2551$ and $M=8191$ points for Mesh~1, Mesh~2, and Mesh~3, respectively. The same sampling points are used for all spline degrees considered on a given mesh.
\begin{figure}[hbt]
\begin{center}
\begin{minipage}{0.32\textwidth}
\centering
\includegraphics[width=1\textwidth]{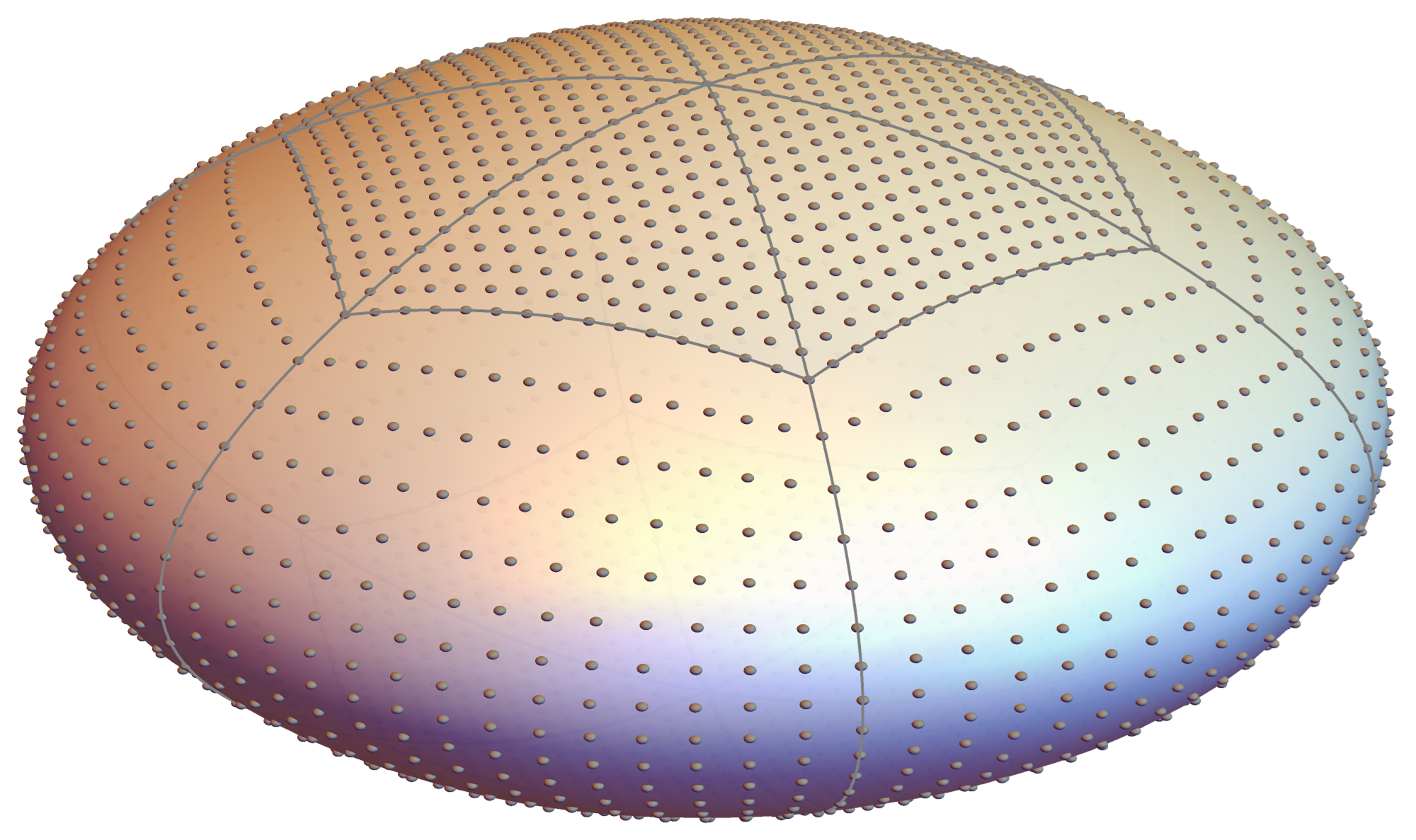}
\end{minipage}
\begin{minipage}{0.32\textwidth}
\centering
\includegraphics[width=1\textwidth]{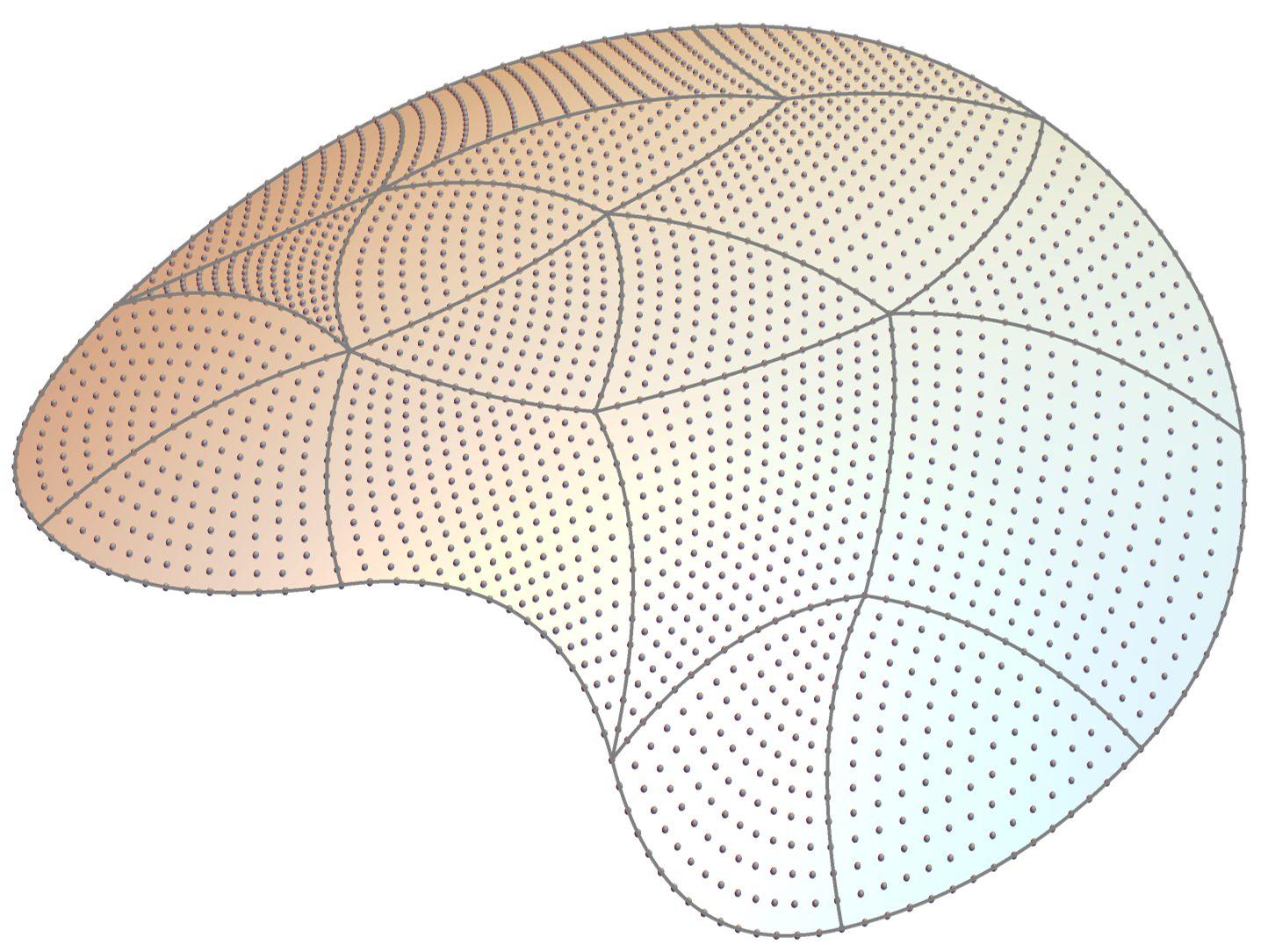}
\end{minipage}
\begin{minipage}{0.32\textwidth}
\centering
\includegraphics[width=1\textwidth]{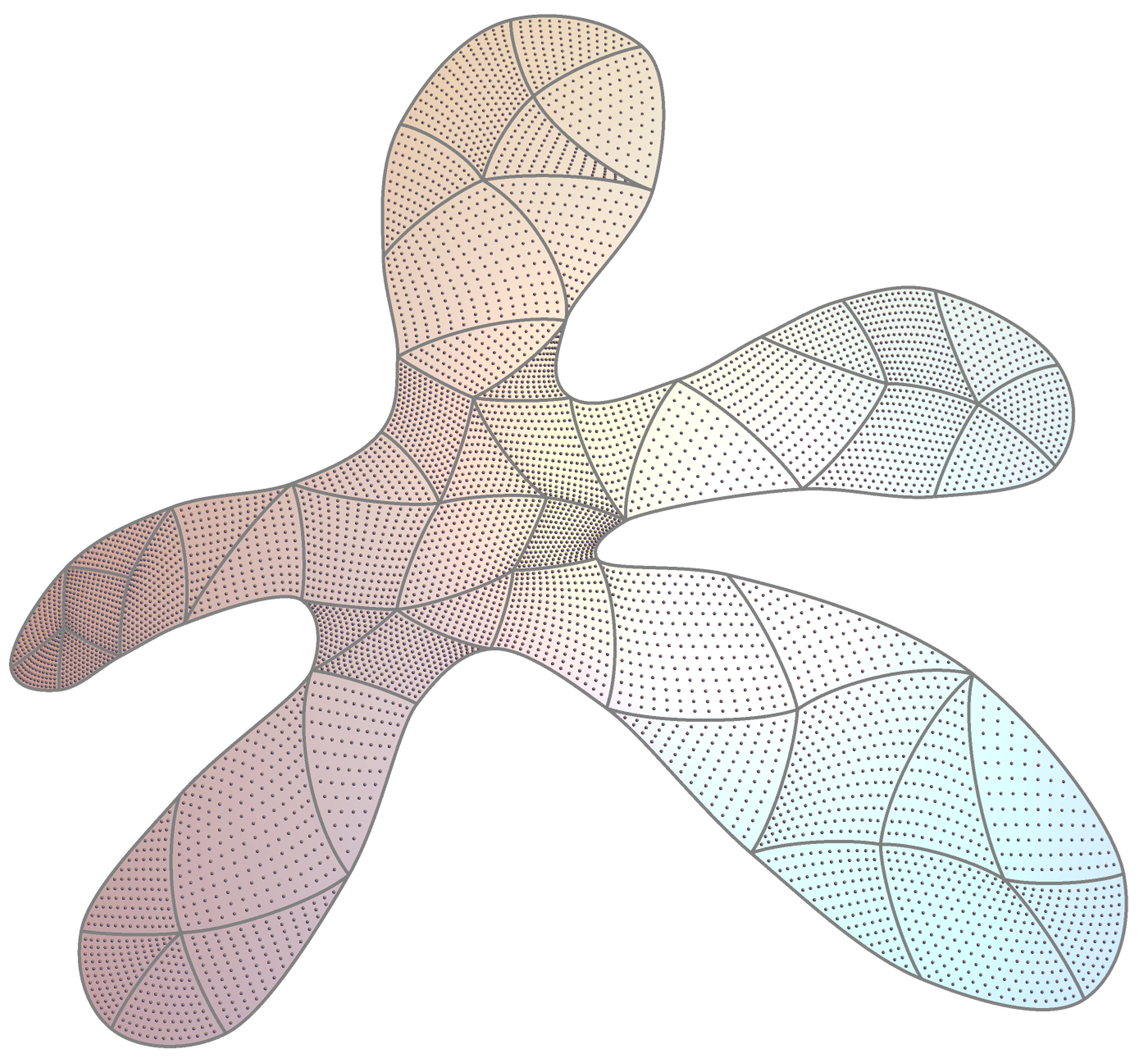}
\end{minipage}
\end{center}
\caption{Sampling points used for the least-squares approximation on Meshes~1--3.}
\label{fig:LS-points}
\end{figure}
For Mesh~1 and Mesh~2, we compute least-squares approximants of degrees
$d=6,7,8,9,10$, while for Mesh~3 we consider degrees $d=7,8,9,10$.
The quality of the approximation is measured by both the discrete
least-squares error
$$
{\rm err}_{\rm LS}(\varphi,g)
:=
\left(\frac{1}{M}
\sum_{\ell=1}^{M}
\bigl(\varphi(\bfm{x}_\ell)-g(\bfm{x}_\ell)\bigr)^2
\right)^{1/2},
$$
and the discrete maximum error \eqref{eq-def-maxError} using the same evaluation sets $\Xi_{50}$ as in the interpolation
experiments.
The obtained errors are reported in Table~\ref{tab:LS-errors}. Figure~\ref{fig:err-exm-3} shows the discrete least-squares error and the discrete maximum error with respect to the spline degree and the square root of the number of degrees of freedom.
For all three meshes, the discrete least-squares error decreases monotonically as the spline degree increases. We numerically observe the exponential convergence as expected for $p$-refinement.
\begin{table}[h!]
\centering
\renewcommand{\arraystretch}{1.1}
\setlength{\tabcolsep}{4pt}
$$
\begin{array}{|c|c|c|c|c|c|c|}
\hline
& \multicolumn{2}{c|}{{\tt Mesh}\,1}
& \multicolumn{2}{c|}{{\tt Mesh}\,2}
& \multicolumn{2}{c|}{{\tt Mesh}\,3} \\
\hline
d
& {\rm err}_{\rm LS} & {\rm err}_{\infty,50}
& {\rm err}_{\rm LS} & {\rm err}_{\infty,50}
& {\rm err}_{\rm LS} & {\rm err}_{\infty,50} \\
\hline
6
& 4.96447\cdot 10^{-6} & 2.16064\cdot 10^{-5}
& 4.19816\cdot 10^{-4} & 2.09968\cdot 10^{-3}
& - & - \\
7
& 4.22751\cdot 10^{-6} & 1.96060\cdot 10^{-5}
& 7.27960\cdot 10^{-5} & 6.39920\cdot 10^{-4}
& 2.87146\cdot 10^{-5} & 2.50996\cdot 10^{-4} \\
8
& 9.06133\cdot 10^{-8} & 7.33523\cdot 10^{-7}
& 1.00479\cdot 10^{-5} & 9.21550\cdot 10^{-5}
& 6.20028\cdot 10^{-6} & 1.12892\cdot 10^{-4} \\
9
& 8.04526\cdot 10^{-8} & 7.34555\cdot 10^{-7}
& 1.31870\cdot 10^{-6} & 1.55687\cdot 10^{-5}
& 7.22258\cdot 10^{-7} & 9.09598\cdot 10^{-6} \\
10
& 1.90015\cdot 10^{-9} & 3.86193\cdot 10^{-8}
& 3.23998\cdot 10^{-7} & 4.92505\cdot 10^{-6}
& 1.64289\cdot 10^{-7} & 4.74980\cdot 10^{-6} \\
\hline
\end{array}
$$
\caption{least-squares approximation errors for Meshes~1--3. The table
reports the discrete least-squares error ${\rm err}_{\rm LS}$ and the
discrete maximum error ${\rm err}_{\infty,50}$.}
\label{tab:LS-errors}
\end{table}
\end{example}

\section{Conclusions}
\label{sec:Conclusion}
 
A $\C{1}$ super-smooth isogeometric spline space $\mathcal{V}^{1,2}_\pd(\Omega)$ has been constructed over mixed triangular and quadrilateral surface meshes whose elements are polynomial surface patches of an arbitrary degree $\delta$. The starting point is a characterization of $\C{1}$-smoothness of an isogeometric function across an interface of two neighboring elements (Theorem~\ref{thm-C1ForTwoElements}). Expressed in local coordinates, the condition states that the trace function $\tr_\ell$ and the crossing derivative function $\w_\ell$ along the interface are the same from both sides (Theorem~\ref{thm-degree-of-omega}), and it involves six gluing functions determined by the two geometry mappings. Both $\tr_\ell$ and $\w_\ell$ are univariate polynomials, and the smoothness condition reduces to divisibility and degree conditions on polynomials built from them. The construction applies in the same way to meshes consisting only of triangles or only of quadrilaterals, for which geometry mappings of arbitrary degree have not been treated before in this context.

Based on this characterization, a numerical procedure (Algorithm~\ref{alg:DegreesOfFreedom}) computes, for a single interface, a Hermite-like representation of the admissible pairs $(\tr_\ell,\w_\ell)$ and returns the number $\chi_\ell$ of degrees of freedom associated with the interface. The procedure applies to geometry mappings of any degree. Imposing $\C{2}$-smoothness at the vertices in addition leads to a spline space whose dimension is given by formula \eqref{eq-dimension} in terms of the numbers $\chi_\ell$. The vertex contribution to the dimension is independent of the geometry, since the $\C{2}$ data at a vertex determine the values and first and second derivatives of all adjacent pull-backs (Theorem~\ref{thm-3}). The basis functions are of four types with local supports, namely vertex, edge, triangle and quadrilateral functions. An equivalent description of the space by interpolation functionals has been given in Section~\ref{subsec:interpolation}.

The numerical experiments in Section~\ref{sec:Examples} confirm the construction on surface meshes associated with quadratic and cubic geometry mappings and on planar meshes. For quadratic geometry mappings, the numbers of degrees of freedom per interface computed by Algorithm~\ref{alg:DegreesOfFreedom} are consistent with the lower bound $2\pd-1$ proven for planar meshes in \cite{GrKaKnTaVi2024}. The basis functions have local supports, and the interpolation problem of Section~\ref{subsec:interpolation} is uniquely solvable for the chosen nodes. Pull-backs of polynomials of degree up to $\lfloor \pd/\delta \rfloor$ belong to the space, and the interpolation scheme reproduces them exactly on the planar test meshes. Interpolation and least-squares approximation of a smooth function on curved surface domains show that the error decreases with the spline degree.

Several questions remain open. The numbers $\chi_\ell$ and the lowest admissible degree are obtained by running Algorithm~\ref{alg:DegreesOfFreedom} for the given geometry mappings. Sharp lower bounds on $\chi_\ell$ in terms of $\pd$ and $\delta$, extending the bound $2\pd-1$ known for planar quadratic mappings, and a characterization of the lowest admissible degree would make the construction fully explicit. The $\C{2}$ conditions at the vertices keep the vertex contribution to the dimension independent of the geometry. A macro-element approach that avoids the additional vertex smoothness could allow lower degrees. The present paper studies the spline space on a fixed mesh for increasing degree, {where we numerically verify the exponential error convergence with respect to the number of degrees of freedom.} Approximation estimates under mesh refinement remain to be established. Finally, the space is intended for the Galerkin discretization of fourth-order problems on surfaces, such as the biharmonic equation and Kirchhoff--Love shells. Its use in this setting raises further questions, among them the choice of quadrature and the conditioning of the resulting systems.

\section*{Acknowledgments}

This paper was developed within the Scientific and Technological Cooperation ``Smooth splines over mixed triangular and quadrilateral meshes for numerical simulation'' between Austria and Slovenia 2023-24, funded by the OeAD under grant nr. SI 17/2023 and by ARIS bilateral project nr. BI-AT/23-24-018. 

The research of M.~Knez is partially supported by the research program P1-0288 and the research grant N1-0237 from ARIS, Republic of Slovenia. The research of J.~Gro\v{s}elj is partially supported by the research program P1-0294 from ARIS, Republic of Slovenia. The research of V.~Vitrih is partially supported by the research program P1-0404 and research project N1-0296 from ARIS, Republic of Slovenia. The research by T.~Takacs was funded in whole or in part by the Austrian Science Fund (FWF) 10.55776/P37177. This support is gratefully acknowledged.

\end{document}